\pdfoutput=1

\RequirePackage{tikz}

\documentclass[sn-basic,iicol,referee]{sn-jnl}%

\usepackage[colorinlistoftodos,prependcaption,textsize=tiny]{todonotes}
\usepackage{newtxmath}
\usepackage{amsmath}
\usepackage{amsmath,bm}
\usepackage{amssymb}
\usepackage{placeins}
\usepackage{ulem}
\usepackage{float}
\usepackage{amsfonts}
\usepackage{bm}
\usepackage{statmath}
\usepackage{graphicx}
\usepackage{cleveref}
\usepackage{xcolor}

\jyear{2022}%

\usepackage[percent]{overpic}    %
\usepackage{units}
\usepackage{ifthen}               %
\usepackage{tikz}
\usepackage{multirow}

\newcommand{\takeout}[1]{ } 

\usepackage{transparent}
\usepackage{eso-pic}

\makeatletter
\renewcommand*\env@matrix[1][\arraystretch]{%
  \edef\arraystretch{#1}%
  \hskip -\arraycolsep
  \let\@ifnextchar\new@ifnextchar
  \array{*\c@MaxMatrixCols c}}
\makeatother
\newcommand{\gfxpath}{images/}   %
\DeclareRobustCommand\DefDof{\tikz[baseline=-0.6ex] \fill[black] (0.ex,0.ex) circle (0.4ex); $\,\! $}
\DeclareRobustCommand\FluxVecDof{\tikz[baseline=-0.6ex] \filldraw[black,fill=white] (0.ex,0.ex) circle (0.8ex); $\,\! $}
\DeclareRobustCommand\NormalTrDof{\tikz[baseline=-0.6ex] \draw[thick,black] (0.ex,0.ex) -- (1.4ex,0.ex);}
\newboolean{boldeqnitalic}
\setboolean{boldeqnitalic}{false}

\usepackage{defdb1}                 %

\definecolor{dgreen}{rgb}{0,0.39,0}
\definecolor{dred}{rgb}{0.8,0,0}
\definecolor{dyellow}{rgb}{.8,.5,0.0}
\definecolor{dblue}{rgb}{0.05, 0.25, 0.75}
\definecolor{d2}{rgb}{0,.5,.5}
\definecolor{mypurple}{rgb}{0.4940,0.1840,0.5560}

\def\FR{\color{black}}    %
\def\OR{\color{black}}      %
\def\SP{\color{black}}     %

\theoremstyle{thmstyleone}%

\theoremstyle{thmstyletwo}%

\theoremstyle{thmstylethree}%

\begin{document}

\title[Spp2256]{Monolithic parallel overlapping Schwarz methods in fully-coupled nonlinear chemo-mechanics problems}

\author*[1]{\fnm{Bjoern} \sur{Kiefer}}\email{bjoern.kiefer@imfd.tu-freiberg.de}
\equalcont{\tiny These authors contributed equally to this work.}

\author[1]{\fnm{Stefan} \sur{Pr\"uger}}\email{stefan.prueger@imfd.tu-freiberg.de}
\equalcont{\tiny These authors contributed equally to this work.}

\author*[2]{\fnm{Oliver} \sur{Rheinbach}}\email{oliver.rheinbach@math.tu-freiberg.de}
\equalcont{\tiny These authors contributed equally to this work.}

\author[2]{\fnm{Friederike} \sur{R\"over}}\email{friederike.roever@math.tu-freiberg.de} 
\equalcont{These authors contributed equally to this work.}

\affil*[1]{\orgdiv{Institute of Mechanics and Fluid Dynamics}, \orgname{TU Bergakademie Freiberg}, \orgaddress{\street{Lampadiusstr.~4}, \postcode{09599} \city{Freiberg}, \country{Germany}}}

\affil[2]{\orgdiv{Institut f\"ur Numerische Mathematik und Optimierung}, \orgname{TU Bergakademie Freiberg}, \orgaddress{\street{Akademiestr.~6}, \postcode{09599} \city{Freiberg}, \country{Germany}}}

\abstract{
We consider the swelling of hydrogels as an example of a chemo-mechanical problem with strong coupling between the mechanical balance relations and the mass diffusion. The problem is cast into a minimization formulation using a time-explicit approach for the dependency of the dissipation potential on the deformation and the swelling volume fraction to obtain symmetric matrices, which are typically better suited for iterative solvers.
The MPI-parallel implementation uses the software libraries deal.II, p4est and FROSch (Fast of Robust Overlapping Schwarz). FROSch is part of the Trilinos library and is used in fully algebraic mode, i.e., the preconditioner is constructed from the monolithic system matrix without making explicit use of the problem structure.
Strong and weak parallel scalability is studied using up to 512 cores, considering the standard GDSW (Generalized Dryja-Smith-Widlund) coarse space and the newer coarse space with reduced dimension.
The FROSch solver is applicable to the coupled problems within in the range of processor cores considered here, although numerical scalablity cannot be expected (and is not observed) for the fully algebraic mode.  In our strong scalability study, the average number of Krylov iterations per Newton iteration is higher by a factor of up to six compared to a linear elasticity problem.  However, making mild use of the problem structure in the preconditioner, this number can be reduced to a factor of two and, importantly, also numerical scalability can then be achieved experimentally. Nevertheless, the fully algebraic mode is still preferable since a faster time to solution is achieved.
}

\takeout{Abstract:Neu:
We consider the swelling of hydrogels as an example of a chemo-mechanical problem  
where a strong coupling between the mechanical balance relations and the mass diffusion is present. 
The problem is cast in a minimization formulation. 
A time-explicit approach is used for the dependency of the dissipation potential on the deformation and the swelling volume fraction to obtain symmetric finite element matrices, which are typically better suited for iterative solvers.
The MPI-parallel implementation makes use of the deal.II and p4est software libraries and uses the FROSch (Fast of Robust Overlapping Schwarz) solver, which is part of the Trilinos software library.
FROSch is used in fully algebraic mode, i.e., the preconditioner is constructed from the assembled, monolithic system matrix without making explicit use 
of the structure of the problem. 

Strong and weak parallel scalability is studied using up to 512 cores of a standard HPC cluster, 
considering the standard GDSW (Generalized Dryja-Smith-Widlund) coarse space for overlapping Schwarz methods and the newer RGDSW (Reduced GDSW) coarse space with reduced dimension. 

The numerical experiments show that the FROSch solver is applicable to the coupled problems in fully algebraic mode within in the range of processor cores considered here although numerical scalablity cannot be expected (and is not observed) for the fully algebraic mode.
In our strong scalability study, the average number of Krylov iterations per Newton iteration is higher by a factor of up to six compared to a linear elasticity problem.

However, making mild use of the problem structure in the construction of the preconditioner, this number can be reduced to a factor of two and, importantly, also numerical scalability can then be  achieved experimentally. 

Nevertheless, the fully algebraic mode is still preferable in our setting since a fast time to solution is achieved. 
}


\newpage
\hspace*{1ex}
\thispagestyle{empty}
{ \setlength{\unitlength}{1cm}
\begin{picture}(10,30)
 \put(0,0){A}
 \put(10,0){B}
 \put(0,0){C}
 \put(10,30){D}
 \put(-11.05,6){\includegraphics[width=\paperwidth,height=\paperheight]{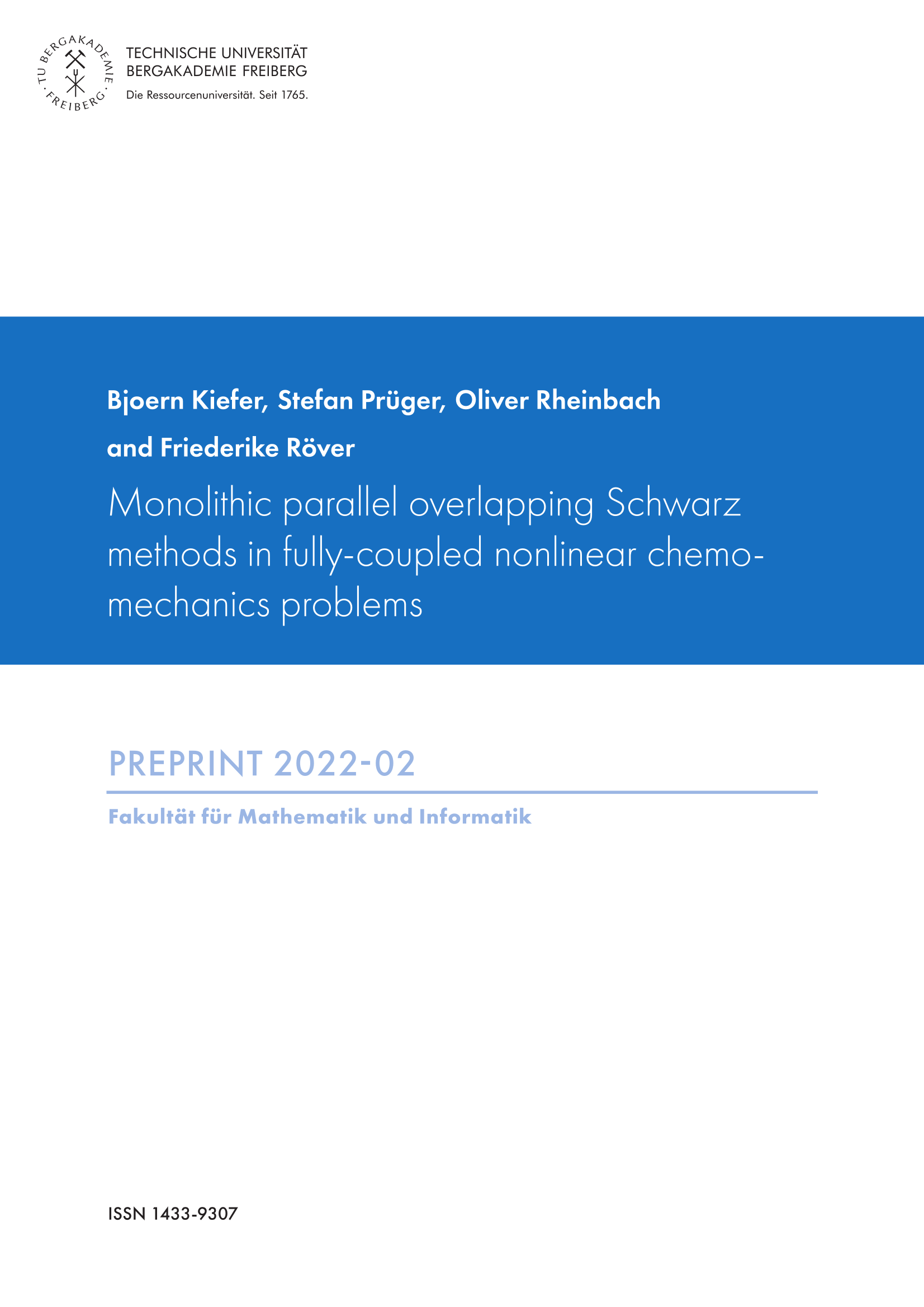}}
\end{picture}
}

\newpage
\hspace*{1ex}
\thispagestyle{empty}
{ \setlength{\unitlength}{1cm}
\begin{picture}(10,30)
 \put(0,0){A}
 \put(10,0){B}
 \put(0,0){C}
 \put(10,30){D}
 \put(-3,6){\includegraphics[width=\paperwidth,height=\paperheight]{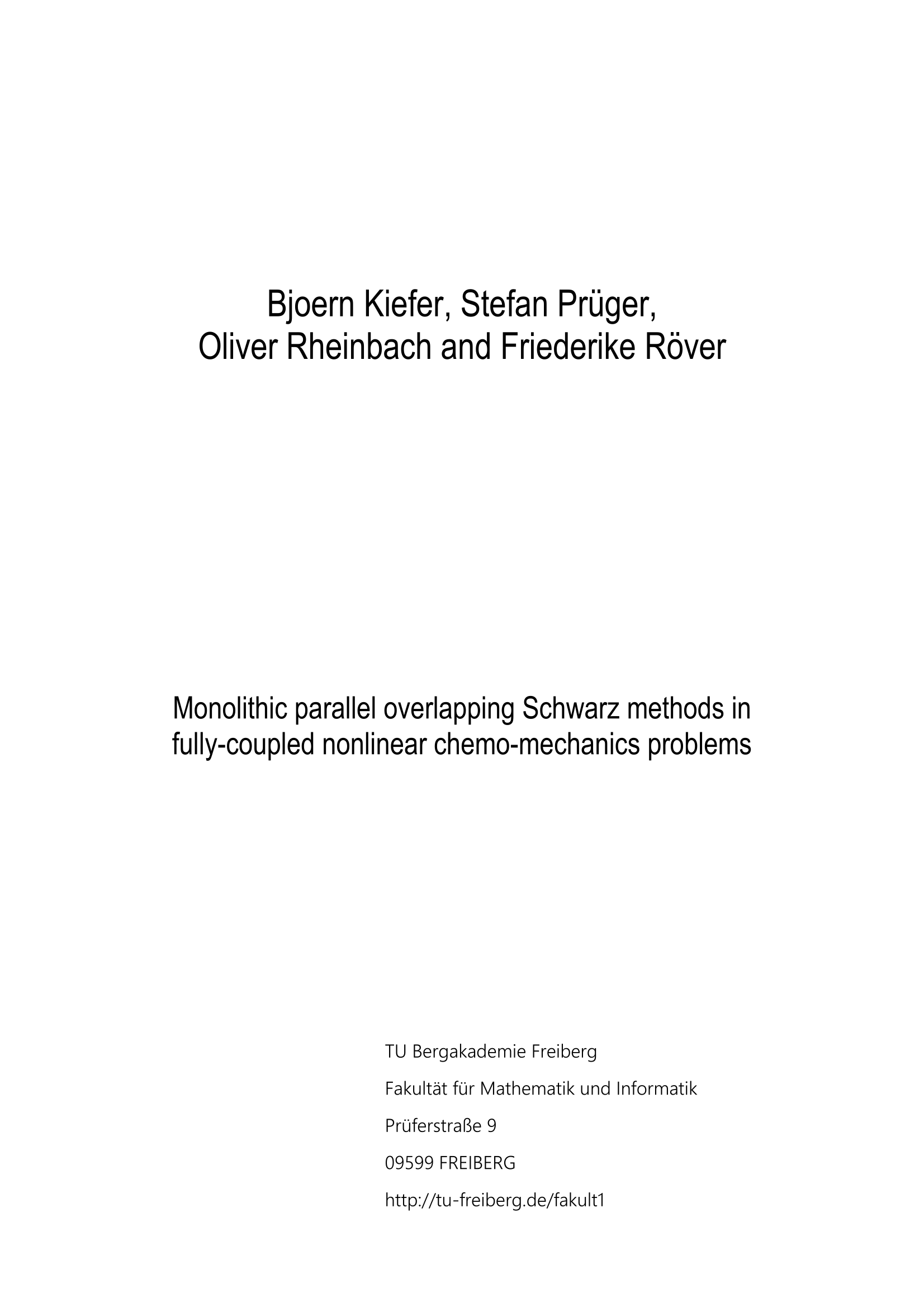}}
\end{picture}
}

\clearpage
\newpage
\hspace*{1ex}
\thispagestyle{empty}
{ \setlength{\unitlength}{1cm}
\begin{picture}(10,30)
 \put(0,0){A}
 \put(10,0){B}
 \put(0,0){C}
 \put(10,30){D}
 \put(-11.05,6) {\includegraphics[width=\paperwidth,height=\paperheight]{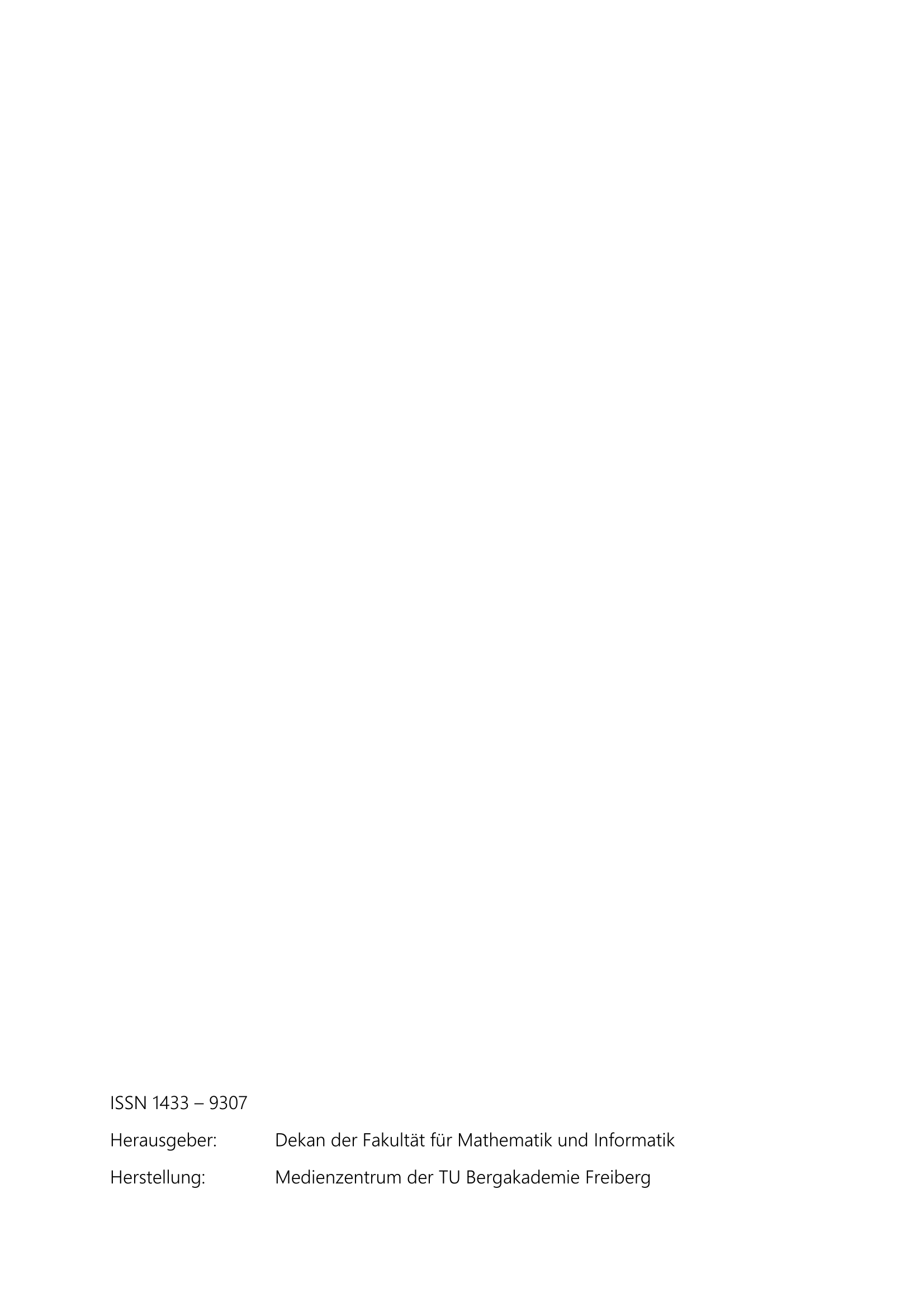}}
\end{picture}
}
\newpage

\clearpage
\newpage
\hspace*{1ex}
\thispagestyle{empty}
\newpage

\setcounter{page}{0}
\maketitle

\section{Introduction}
{\SP Chemo-mechanics problems have gained %
{\OR increasing}
attention in the past decades, as a more refined understanding of processes in man-made and natural materials as well as living tissue can only be obtained by incorporation of mechanical and chemical loading conditions and their mutual interactions.
  Research in various fields of chemo-mechanics have emerged and are concerned for example with the prediction of the transition from pitting corrosion to crack formation in metals~\cite{CuiMaMar:2021:pha,CheJafZhaBob:2021:cou}, hydrogen diffusion~\cite{DiAna:2013:hyd,SalMcMGraMag:2018:cou} and embrittlement~\cite{AnaMaoTal:2019:mod,KriNioMar:2020:pha,AutBroEkh:2022:ful}, functional degradation in Li-ion batteries~\cite{RejDiAna:2015:the,RezAshXu:2021:con}, chemical reaction-induced degradation of concrete~\cite{WuTemWri:2014:mul,NguWalBui:2019:com} and diffusion mediated tumor growth~\cite{XueYinLiFen:2018:bio,FagFenLimOdeYan:2020:cou}.
  As all of these examples involve a strong coupling of mechanic{\OR al} balance relations and mass diffusion, either of Fickian or gradient extended Cahn-Hilliard type, we adopt a simple benchmark problem of swelling of hydrogels~\cite{BouLanHua:2015:non,SprKieMen:2016:com} that at the one hand accounts for this coupling and on the other hand is simple enough to develop efficient, problem specific numerical solution schemes.
  
  {\OR In this paper, we are interested in the model presented in~\cite{BogNatMie:2017:min},}
  which is derived from an incremental variational formulation and can therefore easily be recast into a minimization formulation as well as into a saddle point problem.
  {\OR The different variational formulations also have consequences for the solver algorithms to be applied.}
}

{\OR In this contribution, as a first step, we}
consider the minimization formulation. The discretization of our three-dimensional model problem by finite elements is carried out using the deal.II finite element software library~\cite{dealii2019design}. We solve the arising nonlinear system by means of a monolithic Newton-Raphson scheme; the linearized systems of equations %
are solved using the \textbf{F}ast and \textbf{R}obust \textbf{O}verlapping  \textbf{Sch}warz (FROSch) solver~\cite{Heinlein:2016:PIT,Heinlein:2018:FPI} which is part of the Trilinos software~\cite{trilinosrepo}.
The FROSch framework provides a parallel implementation of the GDSW~\cite{Dohrmann:2008:DDL,Dohrmann:2008:FEM} and RGDSW-type~Overlapping Schwarz preconditioners~\cite{Dohrmann:2017:DSC}. 
These preconditioners have shown a good performance for %
problems ranging from the
fluid-structure interaction of arterial walls and blood flow~\cite{FROSch:FSI:2016} to land ice simulation~\cite{Heinlein:2021:FRO}.
Within this project, it has first been applied in %
\cite{Kiefer:PAM:2021}.
The preconditioner also provides a recent %
extension {\OR using more than two levels}, which has been tested up to 220\,000 cores~\cite{Heinlein:2022:PST}. %

The FROSch preconditioners considered in this paper can be constructed algebraically, i.e., from the assembled finite element matrix, without the use of geometric information. In this paper, we will apply the preconditioners from FROSch in this algebraic mode. 
Note that recent RGDSW methods with adaptive coarse spaces are not fully algebraic; e.g.,~\cite{Heinlein:2022:AGC}.

{\OR For our benchmark problem of the swelling of hydrogels, we consider two sets of boundary conditions.
We also compare the consequences of two different types of finite element discretizations for the flux flow: %
Raviart-Thomas finite elements and standard Lagrangian finite elements.  
We then evaluate the numerical and parallel performance of the iterative solver applied to the monolithic system, discussing strong and weak parallel scalability.}   

\section{Variational framework of fully coupled chemo-mechanics}
\label{sec:variational-framework-fully-coupled-chemo-mech}
In order to evaluate the performance of the FROSch framework in the context of multi-physics problems, the variational framework of chemo-mechanics is adopted as outlined in~\cite{BogNatMie:2017:min}. {\OR This framework} is suitable to model hydrogels.

This setting is employed here to solve some representative model problems involving full coupling between mechanics and mass diffusion in a finite deformation setting.

{\OR The} rate-type potential
\begin{equation}
  \label{eq:gen-time-continuous-potential}
  \Pi\left(\dot{\Bvarphi},\dot{v},\bm{J}_{v}\right)=\dd{}{t} E\left(\dot{\Bvarphi},\dot{v}\right)+D\left(\bm{J}_{v}\right) - P_{\mathrm{ext}}\left(\dot{\Bvarphi},\bm{J}_{v}\right)
\end{equation}
serves as a starting point for our description of the coupled problem, 
{\OR where the deformation is denoted $\Bvarphi$, the swelling volume fraction $v$, and the fluid flux $\bm{J}_{v}$.}
The stored energy functional {\FR$E$} of the body $\Br$ is computed from the free-energy $\psi$ as
\begin{equation}
  \label{eq:stored-energy-functional}
  E\left(\Bvarphi,v\right)=\int\nolimits_{\Br}\widehat{\psi}\left(\nabla\Bvarphi,v\right)\mathrm{d}V \mbox{.}
\end{equation}
Furthermore, the global dissipation potential functional is defined as
\begin{equation}
  \label{eq:dissipation-potential-functional}
  D\left(\bm{J}_{v}\right)=\int\nolimits_{\Br}\widehat{\phi}\left(\bm{J}_{v};\nabla\Bvarphi,v\right)\mathrm{d}V \ \mbox{,}
\end{equation}
involving the local dissipation potential $\widehat{\phi}$.

Note that the dissipation potential possesses an additional dependency on the deformation via its material gradient and the swelling volume fraction.
However, this dependency is not taken into account in the variation of the potential~$\Pi$ when determining the corresponding Euler-Lagrange equations, as indicated by the semicolon in the list of arguments.
Lastly, the external load functional is split into a solely mechanical and solely chemical contribution of the form
\begin{equation}
  \label{eq:external-load-functional}
  P_{\mathrm{ext}}\left(\dot{\Bvarphi},\bm{J}_{v}\right)=P_{\mathrm{ext}}^{\Bvarphi}\left(\dot{\Bvarphi}\right) + P_{\mathrm{ext}}^{\mu}\left(\bm{J}_{v}\right) \ \mbox{,}
\end{equation}
where the former includes the vector of body forces per unit reference volume $\bm{R}_{\Bvarphi}$ and the prescribed traction vector $\bar{\bm{T}}$ at the surface of the body $\partial \Br^{\bm{T}}$ such that
\begin{equation}
  \label{eq:external-load-functional-mechanical}
  P_{\mathrm{ext}}^{\Bvarphi}\left(\dot{\Bvarphi}\right)=\int\nolimits_{\Br} \bm{R}_{\Bvarphi}\cdot\dot{\Bvarphi}\,\mathrm{d}V + \int\nolimits_{\partial \Br^{\bm{T}}}\bar{\bm{T}}\cdot\dot{\Bvarphi}\,\mathrm{d}A \ \mbox{.}
\end{equation}
The latter contribution in~\eqref{eq:external-load-functional} incorporates the prescribed chemical potential $\bar{\mu}$ and the normal component of the fluid flux $H_{v}$ at the surface $\partial \Br^{\mu}$ as
\begin{equation}
  \label{eq:external-load-functional-diffusion}
  P_{\mathrm{ext}}^{\mu}\left(\bm{J}_{v}\right)=- \int\nolimits_{\partial \Br^{\mu}} \bar{\mu}\underbrace{\bm{J}_{v}\cdot\bN}_{H_{v}}\mathrm{d}A \ \mbox{.}
\end{equation}
Along the disjoint counterparts of the mentioned surface, namely $\partial \Br^{\Bvarphi}$ and $\partial \Br^{H_{v}}$, the deformation and the normal component of the fluid flux are prescribed, respectively.\newline
Taking into account the balance of solute volume
\begin{equation}
  \label{eq:field-equation-balance-solute-volume}
  \dot{v}=-\Divb{\bm{J}_{v}}
\end{equation}
in~\eqref{eq:gen-time-continuous-potential} allows one to derive the two-field minimization principle
\begin{multline}
  \label{eq:continuous-two-field-minimization-principle}
  \Pi\left(\dot{\Bvarphi},\bm{J}_{v}\right)=\\
  \int\nolimits_{\Br} \underbrace{\ps{\bF}\widehat{\psi}:\nabla\dot{\Bvarphi}-\ps{v}\widehat{\psi}\Divb{\bm{J}_{v}}+\widehat{\phi}\left(\bm{J}_{v};\nabla\Bvarphi,v\right)}_{\pi\left(\nabla\dot{\Bvarphi},\bm{J}_{v},\Divb{\bm{J}_{v}}\right)}\mathrm{d}V \\
  - P_{\mathrm{ext}}\left(\dot{\Bvarphi},\bm{J}_{v}\right) \ \mbox{,}
\end{multline}
which solely depends on the deformation and the fluid flux.
Herein,~\eqref{eq:field-equation-balance-solute-volume} is accounted for locally to capture the evolution of $v$ and update the corresponding material state.
To summarize, the deformation map and the flux field are determined from
\begin{equation}
  \label{eq:minimization-problem-two-field-formulation}
  \left\{\dot{\Bvarphi},\bm{J}_{v}\right\} = \mathrm{Arg}\left\{\underset{\dot{\Bvarphi}\in \mathcal{W}_{\dot{\Bvarphi}}}{\mathrm{inf}}\ \underset{\bm{J}_{v}\in \mathcal{W}_{\bm{J}_{v}}}{\mathrm{inf}} \Pi\left(\dot{\Bvarphi},\bm{J}_{v}\right) \right\} \, \text{,}
\end{equation}
using the following admissible function spaces.
\begin{align}
  \label{eq:admissible-function-space-deformation}
  \mathcal{W}_{\dot{\Bvarphi}}&=\left\{\dot{\Bvarphi}\in H^{1}\!\left(\Br\right) \vert \ \dot{\Bvarphi}=\dot{\bar{\Bvarphi}} \ \text{on} \ \partial \Br^{\Bvarphi}\right\}\\
  \label{eq:admissible-function-space-fluid-flux}
  \mathcal{W}_{\bm{J}_{v}}&=\left\{\bm{J}_{v}\in H\left(\Div,\Br\right) \vert \ \bm{J}_{v}\!\cdot\bN=H_{v} \ \text{on} \ \partial \Br^{H_{v}}\right\}
\end{align}
\subsection{Specific free-energy function and dissipation potential}
\label{sec:specific-free-energy-funct-dissipation-potential}
The choice of the free-energy function employed in this study is motivated by the fact that it accurately captures the characteristic nonlinear elastic response of a certain class of hydrogels~\cite{ManAroTam:2021:inf} with moderate water content as well as the swelling induced volume changes.
In principle, it also incorporates softening due pre-swelling, despite its rather simple functional form.
The isotropic, Neo-Hookean type free-energy reads as
\begin{multline}
  \label{eq:free-energy-function-isotropic-pre-swollen-state}
  \widehat{\psi}\left(\bF, v\right)=\frac{\upgamma}{2 J_{0}}\left[J_{0}^{2/3}I_{1}^{\bC} -3 -2 \ln\left(JJ_{0}\right)\right] \\ + \frac{\uplambda}{2J_{0}}\left[JJ_{0}-1- v\right]^{2}\\
   +\frac{\upalpha}{J_{0}}\left[ v \ln\left(\frac{ v}{1+ v}\right)+\frac{\upchi v}{1+ v}\right] \ \mbox{,}
\end{multline}
in which the first invariant of the right Cauchy-Green tensor is defined as $I_{1}^{\bC}=\tr\left(\bC\right)=\tr\left(\bF^{\mathrm{T}}\!\cdot\bF\right)$, while the determinant of the deformation gradient $\bF=\nabla\Bvarphi$ is given as $J=\det\left(\bF\right)$.
The underlying assumption of the particular form of this {\FR energy} function is the multiplicative decomposition of the deformation gradient
\begin{equation}
  \label{eq:mulitplicative-split-pre-swelling}
  \bF^{\mathrm{d}}=\bF\cdot\bF_{\mathrm{0}}=J_{0}^{1/3}\bF
\end{equation}
which splits the map from the reference configuration (dry-hydrogel) to the current configuration into a purely volumetric deformation gradient associated with the pre-swelling of the hydrogel $\bF_{\mathrm{0}}$ and the deformation gradient $\bF$, accounting for elastic and diffusion-induced deformations.
Clearly,~\eqref{eq:free-energy-function-isotropic-pre-swollen-state} describes the energy relative to the pre-swollen state of the gel.
In its derivation it is assumed, additionally, that the pre-swelling is stress-free and the energetic state of the dry-state and the pre-swollen state are equivalent, which gives rise to the scaling $J_{\mathrm{0}}^{-1}$ of the individual terms of the energy.
Although the incompressibility of both the polymer, forming the dry hydrogel, and the fluid is widely accepted, its exact enforcement is beyond the scope of the current study and a penalty formulation is employed here, utilizing a quadratic function that approximately enforces {\SP the coupling constraint}
{\SP
  \begin{equation}
    \label{eq:constraint-vol-change-swelling}
    JJ_{0}-1- v=0
  \end{equation}%
}%
for a sufficiently high value of $\uplambda$.
{\SP Thus, in the limit $\uplambda\rightarrow \infty$, the volume change in the hydrogel is solely due to diffusion and determined by the volume fraction $v$, which characterizes the amount of fluid present in the gel.
On the other hand, relaxing the constraint by choosing a small value of $\uplambda$ allows for additional elastic volume changes.}\newline
Energetic contributions due the change in fluid concentration in the hydrogel are accounted for by the Flory-Rehner type energy, in which the affinity between the fluid and the polymer network is controlled by the parameter $\upchi$.
Finally, demanding that the pre-swollen state is stress-free requires the determination of the initial swelling volume fraction from
\begin{equation}
  \label{eq:initialization-swelling-vol-frac}
   v_{0}= \frac{\upgamma}{\uplambda}\left[J_{0}^{-1/3}-\frac{1}{J_{0}}\right]+J_{0}-1\ \mbox{.}
\end{equation}
A convenient choice of the local dissipation potential, in line with~\cite{BogNatMie:2017:min}, is given as
\begin{equation}
  \label{eq:dissipation-potential-minimization}
  \widehat{\phi}\left(\bm{J}_{ v};\bC, v\right) = \frac{1}{2 \mathrm{M} v}\bC:\left(\bm{J}_{ v}\otimes\bm{J}_{ v}\right) \ \mbox{,}
\end{equation}
which is formulated with respect to the pre-swollen configuration.
It is equivalent to an isotropic, linear constitutive relation between the spatial gradient of the chemical potential and the spatial fluid flux, 
in the current configuration.
Again, the state dependence of the dissipation potential through the right Cauchy-Green tensor and the swelling volume fraction is not taken into account in the course of the variation of the total potential.
The material parameters employed in the strong and weak scaling studies in Section~\ref{sec:strong-scaling} and Section~\ref{sec:weak-scaling} are summarized in Table~\ref{tab:material-parameter}.
Note that the chosen value of the pre-swollen Jacobian is problem dependent and indicated in the corresponding section.
\begin{table}[h]
 \centering
 \begin{tabular}[h]{l l l l}
   symbol & physical meaning & value & unit \\ \hline
   $\upgamma$ & shear modulus & $0.1$ & $\unitfrac[]{N}{mm^{2}}$ \\
   $\upalpha$ & mixing modulus & $24.2$ &$\unitfrac[]{N}{mm^2}$ \\
   \multirow{2}{*}{$\upchi$} & mixing control& \multirow{2}{*}{$0.2$} & \multirow{2}{*}{-} \\
   &parameter&&\\
   \multirow{2}{*}{$\mathrm{M}$} & volumetric diffusity & \multirow{2}{*}{$10^{-2}$} & \multirow{2}{*}{$\unitfrac[]{mm^4}{Ns}$} \\
   &parameter&&\\
   $J_0$ & pre-swollen Jacobian  & $1.01 - 4.5$ & -\\
   \multirow{2}{*}{$\uplambda$}& volumetric penalty & \multirow{2}{*}{$0.2$} & \multirow{2}{*}{$\unitfrac[]{N}{mm^2}$} \\
   &parameter&&\\ \hline
 \end{tabular}
 \caption{%
 Material parameters of the coupled hyperelastic model.}
 \label{tab:material-parameter}
\end{table}
\subsection{Incremental two-field potential}
\label{sec:incr-two-field-potential}
Although the rate-type potential~\eqref{eq:gen-time-continuous-potential} allows for valuable insight into the
variational structure of the coupled problem, e.g., a minimization formulation for the model at hand, an incremental framework is typically required for the implementation into a finite element code.
This necessitates the integration of the total potential over a finite time step $\Delta t = t_{n+1}-t_{n}$.
Thus, the incremental potential takes the form
\begin{multline}
  \label{eq:time-discrete-two-field-principle}
  \Pi^{\Delta t}\left(\Bvarphi,\bm{J}_{v}\right)=\\
  \int\nolimits_{\Br} \underbrace{\widehat{\psi}\left(\nabla \Bvarphi,v_{n}-\Delta t \Divb{\bm{J}_{v}}\right) + \Delta t \widehat{\phi}\left(\bm{J}_{v};\nabla\Bvarphi_{n},v_{n}\right)}_{\pi^{\Delta t}\left(\nabla\Bvarphi,\bm{J}_{v},\Divb{\bm{J}_{v}}\right)}\mathrm{d}V \\
  - \int \nolimits_{\Br} \bm{R}_{\Bvarphi}\cdot\left(\Bvarphi-\Bvarphi_{n}\right) \mathrm{d}V - \int \nolimits_{\partial \Br^{\bm{T}}} \bar{\bm{T}}\cdot\left(\Bvarphi-\Bvarphi_{n}\right)\mathrm{d}A \\
  + \int \nolimits_{\partial\Br^{\mu}} \Delta t\, \bar{\mu}\, \bm{J}_{v}\cdot\bm{N} \,\mathrm{d}A \ \mbox{,}
\end{multline}
in which an Euler implicit time integration scheme is applied to approximate the global dissipation potential~\eqref{eq:dissipation-potential-functional} as well as the external load functional~\eqref{eq:external-load-functional}.
Furthermore, the balance of solute volume is also integrated numerically by means of the implicit backward Euler scheme yielding an update formula for the swelling volume fraction
\begin{equation}
  \label{eq:time-discrete-balance-solute-volume}
  v=v_{n}-\Delta t\Divb{\bm{J}_{v}} \ \mbox{,}
\end{equation}
which is employed to evaluate the stored energy functional~\eqref{eq:stored-energy-functional} at $t_{n+1}$.
Note that, quantities given at the time step $t_{n}$ are indicated by subscript ${}_{n}$, while the subscript is dropped for all quantities at $t_{n+1}$ to improve readability.
Additionally, it is remarked that the stored energy functional at $t_{n}$ is excluded from~\eqref{eq:time-discrete-two-field-principle} as it only changes its absolute value, while it does not appear in the first variation of $\Pi^{\Delta t}$, because it is exclusively dependent on quantities at $t_{n}$.
Finally, the state dependence of the local dissipation potential is only accounted for in an explicit manner, in order to ensure consistency with the rate-type potential~\eqref{eq:gen-time-continuous-potential} and thus guarantee the symmetry of the tangent operators in the finite element implementation.

{\OR For symmetric systems, we can hope for a better convergence of the Krylov methods applied to the preconditioned system. On the other hand, we are restricted by small time steps.}
\section{Fast and Robust Overlapping Schwarz (FROSch) Preconditioner}
\label{sec:frosch}

Domain decomposition solvers~\cite{Toselli:2005:DDM} %
are based on the idea to construct an approximate solution to a problem, defined on a computational domain,
from the solutions of parallel problems on small subdomains and, typically, of an additional coarse problem, which introduces the global coupling.
Traditionally, the coarse problem is defined on a coarse mesh. In the FROSch software, however, we only consider methods where no explicit coarse mesh needs to be provided.

Domain decomposition methods are typically used as a preconditioner in combination with Krylov subspace methods such as conjugate gradients or GMRES.

In overlapping Schwarz domain decomposition methods~\cite{Toselli:2005:DDM}, the subdomains have some overlap. A large overlap increases the size of the subdomains but typically improves the speed of convergence.

The C++ library FROSch~\cite{Heinlein:2016:PIT,Heinlein:2018:FPI},
which is part of the Trilinos software library~\cite{trilinosrepo},
implements versions of the \textit{Generalized Dryja-Smith-Widlund (GDSW)} preconditioner, which is a two-level overlapping Schwarz domain decomposition preconditioner~\cite{Toselli:2005:DDM} using an energy-minimizing coarse space introduced in~\cite{Dohrmann:2008:DDL,Dohrmann:2008:FEM}.
This coarse space is inspired by iterative substructuring methods such as FETI-DP and BDDC methods~\cite{Toselli:2005:DDM}.
An advantage of GDSW-type preconditioners, compared to iterative substructuring methods and classical two-level Schwarz domain decomposition preconditioners, is that they can be constructed in an algebraic fashion from the fully assembled stiffness matrix.
\begin{figure}
\centering
\includegraphics[width=0.4\textwidth, trim={40mm 0mm 40mm 15mm},clip]{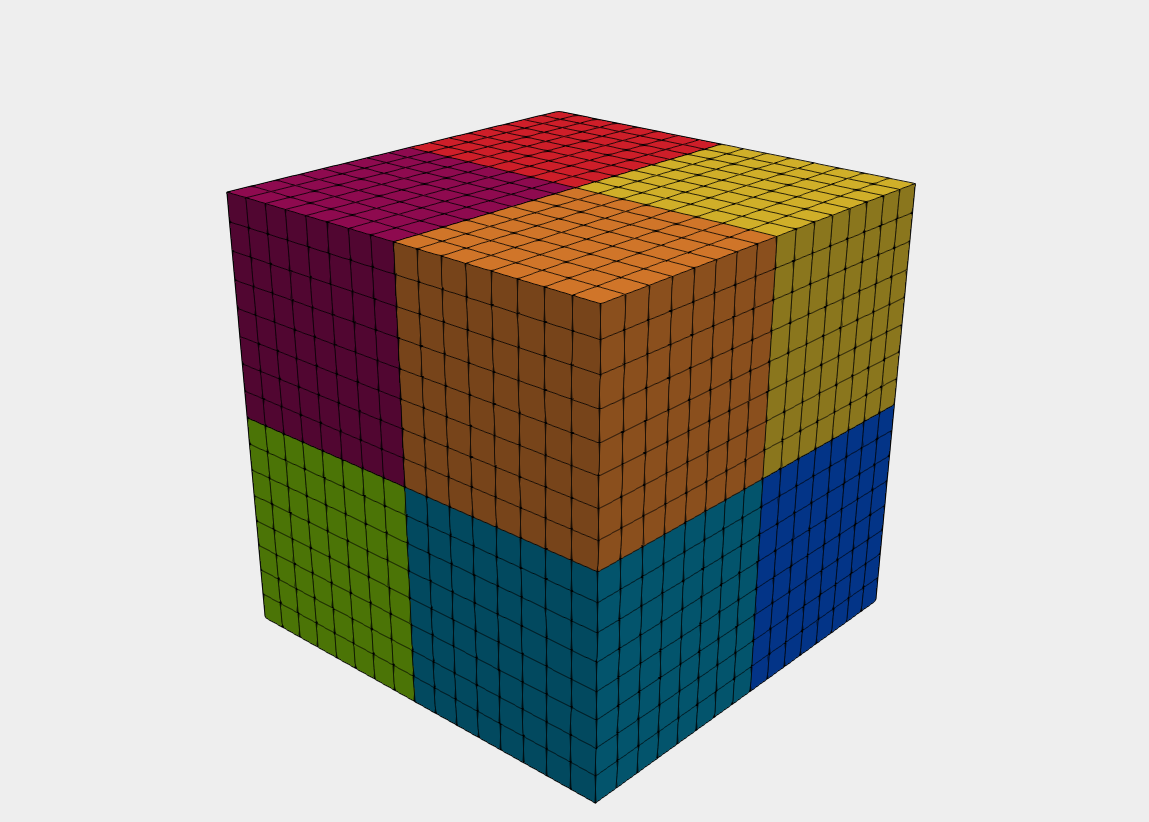} \\
\includegraphics[width=0.4\textwidth, trim={40mm 0mm 40mm 15mm},clip]{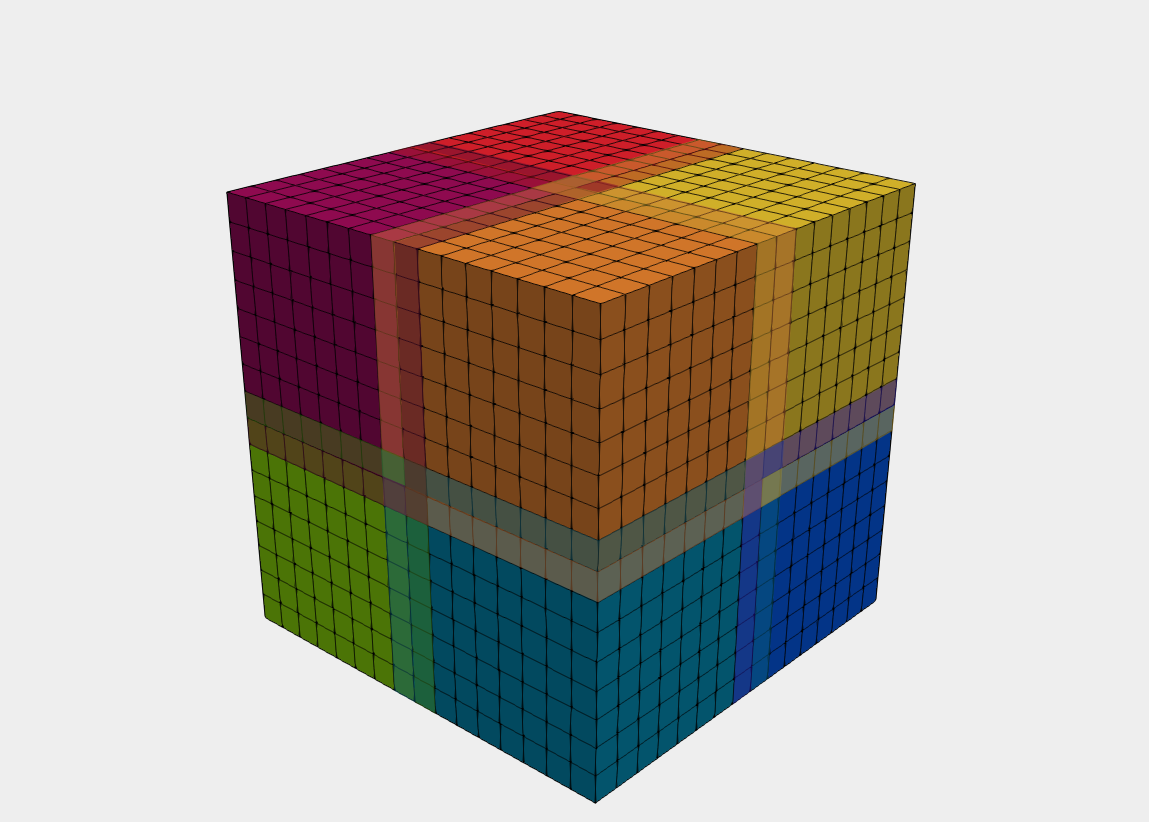}
\caption{Decomposition of a cube into non-overlapping (top) and overlapping subdomains (bottom) on a structured grid with $\delta=1h$.
{\OR The GDSW preconditioner uses the overlapping subdomains to define the local solvers and the non-overlapping subdomains to construct the second level, which  ensures global transport of information.}
}
\label{fig:decomposition}
\end{figure}
Therefore, they do not require a coarse triangulation (as in classical two-level Schwarz methods) nor access to local Neumann matrices for the subproblems (as in FETI-DP and BDDC domain decomposition methods). 

{\OR For simplicity, we will describe the construction of the preconditioner in terms of the computational domain,
although the construction is fully algebraic in FROSch, i.e., subdomains arise only implicitly from the algebraic construction:} %
the computational domain $\Omega$ is decomposed into non-overlapping subdomains $\lbrace \Omega_i\rbrace_{i = 1 \ldots N}$; see Figure~\ref{fig:decomposition}. Extending each subdomain by $k$-layers of elements we obtain the overlapping subdomains $\lbrace \Omega_i' \rbrace_{i = 1 \ldots N}$ with an overlap $\delta = k h$, where $h$ is the size of the finite elements. We denote the size of a non-overlapping subdomain by $H$. The GDSW preconditioner can be written in the %
form
\begin{equation}
M_{\rm  GDSW}^{-1} ={\Phi K_0^{-1} \Phi^T}+\sum\nolimits_{i = 1}^N R_i^T K_i^{-1} R_i\,,
\label{eq:gdsw}
\end{equation}
where $K_i = R_i K R_i^T, i = 1, \ldots N$ {\FR represent the local overlapping subdomain problems.}
The coarse problem is given by the Galerkin product $K_0 = \Phi^T K \Phi$. {\FR The matrix $\Phi$ contains the coarse basis functions spanning the coarse space $V^0$.} For the classical two-level overlapping Schwarz method these functions would be nodal finite elements functions on a coarse triangulation. 

\begin{figure}
\centering
\includegraphics[width=0.3\textwidth, trim={20mm 80mm 20mm 15mm},clip]{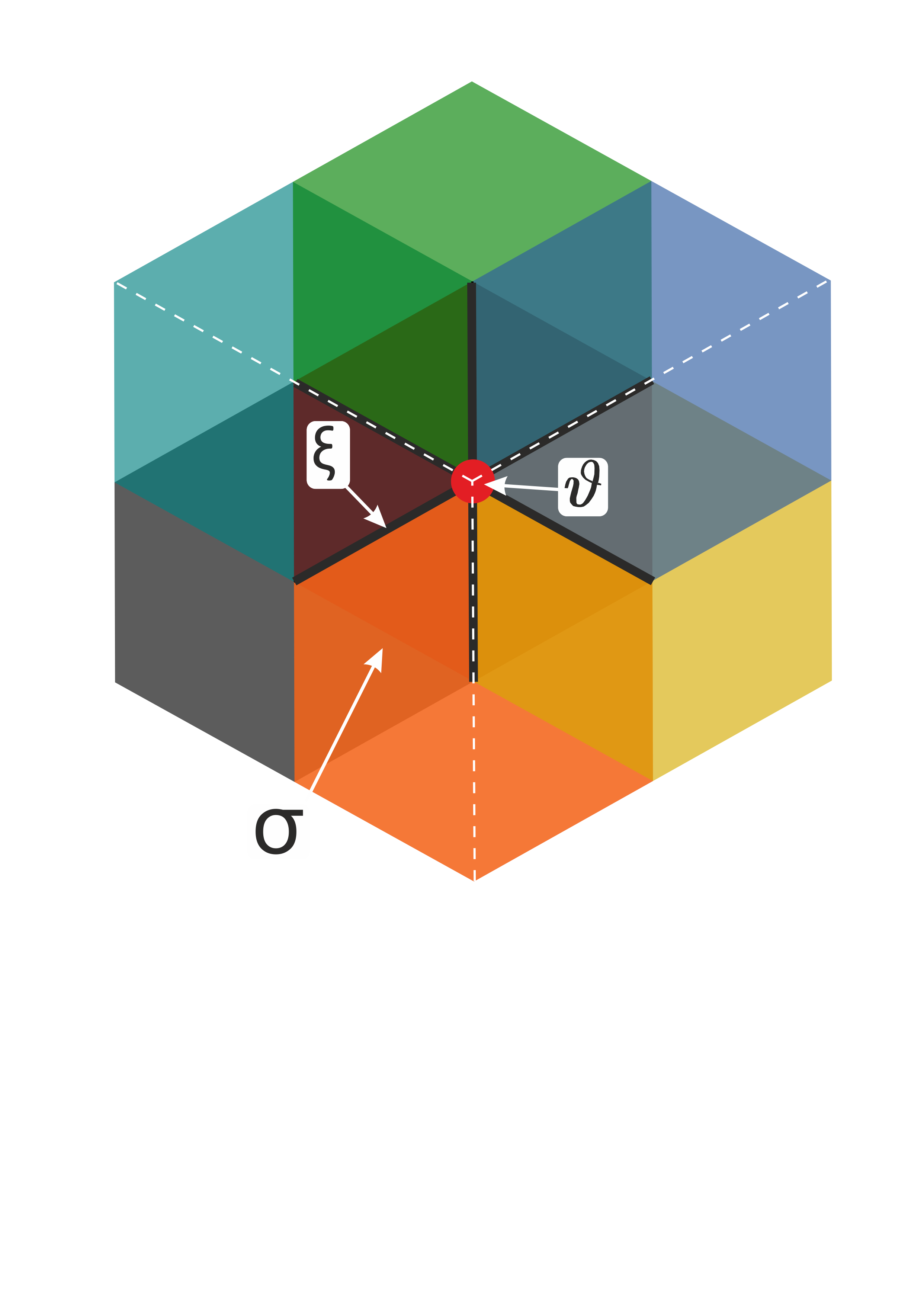}
\caption{Illustration of the interface components of the non-overlapping decomposition into eight subdomains.}
\label{fig:interface}
\end{figure}

The GDSW coarse basis functions are chosen as energy-minimizing extensions of the interface functions $\Phi_{\Gamma}$ to the {\OR interior of the non-overlapping subdomains}. {\OR These extensions can be computed from the assembled finite element matrix.} %
{\FR The interface functions are typically chosen as restrictions of the nullspace of the global Neumann matrix to the vertices $\vartheta$, edges $\xi$, and %
faces $\sigma$ of the non-overlapping decomposition, forming a partition of unity.   Figure \ref{fig:interface} illustrates the interface components for a small 3D decomposition of a cube into eight subdomains.}
{\FR In terms of the interior degrees of freedom ($I$) and the interface degrees of freedom ($\Gamma$) the coarse basis functions can be written as}
\begin{equation} 
\label{eq:phi}
	\Phi
	=
	\begin{bmatrix}
		\Phi_I \\
		\Phi_{\Gamma}
	\end{bmatrix}
	=
	\begin{bmatrix}
		-K_{II}^{-1}K_{I\Gamma}\Phi_\Gamma \\
		\Phi_\Gamma
	\end{bmatrix},
\end{equation}
where $K_{II}$ and $K_{I\Gamma}$ are submatrices of $K${\SP. Here} {\OR $\Gamma$ corresponds to degrees of freedom on the interface of the non-overlapping subdomains 
$\lbrace \Omega_i\rbrace_{i = 1 \ldots N}$
and $I$ corresponds to degrees of freedom in the interior. 
{\OR The algebraic construction of the extensions is based on the partitioning of the system matrix $K$ according to the $\Gamma$ and $I$ degrees of freedom, i.e.,}} %
\begin{equation*}
	K
	=
	\begin{bmatrix}
		K_{II} & K_{I\Gamma} \\
		K_{\Gamma I} & K_{\Gamma\Gamma}
	\end{bmatrix}.
\end{equation*}
Here, $K_{II} = \diag (K_{II}^{(i)})$ is a block-diagonal matrix, where $K_{II}^{(i)}$ defines the $i$-th non-overlapping subdomain.
{\FR The computation of its inverse
{\OR $K_{II}^{-1}$}
can thus be performed independently and {\OR in parallel} %
for all subdomains.} 

{\FR By construction, the number of interface components determines the size of the coarse problem. This number is smaller for the more recent RGDSW methods, which use a reduced coarse space~\cite{Dohrmann:2017:DSC, Rheinbach2018, Heinlein:2022:PST}.
}

For scalar elliptic problems and under certain regularity conditions the GDSW preconditioner allows for a condition number bound
\begin{equation}
	\kappa (M_{\rm GDSW}^{-1} K) \le C \left(1+ \frac{H}{\delta}\right)\left(1+ \log\left(\frac{H}{h}\right)\right),%
\label{eq:cond}
\end{equation}
where $C$ is a constant independent of the other problem parameters; cf.~\cite{Dohrmann:2008:FEM,Dohrmann:2008:DDL}; also
cf.~\cite{dohrmann:2010:HDD} for three-dimensional compressible elasticity.
{\OR Here, $H$ is the diameter of a subdomain, $h$ the diameter of a finite element, and $\delta$ the overlap.}

For three-dimensional almost incompressible elasticity, using adapted coarse spaces,
and
a bound of the form
$$
  \kappa %
  \le C  
  \left(1+\frac{H}{\delta}\right)^3
  \left(1+ \log\left(\frac{H}{h}\right)\right)^2
$$
was established for the GDSW coarse space~\cite{dohrmann:OSA:2009} and also
for a reduced dimensional coarse space~\cite{dohrmann:2010:HDD}.

The more recent reduced dimensional GDSW (RGDSW) coarse spaces~\cite{Dohrmann:2017:DSC} are constructed from nodal interface function, forming a different partition of unity on the interface. The parallel implementation of the RGDSW coarse spaces is {\FR also part of} the FROSch framework; cf.~\cite{Heinlein:2018:IPP}. 
The RGDSW basis function can be computed in different ways.
Here, we use the fully algebraic approach (\textit{Option~1} in~\cite{Dohrmann:2017:DSC}), where the interface values are determined by the multiplicity~\cite{Dohrmann:2017:DSC}; see~\cite{Heinlein:2018:IPP} for a visualization.
Alternatives %
can lead to a slightly lower number of iterations and a faster time to solution~\cite{Heinlein:2018:IPP}, but these
use geometric information~\cite{Dohrmann:2017:DSC,Heinlein:2018:IPP}. 

For problems up to $1000$ cores the GDSW preconditioner {\OR with an exact coarse solver} is a suitable choice. The RGDSW method is able the scale up to $10\,000$ cores. For even larger numbers {\FR of cores and %
subdomains}, a multi-level extension~\cite{Roever:2019,Roever:Three} is available in the FROSch framework.
Although it is not covered by theory, the FROSch preconditioner is sometimes able to scale even if certain dimension of the coarse space
are neglected~\cite{Heinlein:2016:PIT,kupsAlg}, i.e., for linear elasticity the linearized rotations can sometimes be neglected. %
\section{Parallel Software Environment}
\label{sec:implementation}

\subsection{Finite element implementation}
\label{sec:finite-element-implement}
The implementation of the coupled problem by means of the finite element method is based on the incremental two-field potential~\eqref{eq:time-discrete-two-field-principle}, in which the arguments of the local incremental potential $\pi^{\Delta t}$ and the external load functional are expressed by the corresponding finite element approximations.
Introducing the generalized B- and N-matrix in the following manner
\begin{equation}
  \label{eq:gen-b-matrix-minimization}
  \mathcal{Q}=
  \begin{bmatrix}[1.5]
    \nabla\Bvarphi\\
    \bm{J}_{v}\\
    \Divb{\bm{J}_{v}}
  \end{bmatrix}
  =
  \begin{bmatrix}[1.5]
    \underline{\bB}^{\Bvarphi} & \underline{\bzero}\\
    \underline{\bzero} & \underline{\bN}^{\bm{J}_{v}}\\
    \underline{\bzero} & \underline{\bB}^{\Divb{\bm{J}_{v}}}
  \end{bmatrix}
  \begin{bmatrix}[1.5]
    \underline{\widetilde{\Bvarphi}}\\
    \underline{\widetilde{\bm{J}}}_{v}
  \end{bmatrix}
  =\underline{\bB}\,\underline{\bd}
\end{equation}
\begin{equation}
  \label{eq:gen-n-matrix-minimization}
  \mathcal{R}=
  \begin{bmatrix}[1.5]
    \Bvarphi\\
    \bm{J}_{v}
  \end{bmatrix}=
  \begin{bmatrix}[1.5]
    \underline{\bN}^{\Bvarphi}&\underline{\bzero}\\
    \underline{\bzero}&\underline{\bN}^{\bm{J}_{v}}
  \end{bmatrix}
  \begin{bmatrix}[1.5]
    \underline{\widetilde{\Bvarphi}}\\
    \underline{\widetilde{\bm{J}}}_{v}
  \end{bmatrix}=\underline{\bN}\,\underline{\bd}
\end{equation}
and denoting the degrees of freedom of the finite elements by $\widetilde{(\phantom{i})}$, gives rise to the rather compact notation of~\eqref{eq:time-discrete-two-field-principle}
\begin{equation}
  \label{eq:time-discrete-two-field-principle-fem}
  \Pi^{\Delta t,h}\left(\underline{\bd}\right)=\int\nolimits_{\Br} \pi^{\Delta t}\left(\underline{\bB}\,\underline{\bd}\right)\mathrm{d}V  - P^{\Delta t}_{\mathrm{ext}}\left(\underline{\bN}\,\underline{\bd}\right)  \ \mbox{.}
\end{equation}
Upon the subdivision of the domain $\Br$ into finite elements and the inclusion of the assembly operator \textsf{A}, the necessary condition to find a stationary value of the incremental potential is expressed as
\begin{equation}
  \label{eq:residual-time-discrete-two-field-principle-minimization}
  \Pi^{\Delta t, h}_{,\underline{\bd}}= \underline{\bzero} \ \mbox{,}
\end{equation}
which represents a system of nonlinear equations
\begin{equation}
  \label{eq:residual-time-discrete-two-field-principle}
  \underline{\bR}\left(\underline{\bd}\right)=  \begin{bmatrix}[1.5]
    \underline{\bbr}_{\Bvarphi}\left(\underline{\widetilde{\Bvarphi}},\underline{\widetilde{\bm{J}}}_{v}\right)\\
    \underline{\bbr}_{\bm{J}_{v}}\left(\underline{\widetilde{\Bvarphi}},\underline{\widetilde{\bm{J}}}_{v}\right)
  \end{bmatrix}
  \stackrel{!}{=} \underline{\bzero}
\end{equation}
with
\begin{multline}
  \label{eq:mech-residual-time-discrete-two-field-principle}
  \underline{\bbr}_{\Bvarphi}=\Assem\int\limits_{\Br_{e}}\underline{\bB}^{\Bvarphi\mathrm{T}}\ps{\bF}\widehat{\psi} - \underline{\bN}^{\Bvarphi\mathrm{T}} \bm{R}_{\Bvarphi}\mathrm{d}V \\
  - \Assemarg{e=1}{k_{\mathrm{ele}}} \int \limits_{\partial \Br^{\bm{T}}_{e}} \underline{\bN}^{\Bvarphi\mathrm{T}}\bar{\bm{T}}\,\mathrm{d}A
\end{multline}
and
\begin{multline}
  \label{eq:diffusion-residual-time-discrete-two-field-principle}
  \underline{\bbr}_{\bm{J}_{v}}=\Assem\int\limits_{\Br_{e}}\!\Delta t \left[\underline{\bN}^{\bm{J}_{v}\mathrm{T}} \ps{\bm{J}_{v}}\widehat{\phi} - \underline{\bB}^{\Divb{\bm{J}_{v}}\mathrm{T}} \ps{v}\widehat{\psi}\right]\mathrm{d}V\\ + \Assemarg{e=1}{l_{\mathrm{ele}}} \int \limits_{\partial\Br^{\mu}_{e}} \!\Delta t\, \bar{\mu}\, \underline{\bN}^{\bm{J}_{v}\mathrm{T}}\bm{N} \,\mathrm{d}A \ \mbox{.}
\end{multline}
Equation~\eqref{eq:residual-time-discrete-two-field-principle} is solved efficiently by means of a monolithic Newton-Raphson scheme.
The corresponding linearization is inherently symmetric and is computed as
\begin{equation}
  \label{eq:stiffness-matrix-minimization}
  \underline{\bK}=\Pi^{\Delta t, h}_{,\underline{\bd}\,\underline{\bd}}=
  \begin{bmatrix}[1.5]
    \underline{\bK}_{\Bvarphi\,\Bvarphi}&\underline{\bK}_{\Bvarphi\,\bm{J}_{v}}\\
    \underline{\bK}_{\bm{J}_{v}\,\Bvarphi}&\underline{\bK}_{\bm{J}_{v}\,\bm{J}_{v}}
  \end{bmatrix} \, \mbox{,}
\end{equation}
in which the individual contributions take the form
\begin{align}
  \label{eq:stiffness-matrix-minimization-mech-mech}
  \underline{\bK}_{\Bvarphi\,\Bvarphi}&=\Assem\int\limits_{\Br_{e}}\underline{\bB}^{\Bvarphi\mathrm{T}}\psq{\bF\bF}\widehat{\psi}\,\underline{\bB}^{\Bvarphi}\mathrm{d}V
  \\
  \label{eq:stiffness-matrix-minimization-mech-flux}
  \underline{\bK}_{\Bvarphi\,\bm{J}_{v}}&=\Assem\int\limits_{\Br_{e}}-\Delta t
                                          \,\underline{\bB}^{\Bvarphi\mathrm{T}}\psq{\bF
                                          v}\widehat{\psi}\,\underline{\bB}^{\Divb{\bm{J}_{v}}}\mathrm{d}V
  \\
  \label{eq:stiffness-matrix-minimization-flux-mech}
  \underline{\bK}_{\bm{J}_{v}\,\Bvarphi}&=\Assem\int\limits_{\Br_{e}}-  \Delta t \,\underline{\bB}^{\Divb{\bm{J}_{v}}\mathrm{T}} \psq{v\bF}\widehat{\psi}\,\underline{\bB}^{\Bvarphi}\mathrm{d}V\\
  \label{eq:stiffness-matrix-minimization-flux-flux}
  \underline{\bK}_{\bm{J}_{v}\,\bm{J}_{v}}&=\Assem\int\limits_{\Br_{e}}\!\Delta t
                                            \left[\underline{\bN}^{\bm{J}_{v}\mathrm{T}}
                                            \psq{\bm{J}_{v}\bm{J}_{v}}\widehat{\phi}\,\underline{\bN}^{\bm{J}_{v}} \right.\nonumber\\
  &\quad \left.+ \Delta t \,
    \underline{\bB}^{\Divb{\bm{J}_{v}}\mathrm{T}}
    \psq{vv}\widehat{\psi}
    \,\underline{\bB}^{\Divb{\bm{J}_{v}}}\right]\mathrm{d}V
    \ \mbox{,}
\end{align}
{\OR where $\widehat{\psi}$ is the hyperelastic energy associated with %
the structural problem
and $\widehat{\phi}$ the dissipation potential corresponding to the diffusion problem; see Section~\ref{sec:specific-free-energy-funct-dissipation-potential}}.
The implementation of the model is carried out {\OR using} the finite element library deal.II~\cite{dealii2019design}, employing tri-linear Lagrange ansatz functions for the deformation, while two different approaches have been chosen to approximate the fluid flux:
First, also a tri-linear Lagrange ansatz function is used for the flux variable, which 
{\OR is not the standard conforming discretization}
{\SP but was nevertheless successfully applied in the context of diffusion-induced fracture of hydrogels~\cite{BogKeiMie:2017:min}}.
Second, the lowest order, conforming Raviart-Thomas ansatz is selected, ensuring the continuity of the normal trace of the flux field across element boundaries.

\begin{figure}[h]
  \centering
  \begin{overpic}[scale=1.0]
    {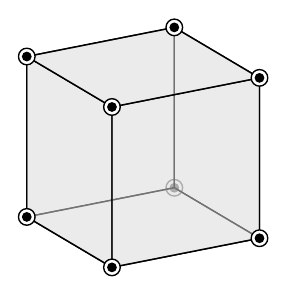}
    \put(0,0){(i)}
  \end{overpic}
  \hspace{8ex}
  \begin{overpic}[scale=1.0]
    {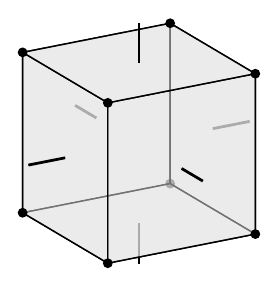}
    \put(0,0){(ii)}
  \end{overpic}
  \caption{Two different finite element combinations employed in the current study. (i) Lagrange/Lagrange ansatz functions {\SP \textit{$Q_{1} Q_{1}$}} and (ii) Lagrange/Raviart-Thomas ansatz functions {\SP \textit{$Q_{1} RT_{0}$}}. Note that, vectorial degree of freedom associated with the deformation are indicated by \DefDof, while vectorial fluid flux degree of freedom are illustrated as \FluxVecDof and scalar normal traces of the flux field are shown as thick solid lines (\NormalTrDof).}
  \label{fig:two-different-finite-element-combinations-lagr-lagr-raviart-thomas}
\end{figure}

In the following, we denote the tri-linear Lagrange ansatz functions by \textit{{\FR $Q_1$}} and the Raviart-Thomas ansatz functions of lowest order by \textit{{\FR $RT_0$}}. {\FR %
We {\OR then} denote the combination for the deformation and flux elements by $Q_1 Q_1$ {\OR and} $Q_1 RT_0$.}
Both element combinations are fully integrated numerically by means of a Gauss quadrature.
They are depicted in Figure~\ref{fig:two-different-finite-element-combinations-lagr-lagr-raviart-thomas}.
\subsection{Linearized monolithic system}
\label{sec:linear-monolith-system}
{\SP
For completeness we state the linearized monolithic system of equations that has to be solved at each iteration $k$ of the Newton-Raphson scheme as
\begin{equation}
  \label{eq:monolithic-linear-system-NR-iter}
    \begin{bmatrix}[1.5]
    \underline{\bK}_{\Bvarphi\,\Bvarphi}&\underline{\bK}_{\Bvarphi\,\bm{J}_{v}}\\
    \underline{\bK}_{\bm{J}_{v}\,\Bvarphi}&\underline{\bK}_{\bm{J}_{v}\,\bm{J}_{v}}
  \end{bmatrix}_{k}
  \begin{bmatrix}[1.5]
    \vartriangle\!\Bvarphi\\
    \vartriangle\!\bm{J}_{v}
  \end{bmatrix}
  =-
  \begin{bmatrix}[1.5]
    \underline{\bbr}_{\Bvarphi}\\
    \underline{\bbr}_{\bm{J}_{v}},
  \end{bmatrix}_{k}
\end{equation}
where $\underline{\bK}_{\bm{J}_{v}\,\Bvarphi}\!=\underline{\bK}^{\mathrm{T}}_{\Bvarphi\,\bm{J}_{v}}$,
to update the degrees of freedom associated with the deformation as well as the flux field according to
\begin{equation}
  \label{eq:NR-update}
  \begin{split}
    \Bvarphi_{k+1} &= \Bvarphi_{k} + \vartriangle\!\Bvarphi\\
    \bm{J}_{v \,k+1} & = \bm{J}_{v \, k} + \vartriangle\!\bm{J}_{v} \ \mbox{.}
  \end{split}
\end{equation}
The convergence criteria employed in this study are outlined in Table~\ref{tab:tol-newton}.
}

{\OR Since our preconditioner is constructed algebraically, it is important to note that 
in deal.II the ordering of the degrees of freedom for the two-field problem is different for different discretizations.

  In the case of a $Q_{1} Q_{1}$ discretization, a node-wise numbering is used, %
  and the %
  global vector thus has the form %
  \begin{equation}
    \label{eq:structure-dof-vector-q1q1}
    \underline{\bd} = [\dots,\underbrace{\varphi_{1},\varphi_{2},\varphi_{3}, J_{v1}, J_{v2}, J_{v3}}_{\mbox{node $p$}} ,\dots]^{\mathrm{T}} \ \mbox{.}
  \end{equation}
  On the contrary, in the $Q_{1} RT_{0}$ discretization, all degrees of freedom associated with the deformation %
  are arranged first, followed by the flux degrees of freedom. %
  Thus the %
  global vector thus takes the form
  \begin{equation}
    \label{eq:structure-dof-vector-q1rt0}
    \underline{\bd}=[\dots,\underbrace{\varphi_{1},\varphi_{2},\varphi_{3}}_{\mbox{node $p$}},\dots,\underbrace{H_{v}}_{\mbox{face $q$}},\dots]^{\mathrm{T}} \ \mbox{.}
  \end{equation}
}
\subsubsection{Free-swelling boundary value problem}
\label{sec:free-swell-bvp}
\begin{figure}[h]
  \centering
  \includegraphics[scale=1.0]{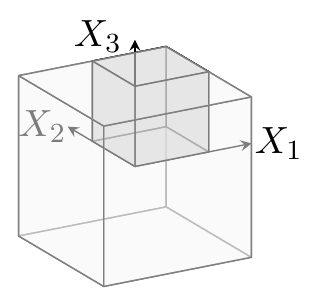}
  \caption{One eighth of the cube domain, highlighted in dark gray, with edge length $\mathrm{L}=\unit[1]{mm}$ employed in the parallel scalability study.}
  \label{fig:one-eighth-cube-free-swelling}
\end{figure}

\begin{figure}[h]
  \centering
  \includegraphics[scale=0.8]{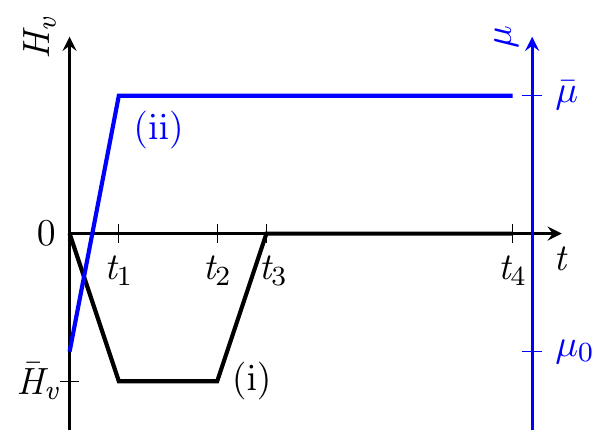}
  \caption{Time-dependent boundary conditions for the free-swelling problem: (i) flux control and (ii) control through chemical potential.}
  \label{fig:time-depend-bound-cond-free-swelling}
\end{figure}

The boundary value problem of a free-swelling cube is studied in the literature, cf.~\cite{BogNatMie:2017:min} for 2D and~\cite{Mau:2017:var} for 3D results, and adopted here as a benchmark problem for the different finite element combinations.
{\SP
Considering a cube with edge length $2\mathrm{L}$, the actual simulation of the coupled problem is carried out employing only one eighth of the domain, as shown in Figure~\ref{fig:one-eighth-cube-free-swelling}, due to the intrinsic symmetry of the problem.
Therefore, symmetry conditions are prescribed along the three symmetry planes, i.e. $X_{1}=0, X_{2}=0, X_{3}=0$, which correspond to a vanishing normal component of the displacement vector and the fluid flux.
At the outer surface the mechanical boundary conditions are assumed as homogeneous Neumann conditions, i.e. $\bar{\bm{T}}=\bm{0}$, while two different boundary conditions are used for the diffusion problem, namely
\begin{itemize}
\item[(i)] Dirichlet conditions, i.e. the normal component of the fluid flux $H_{v}$, are prescribed or
\item[(ii)] Neumann conditions, i.e. the chemical potential $\bar{\mu}$, are specified as shown in Figure~\ref{fig:time-depend-bound-cond-free-swelling} and Table~\ref{tab:param-bound-cond-free-swelling}.
\end{itemize}
Note that, due to the coupling of mechanical and diffusion fields, the boundary conditions (i) and (ii) result in different temporal evolution of the two fields.
However, in both cases a homogeneous, stress-free state is reached under steady state conditions.
}
{\SP
  Type (i) boundary conditions are used for {\FR strong} scalability study outlined in Section~\ref{subsec:free-strong}, while type (ii) boundary conditions are employed in the weak parallel scalability study described in Section~\ref{sec:weak-scaling}.
}
\begin{table}[h]
  \centering
  \begin{tabular}[h]{ccccc}
    type&$\bar{H}_{v}$ / $\bar{\mu}$&$\mu_{0}$&$t_{2}$&$t_{3}$\\
    \hline
    (i)&$\unitfrac[-0.02]{mm}{s}$&-&$\unit[0.75]{s}$&$\unit[1]{s}$\\
    (ii)&$\unit[-40.31]{MPa}$&$\unit[-80.31]{MPa}$&-&-
  \end{tabular}
  \caption{Problem specific parameters associated with the boundary conditions in the free-swelling problem illustrated in Figure~\ref{fig:time-depend-bound-cond-free-swelling}. The common parameters are $t_{1}=\unit[0.25]{s}$ and $t_{4}=\unit[4]{s}$}
  \label{tab:param-bound-cond-free-swelling}
\end{table}
\subsubsection{Mechanically induced diffusion boundary value problem}
\label{sec:mech-induced-diffusion-bvp}
\begin{figure}[h]
  \centering
  \includegraphics[scale=1.0]{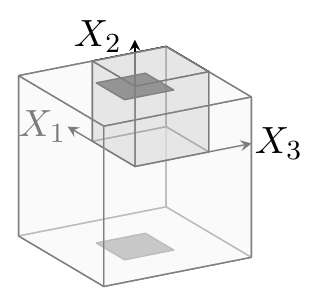}
  \caption{One eighth of the cube domain, highlighted in dark gray, with edge length $\mathrm{L}=\unit[1]{mm}$ employed in the parallel scalability study. The shaded dark gray surface indicates the area with prescribed vertical displacement, mimicking the indentation of the body with a rigid flat punch.}
  \label{fig:one-eighth-cube-mech-induced-diffusion}
\end{figure}

{\SP
  Similar to the free-swelling problem, the mechanically induced diffusion problem is also solved on a unit cube domain with appropriate symmetry conditions applied along the planes $X_{1}=0, X_{2}=0, X_{3}=0$, as shown in Figure~\ref{fig:one-eighth-cube-mech-induced-diffusion}.
  Along the subset $\left(X_{1},X_{3}\right)\in\left[-\frac{\mathrm{L}}{3},\frac{\mathrm{L}}{3}\right]\times\left[-\frac{\mathrm{L}}{3},\frac{\mathrm{L}}{3}\right]$ at the plane $X_{2}=\mathrm{L}$ the coefficients of the displacement vector are prescribed as $u_{i}=[0,-\hat{u},0]$, mimicking the indentation of the body with a rigid flat punch under very high friction condition.
  The non-vanishing displacement coefficient is increased incrementally and subsequently held constant, similar to the function for the chemical potential illustrated in Figure~\ref{fig:time-depend-bound-cond-free-swelling}.
  The corresponding parameters read as $t_{1}=\unit[1]{s}$, $t_{4}=\unit[6]{s}$ and $\hat{u}=\unit[0.4]{mm}$.
  Additionally, the normal component of the fluid flux is set to zero at the complete outer surface of the cube together with traction free conditions at the remaining part of the outer boundary.
}

\subsection{Distributed memory parallelization using deal.ii, p4est, and Trilinos}
In early versions of our software, the assembly of the system matrix and
the Newton steps were performed on a single core, and after distribution of the system matrix to all cores, it was solved by
FROSch in fully algebraic mode~\cite{Kiefer:PAM:2021}.

In this work, the simulation is fully MPI-parallel.
We assemble the system matrix in parallel using the deal.II classes from the \texttt{parallel::distributed} namespace. %
In deal.II, parallel meshes for distributed memory machines %
are handled by a \texttt{parallel::distributed::Triangulation} object, which  %
calls the external library \textit{p4est}~\cite{p4est} to determine the parallel layout of the mesh data.

As a result, each process owns a portion of cells (called \texttt{locally owned cells} in deal.II) of the global mesh.
Each process stores one additional layer of cells surrounding the \texttt{locally owned cells}, which are denoted as \textit{ghost cells}.
{\OR Using the ghost cells two MPI ranks corresponding to neighboring nonoverlapping subdomains can, both, access (global) degrees of freedom on the interface.}

Each local stiffness matrix is assembled by the process which owns the associated cell (i.e., the finite element), thus the processes work independently and concurrently. The handling of the parallel data (cells and degrees of freedom) distribution is performed by a \texttt{DofHandler} object. A more detailed describtion of the MPI parallelization in deal.II can be found in~\cite{dealiiParallel}.

For the parallel linear algebra, deal.II interfaces to either the PETSc~\cite{petsc-web-page} or the Trilinos~\cite{trilinosrepo}. In this work, we make use of the classes in the \texttt{dealii::LinearAlgebraTrilinos::MPI} namespace, 
such that we obtain Trilinos \texttt{Epetra} vectors and matrices, which can be processed by FROSch.
Similarly to the \texttt{DofHandler} the Trilinos \texttt{Map} object handles the data distribition of the parallel linear algebra objects. 

To construct the coarse level, FROSch needs information on the interface between the subdomains. The FROSch framework uses a repeatedly decomposed \texttt{Map} to identify the interface components. In this \texttt{Map} the degrees of freedom on the interface are shared among the relevant processes. This \texttt{Map} can be provided as an input by the user. However, FROSch also provides a fully algebraic construction of the repeated Map~\cite{Heinlein:2018:FPI}, which is what we use here.

\subsection{Solver settings}
\label{sec:solver-settings}
\begin{figure}[]
\centering
\includegraphics[width=0.3\textwidth, trim={40mm 0mm 40mm 0mm},clip]{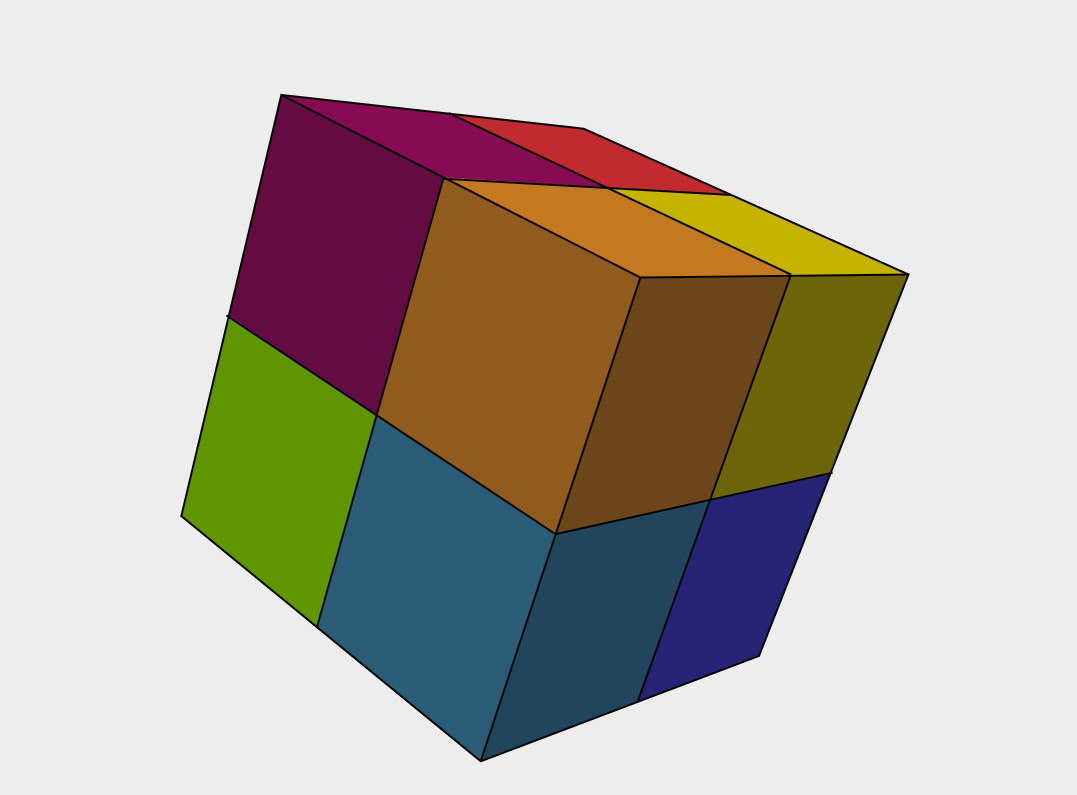}
\includegraphics[width=0.3\textwidth, trim={40mm 0mm 40mm 0mm},clip]{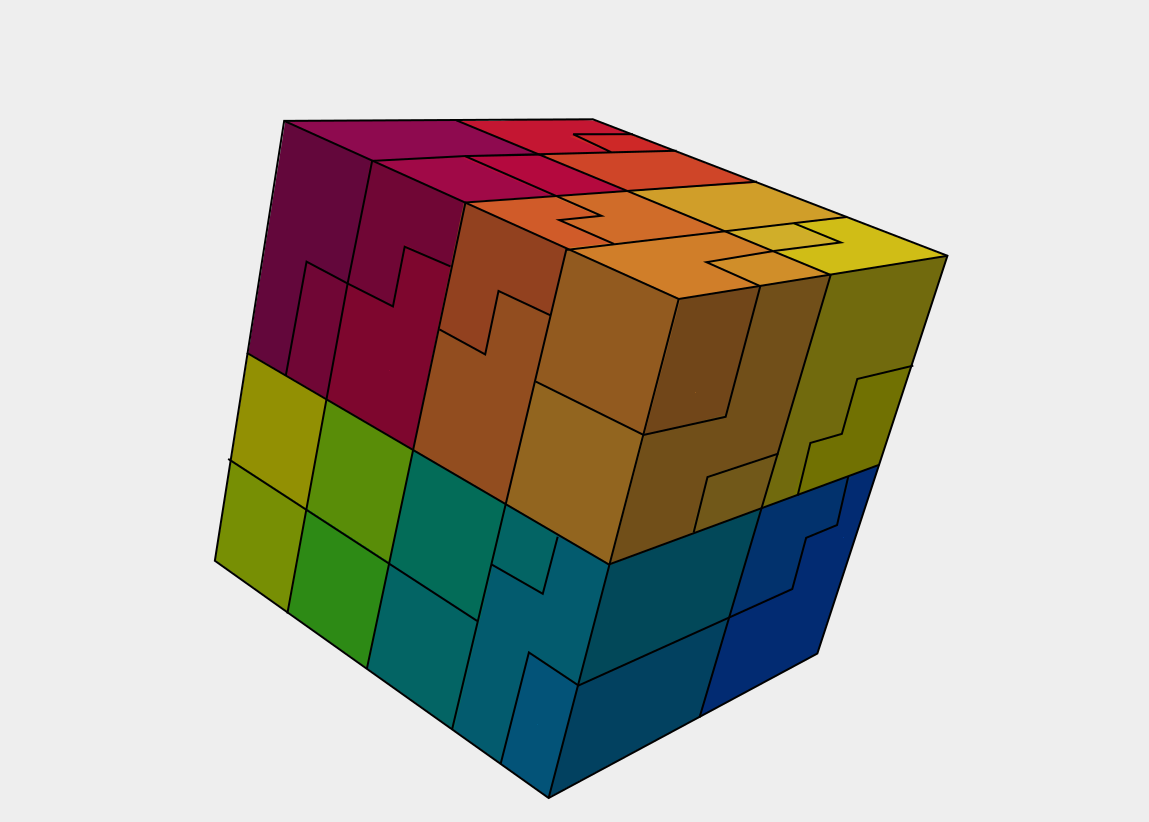}
\includegraphics[width=0.3\textwidth, trim={40mm 0mm 40mm 0mm},clip]{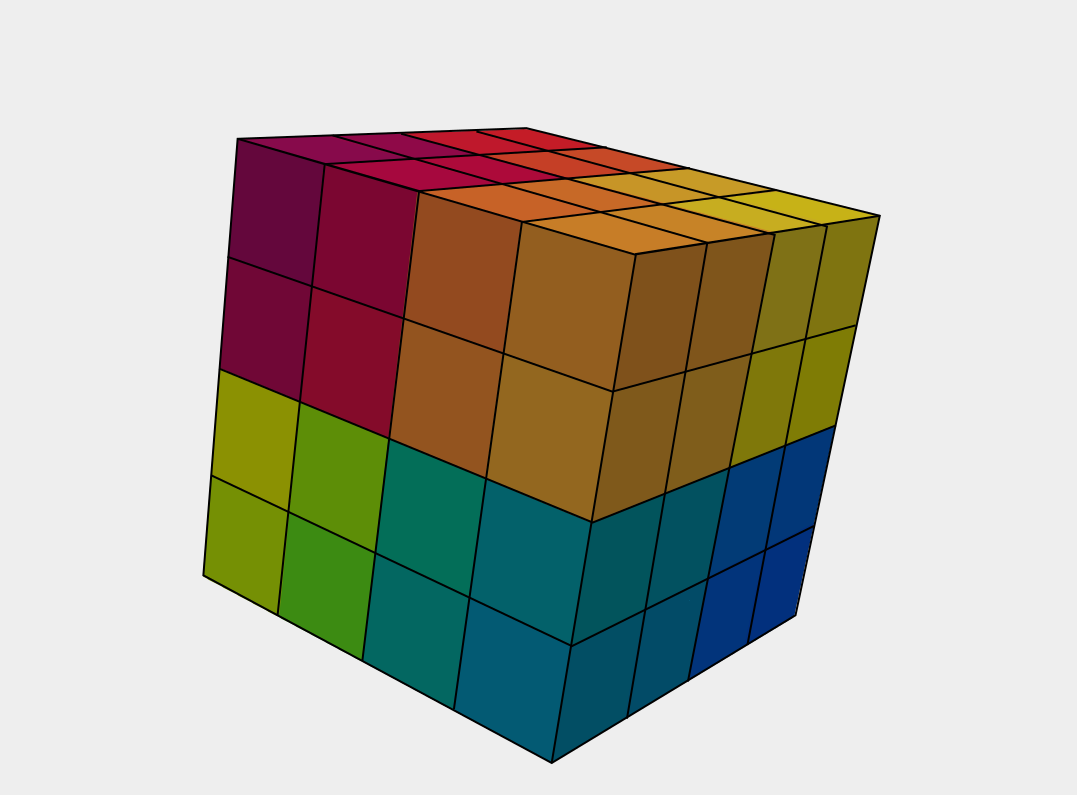}
\caption{\FR Decomposition of the computational domain with 4\,096 finite elements %
into 8 (top), 27 (middle) and 64 (bottom) non-overlapping subdomains.}
\label{fig:diffDecom}
\end{figure}

{\OR We use the deal.II software library (version 9.2.0)~\cite{dealII92, dealii2019design} to implement the model in variational form, and to perform the finite element assembly in parallel.}
The parallel decomposition of the computational domain is performed in deal.II by using the \textit{p4est} software library~\cite{p4est}.
We remark  that using \textit{p4est} %
{\OR small changes in}
the number of {\OR finite elements} %
and 
the number of subdomains %
{\OR may result in decompositions with very different %
subdomain shapes}; see %
Figure~\ref{fig:diffDecom}. A bad subdomain shape will typically degrade the convergence of the domain decomposition solver. %
{\FR 
  We always choose an overlap of two elements.
} %
{\OR   However, since the overlap is constructed algebraically in some positions there can be deviations from a geometric overlap of $\delta=2h$.}

{\FR
{\OR In %
the Newton-Raphson scheme,} we use absolute and relative tolerances for the deformation  $\Bvarphi$ and the fluid flux $\bm{J}_{v}$  according to Table~\ref{tab:tol-newton}, where $\Vert r_k \Vert$ is the residual at the k-th Newton step and $\Vert r_0 \Vert$ the initial residual.
\begin{table}
\centering
\begin{tabular}{l | c | c}
 & $\Vert r_k \Vert$ & $\Vert r_k \Vert/\Vert r_0 \Vert$ \\ \hline
 $\Bvarphi$ & $10^{-9}$ & $10^{-6}$ \\
 $\bm{J}_{v}$ & $5 *10 ^{-12}$ & $10^{-9}$ \\ 
\end{tabular}
\caption{%
\OR Tolerances for %
the Newton-Raphson scheme. Here, $r_k$ is the $k$-th residual.}
\label{tab:tol-newton}
\end{table}
}

{\OR FROSch is part of Trilinos~\cite{trilinosrepo} and makes heavy use of the parallel Trilinos infrastructure.}
The Trilinos software library is applied using the master branch of October 2021~\cite{trilinosrepo}.

{\FR %

On the first level of overlapping subdomains we always apply the restrictive additive Schwarz method.
A one-to-one correspondence for the subdomains and cores is employed. 

{\OR The linearized systems are solved} using the parallel GMRES implementation provided by the Trilinos package Belos using the relative stopping criterion of $\Vert r_k \Vert/\Vert r_0 \Vert \leq 10^{-8}$. 
{\FR We use a vector of all zeros as the initial vector for the iterations.}

The arising subproblems in the FROSch framework are solved by %
Trilinos' {\OR built-in} KLU sparse direct linear solver. %
 
{\OR All parallel experiments are performed on the Compute Cluster of the Fakult\"at f\"ur Mathematik und Informatik at Technische Universit\"at Freiberg. A cluster node has two Intel Xeon Gold 6248 processors (20 cores, 2.50 GHz).}

\label{sec:numRes}

{\OR
\section{Limitations of this study}
This work is a first step towards a co-design of variational formulation, finite element discretization, and parallel iterative solver environment. However, this work has some limitations, which we now briefly discuss.
\subsection{Coupling constraint}
\label{sec:coupling-constraint}
Our model incorporates a coupling constraint 
$
J J_0 -1 - v = 0
$
which couples the volumetric deformation of the structure to the fluid flux $v$; {\SP see Section~\ref{sec:specific-free-energy-funct-dissipation-potential};}
this constraint is implemented using a quadratic penalty
$
\frac{\uplambda}{2J_{0}}\left[JJ_{0}-1- v\right]^{2},
$
where $\uplambda$ is the penalty parameter. For $v=0$ this constraint corresponds to an (almost) incompressibility constraint.
It is well known that, for almost incompressible elasticity, low-order standard Lagrange finite elements can result in locking, and stable discretizations should be used if the penalty parameter $\uplambda$ is high. A standard technique is the use of a three-field formulation known also as the ${\bf \bar{F}}$-method. In the ${\bf \bar{F}}$-method, the point-wise constraint is relaxed, and the penalty enforces the constraint only in a mean sense. In Section~\ref{sec:stability-finite-element}, where the penalty parameter $\uplambda$ is chosen as 
$\uplambda=\unitfrac[10]{N}{mm^2}$, we observe oscillations which indicate the instability of the standard $Q1$-Lagrange discretization for the structure. However, for our numerical results in Section~\ref{sec:perf}, a much lower penalty parameter $\uplambda=\unitfrac[0.2]{N}{mm^2}$ is used. 
For $v=0$, this corresponds to a fairly compressible material.
For this value, no stability problems were observed experimentally.
\subsection{Penalty formulation for the coupling constraint}
\label{sec:penalty-formulation-coupling-constraint}
A penalty formulation affects the conditioning of the stiffness matrix which can degrade the performance of iterative solvers. When direct solvers are used the quality of the solution will degrade if the penalty parameter is high.
The penalty parameter used in our numerical experiments in Section~\ref{sec:numer} 
is mild.
For a larger penalty parameter ill-conditioning of the stiffness matrix could be avoided, e.g., by using an Augmented Lagrange approach as in~\cite{brinkhues:2013:ALM}, at the cost of an additional outer iteration.
\subsection{Discretization of the fluid flux}
\label{sec:discr-fluid-flux}
In our experiments using a $Q_1Q_1$ finite element discretization, the fluid flow is approximated by standard $Q1$ elements. %
Due to the stronger patching condition compared to the $H(\div)$-conforming Raviart-Thomas discretizations certain solutions can not be approximated well using $Q_1$ elements; see, e.g., \cite{steeger:2017:phd} for a detailled discussion.
However, such discretizations have been considered also elsewhere, e.g., in the context of least-squares formulations~\cite{schwarz:2022:3MG}, where no LBB condition has to be fulfilled; see also~\cite{schwarz:2014:WOL}, where different combinations of Raviart-Thomas elements with Lagrange elements have been considered for least-squares formulations}
{\SP of the Navier-Stokes problem in the almost incompressible case.
  Considering coupled problems involving mechanics and diffusion, the hybridization technique introduced in~\cite{YuMalHuy:2020:com} is promising, where the incompressibility constraint of the fluid is exactly enforced by altering the balance of solute volume~\eqref{eq:field-equation-balance-solute-volume} as a consequence of the constant volume of the solid, i.e. the fluid flux is directly coupled to the overall volume change. Furthermore, incorporating refined material models based on the theory of porous media, c.f.~\cite{TeiMauMie:2019:asp}, in principle allows one to enforce the incompressibility of the solid and the fluid independently.
  However, these types of models and the corresponding finite element formulations are beyond the scope of the current contribution, and we leave it for future work.
}
{\OR
\subsection{Sparse direct solver}
\label{sec:sparse-direct-solver}
We apply the KLU sparse direct solver~\cite{KLU} which is the default solver of the Trilinos Amesos2 package~\cite{amesos2}. It is used as a subdomain local sparse solver and as solver for the coarse problem. KLU was originally developed for problems from circuit simulation and is therefore often not the fastest choice available for the factorization of finite element matrices, especially if they are large. However, it was found that it can sometimes outperform other established sparse direct solvers also for matrices from finite element problems; see~\cite{benzi:2016:PER}.
In the future, other well-known fast sparse direct solvers should, however, also be tested in our context.
\subsection{Fully algebraic construction of the preconditioner}
\label{sec:limitalg}
In this work, it is our approach to use the two-level version of the preconditioners in FROSch in fully algebraic mode, i.e., without making use of the structure of %
the problem when constructing the preconditioner from the monolithic finite element matrix.
Therefore, the second level is constructed in the same way as for a scalar elliptic problem. 
It is known that this coarse level is not the correct choice for elasticity.

As a result, numerical scalability in the strict sense should not be expected, i.e., the number of iterations can be expected to grow with the number of cores. 
We will investigate how fast the number of iterations grow, and we will show experiments which indicate that we may be able to achieve numerical scalability by using some additional information on the problem in the construction of the preconditioner.
} %

\subsection{Stability of the finite element formulations}
\label{sec:stability-finite-element}
{\SP
  {\OR We briefly discuss the stability of the %
  finite element discretization, %
  namely the $Q_{1} Q_{1}$ and $Q_{1} RT_{0}$ ansatz; see also Section~\ref{sec:discr-fluid-flux}} %
  {\OR Here,} the mechanically induced diffusion problem, described in Section~\ref{sec:mech-induced-diffusion-bvp}, is solved for the complete loading history, %
  employing a discretization with {$24^{3}$} finite elements, resulting in $93\,750$ degrees of freedom for the $Q_{1} Q_{1}$ and $90\,075$ degrees of freedom for the $Q_{1} RT_{0}$ ansatz function.
  Note that this discretization corresponds to only one uniform refinement step less than the discretizations considered in Section~\ref{sec:strong-scaling} and \ref{sec:weak-scaling}.
  The material parameters in this study are taken as $\uplambda=\unitfrac[10]{N}{mm^2}$ and $J_{0}=4.5$, while all the remaining parameter are chosen %
  {\OR according to} Table~\ref{tab:material-parameter}. The penalty parameter $\uplambda$ is thus larger by a factor of $50$ compared to Table~\ref{tab:material-parameter}.
    \begin{figure*}[h]
    \centering
    \begin{overpic}[scale=0.13]
      {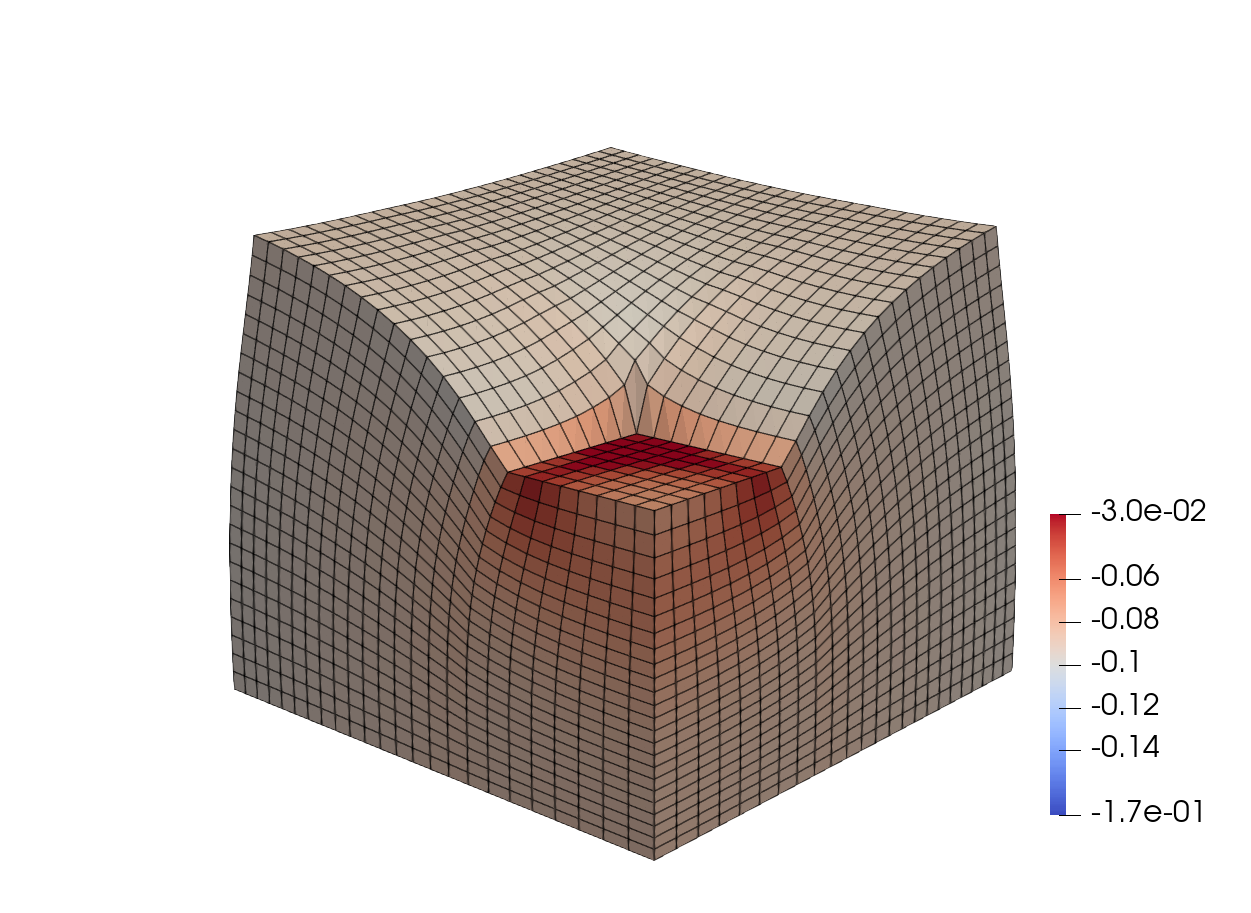}
      \put(5,5){(i)}
    \end{overpic}
    \hspace{-8ex}
    \begin{overpic}[scale=0.13]
      {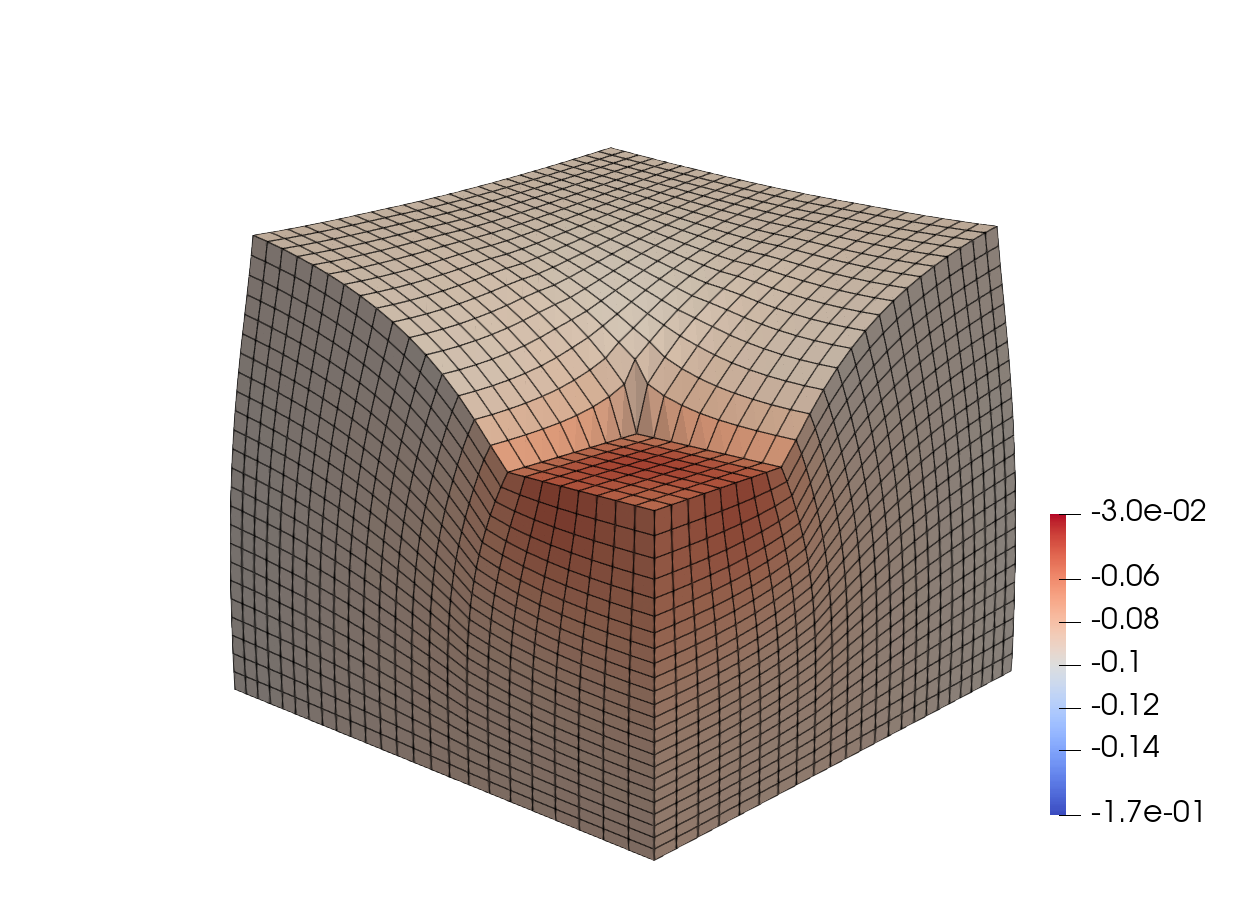}
    \end{overpic}
    \hspace{-8ex}
    \begin{overpic}[scale=0.13]
      {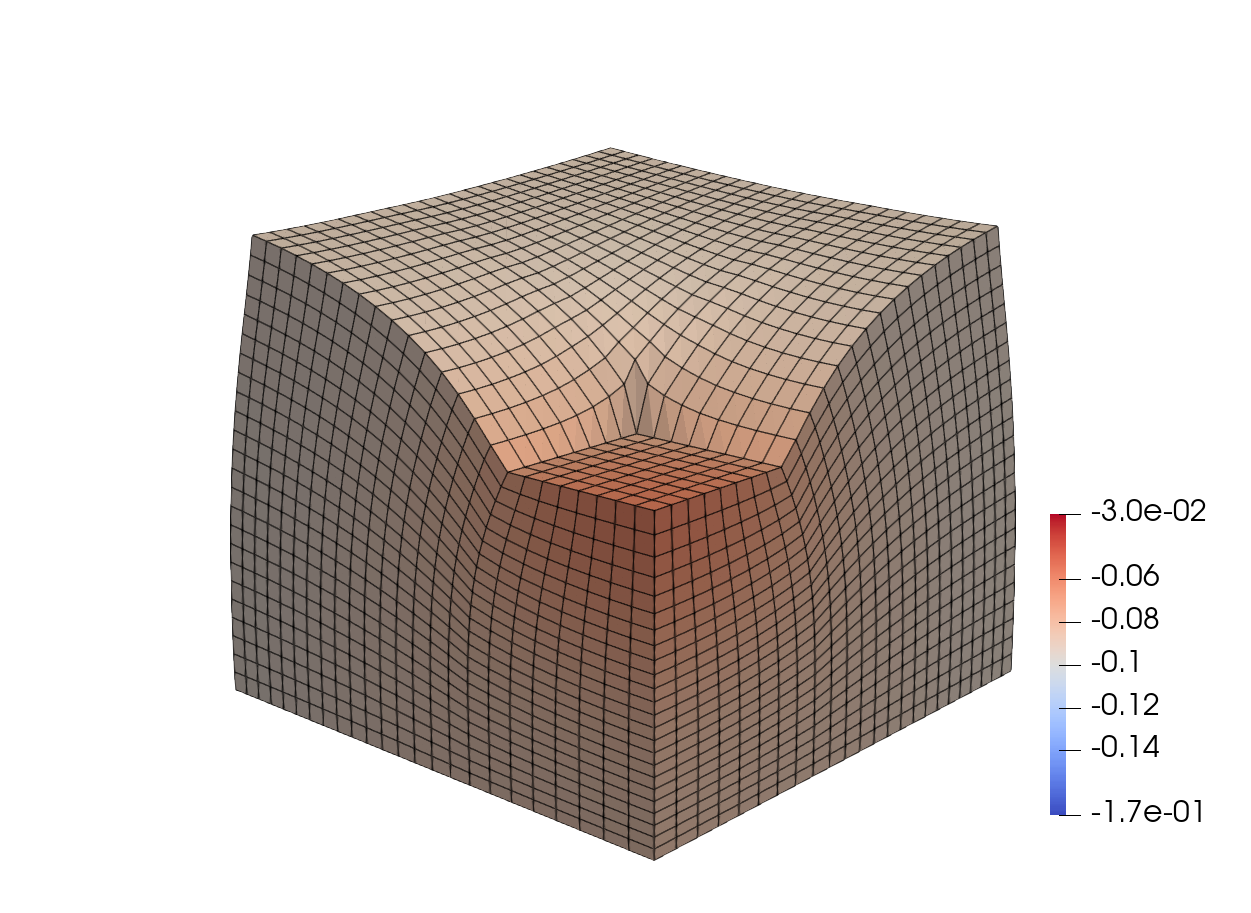}
      \put(95,20){\tiny $\mu$}
    \end{overpic}
    \begin{overpic}[scale=0.13]
      {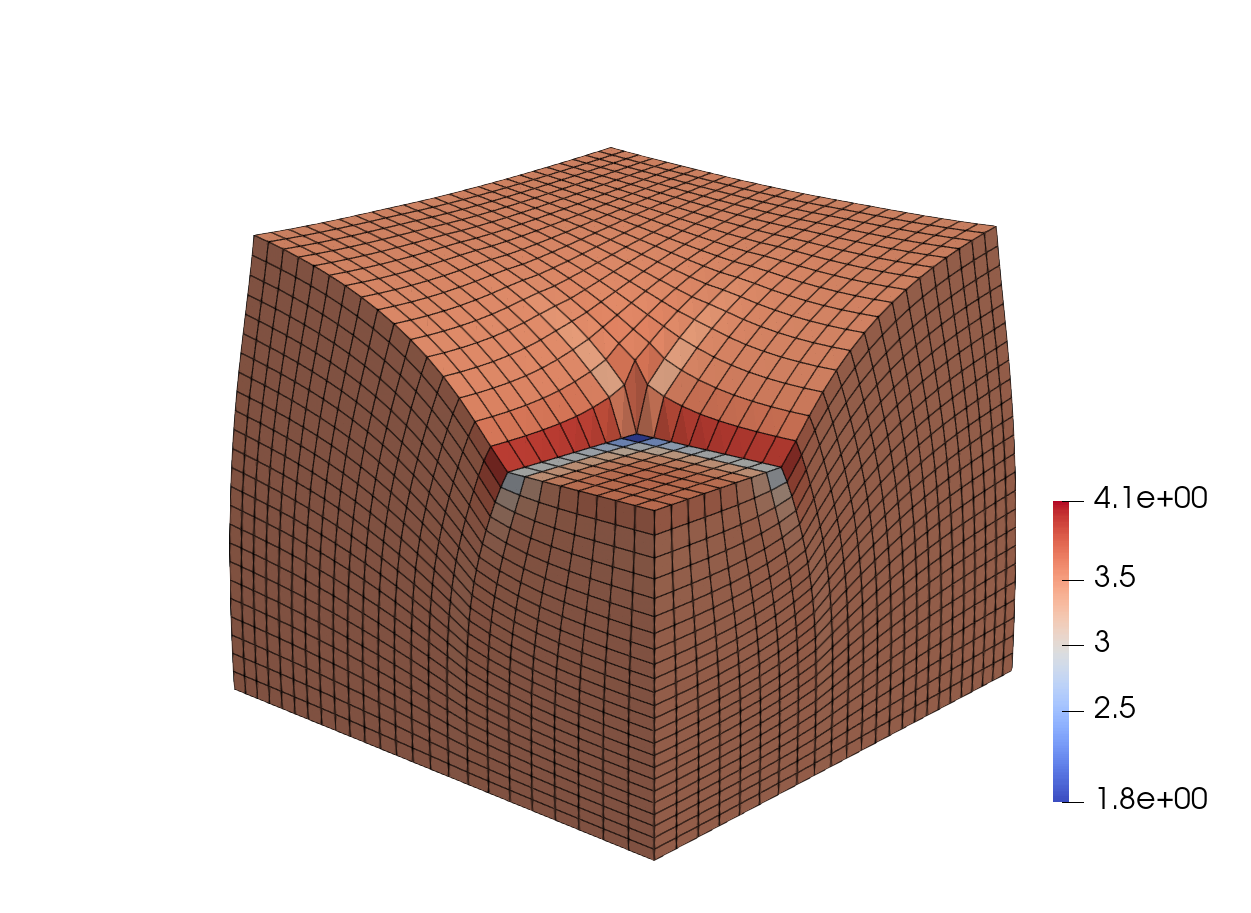}
      \put(5,5){(ii)}
    \end{overpic}
    \hspace{-8ex}
    \begin{overpic}[scale=0.13]
      {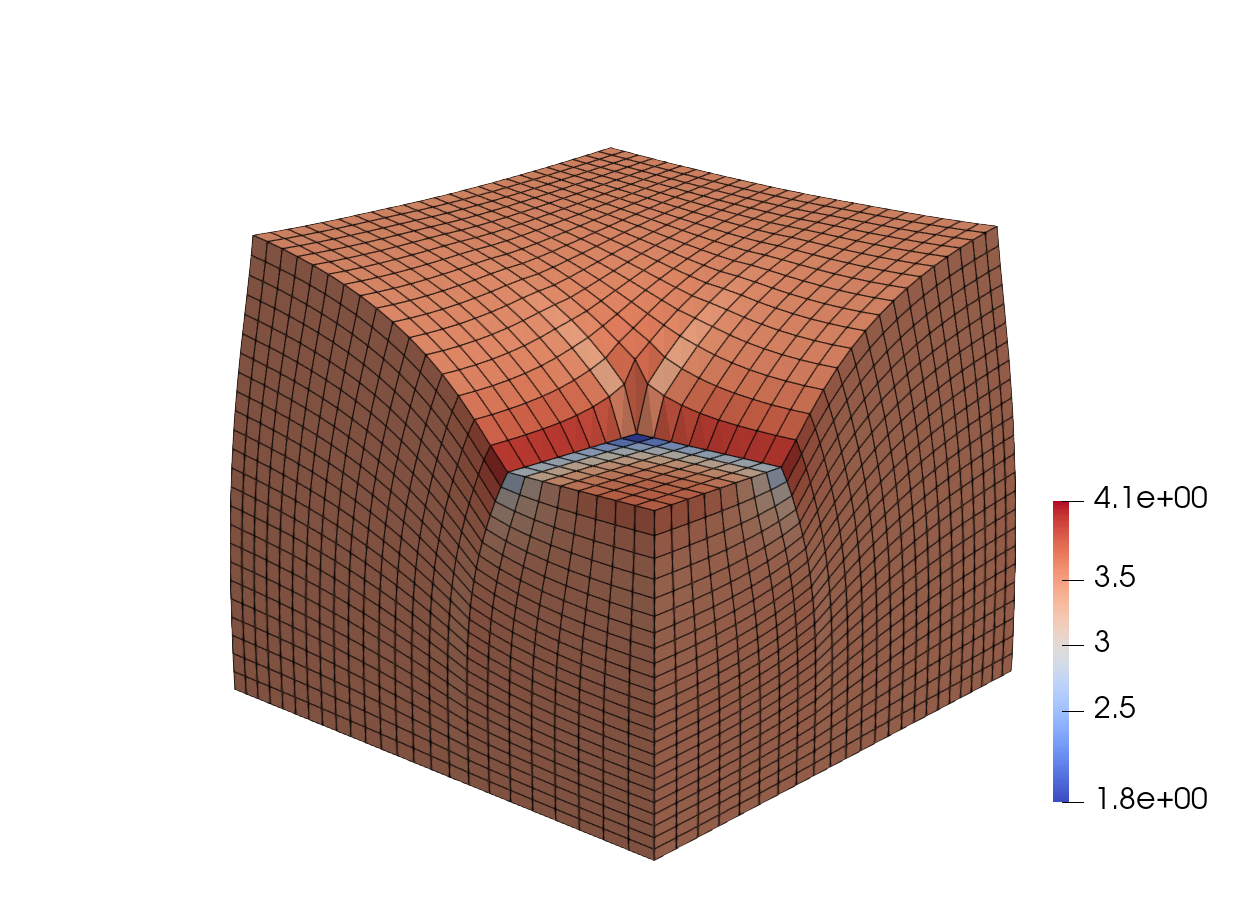}
    \end{overpic}
    \hspace{-8ex}
    \begin{overpic}[scale=0.13]
      {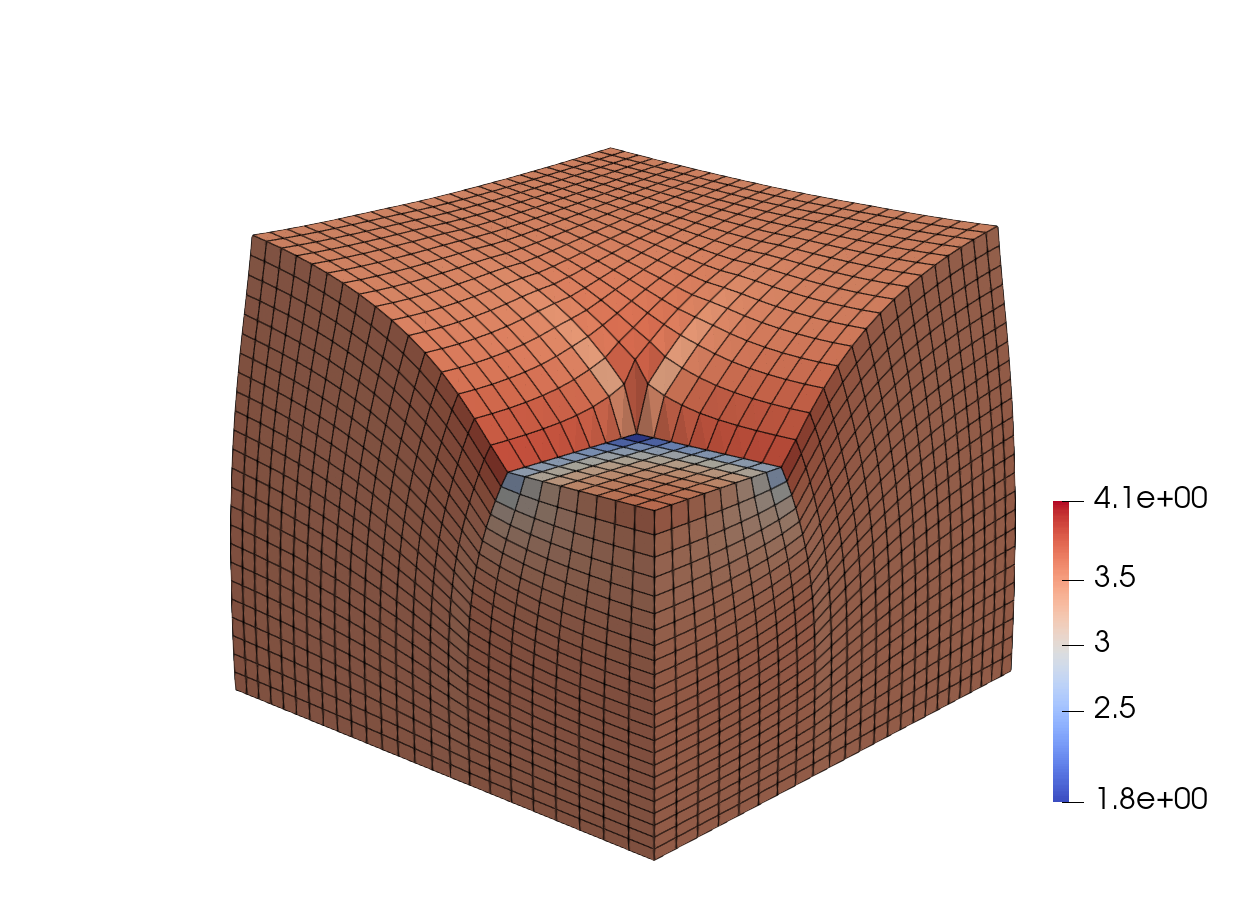}
      \put(95,20){\tiny $v$}
    \end{overpic}
    \begin{overpic}[scale=0.13]
      {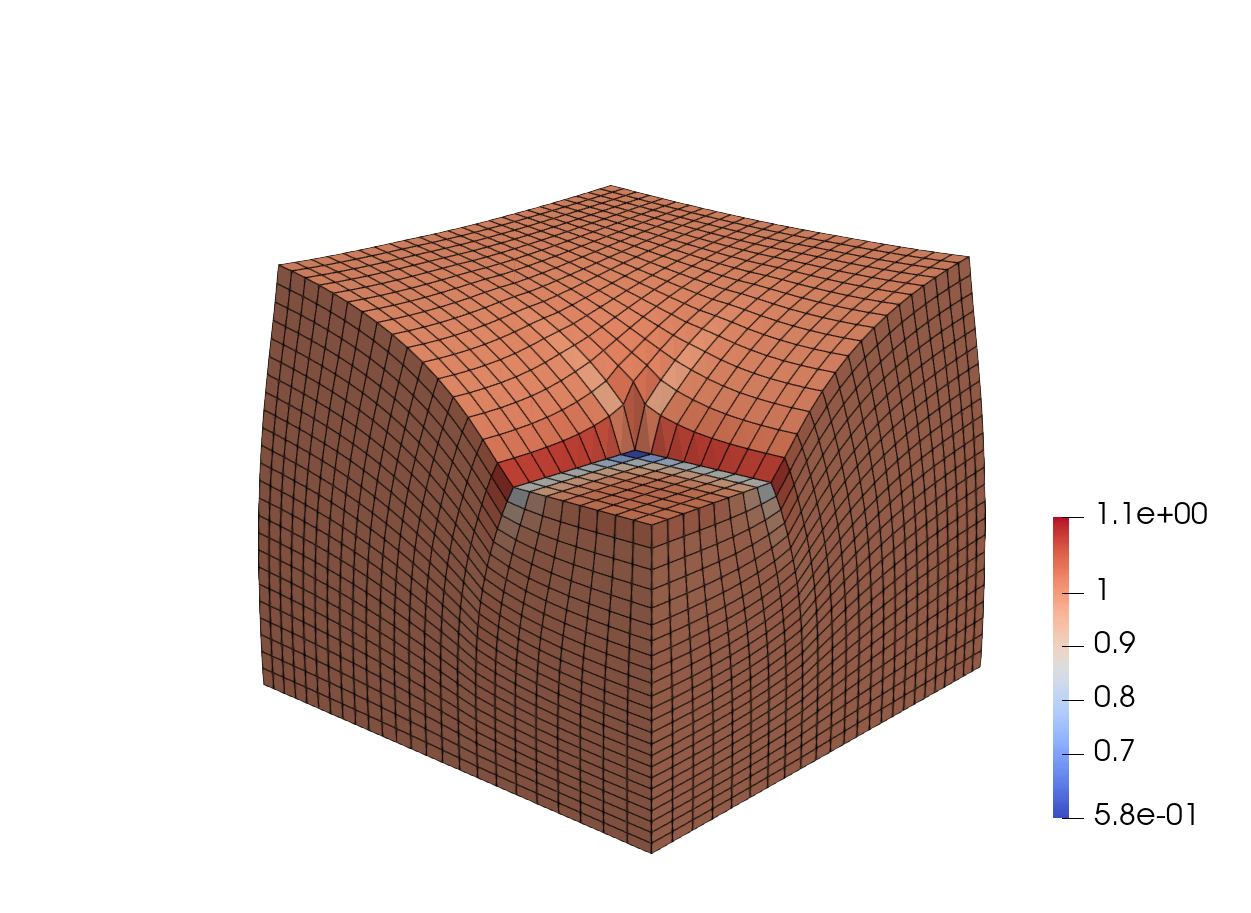}
      \put(5,5){(iii)}
    \end{overpic}
    \hspace{-8ex}
    \begin{overpic}[scale=0.13]
      {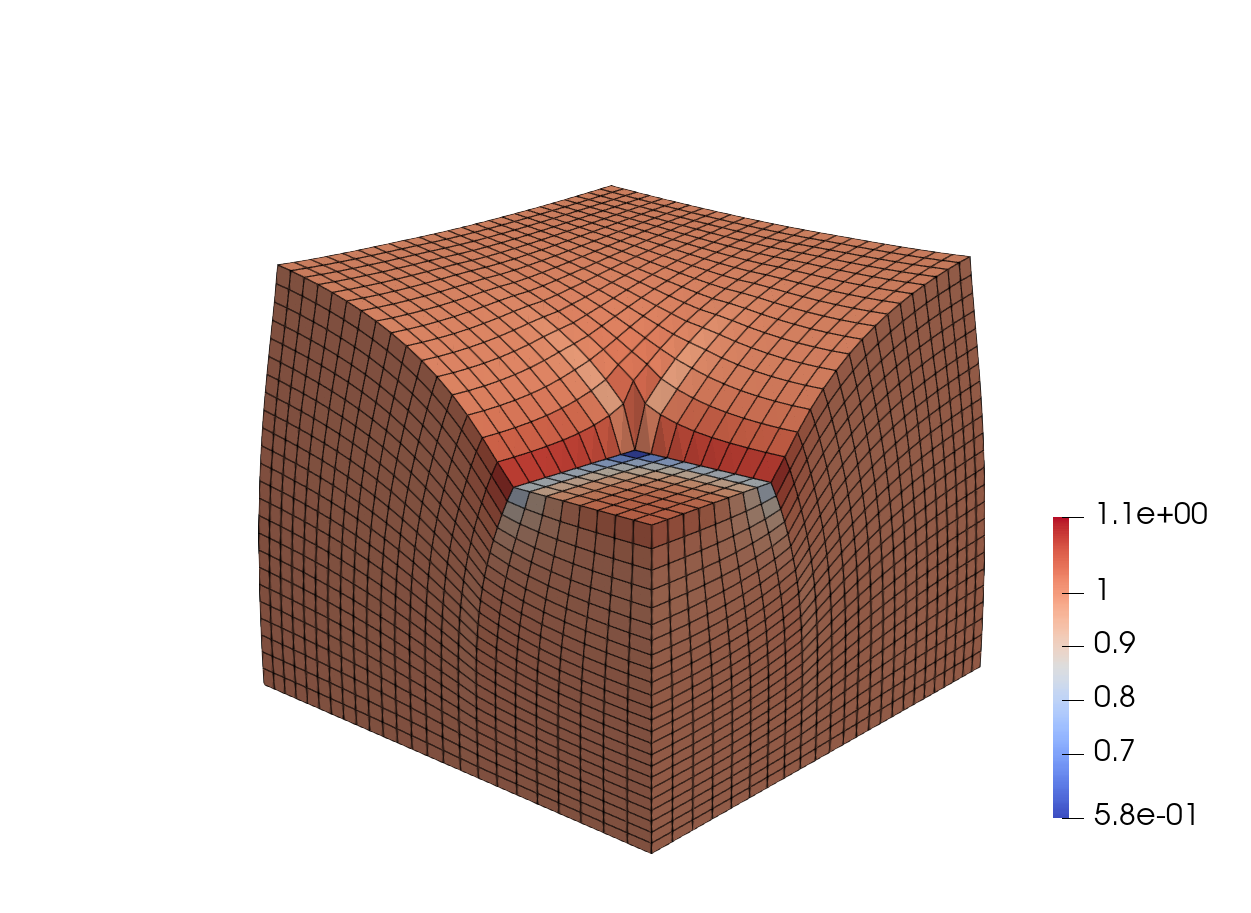}
    \end{overpic}
    \hspace{-8ex}
    \begin{overpic}[scale=0.13]
      {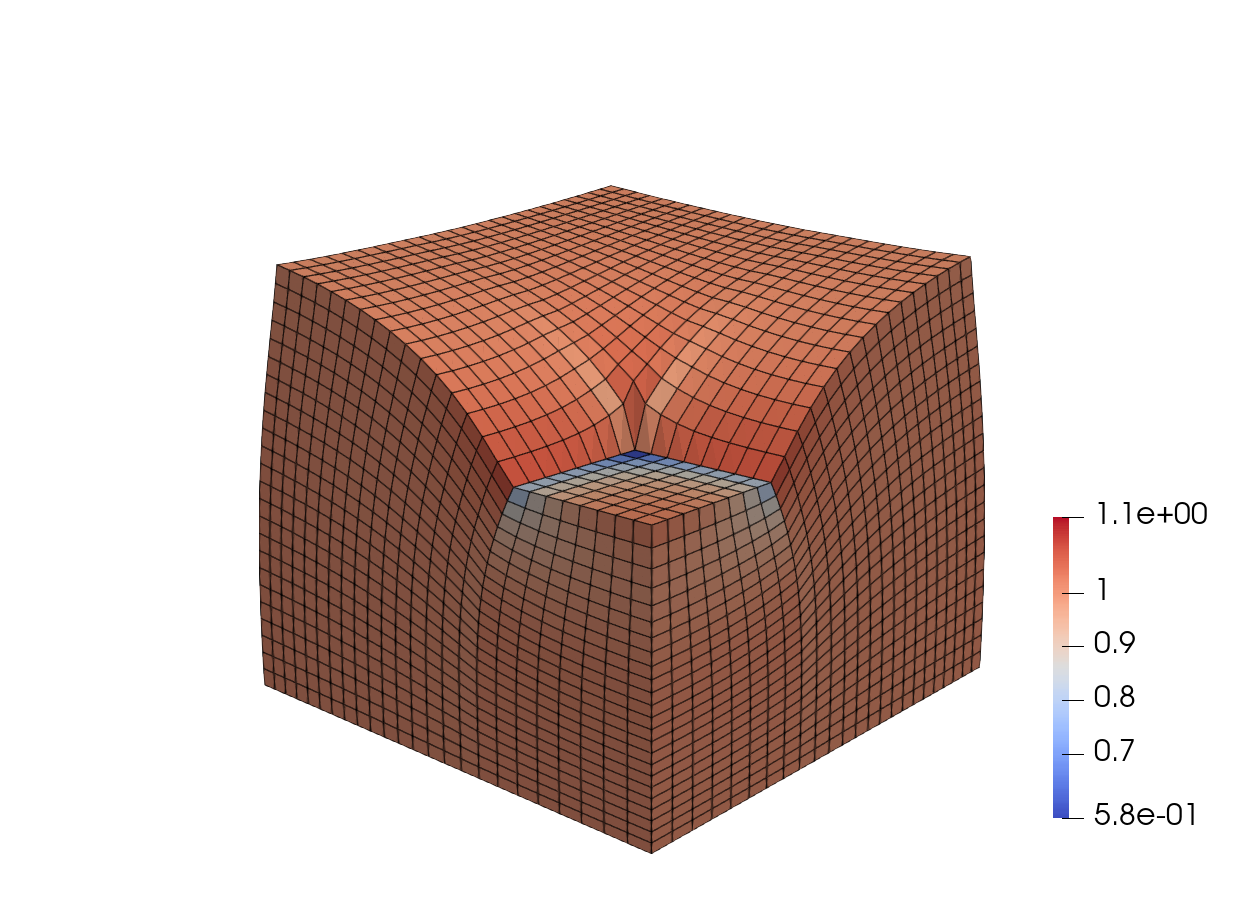}
      \put(95,20){\tiny $J$}
    \end{overpic}
    \caption{Evolution of the spatial distribution of (i) the chemical potential $\mu$, (ii) the swelling volume fraction $v$ and (iii) Jacobian $J$ obtained with the $Q_{1} RT_{0}$ discretization at $t=\unit[1]{s}$, $t=\unit[2]{s}$ and $t=\unit[6]{s}$ (from left to right)
    }
    \label{fig:evol-spat-distr-field-mech-induced-diffusion-lagr-raviart-thomas}
  \end{figure*}
  
  Inspecting the spatial distributions of the chemical potential, the swelling volume fraction and the Jacobian in Figure~\ref{fig:evol-spat-distr-field-mech-induced-diffusion-lagr-raviart-thomas} obtained with the $Q_{1} RT_{0}$ ansatz function, it becomes apparent that the mechanical deformation leads to a significant redistribution of the fluid inside the body.
  In particular, it can be seen that after the initial stage of the loading history ($0\le t \le \unit[1]{s}$), the chemical potential just below the flat punch has increased considerably due to the rather strict enforcement of the %
  {\OR penalty} constraint {\OR by} %
  the choice of material parameters $\frac{\uplambda}{\upgamma}=100$.
  Given the definition of the chemical potential, which specializes to
  \begin{multline}
    \label{eq:chemical-potential-neo-hookean}
    \mu\coloneq\ps{v}\widehat{\psi} = -\frac{\uplambda}{J_{0}}\left[JJ_{0}-1-v\right] + \\ \frac{\alpha}{J_{0}}\left[\ln\left(\frac{v}{1+v}\right)+\frac{1}{1+v}+\frac{\chi}{\left[1+v\right]^{2}}\right]
  \end{multline}
  for the free-energy function given in~\eqref{eq:free-energy-function-isotropic-pre-swollen-state}, the contribution associated with the constraint can readily {\OR be} identified as the first term on the right hand side of~\eqref{eq:chemical-potential-neo-hookean}.

  During the subsequent holding stage of the loading history, i.e., the displacement coefficient $\hat{u}$ is constant, a relaxation of the body can be observed, which results in a balanced chemical potential field alongside with a strong reduction of the swelling volume fraction below the flat punch.
  The spatial distribution of the Jacobian depicted in Figure~\ref{fig:evol-spat-distr-field-mech-induced-diffusion-lagr-raviart-thomas} (iii) is closely tied to the distribution of the swelling volume fraction.
}

{\SP
  In principle, similar observations are made in the simulation of the deformation induced diffusion problem, which employs the $Q_{1} Q_{1}$ ansatz functions.
  However, significant differences occur during the holding stage of the loading history, in which deformations are due to the diffusion of the fluid.
  In particular, a checker board pattern develops below the flat punch, which is clearly visible in all three fields depicted in Fig.~\ref{fig:evol-spat-distr-field-mech-induced-diffusion-lagr-lagr} at the end of the simulation at $t=\unit[6]{s}$.
  {\OR This is a result of the $Q_{1}$ discretization; see the brief discussion in Section~\ref{sec:discr-fluid-flux}}.
  
  {\OR For the $Q_{1}$ ansatz for the fluid flux $\bm{J}_{v}$ the swelling volume fraction is not constant within a finite element.}
  Due to the rather strict enforcement of the incompressibility constraint this heterogeneity is further amplified.
  Of course, a selective reduced integration technique is able to cure this problem, as shown in~\cite{BogNatMie:2017:min} for the two-dimensional case.
  Herein, \eqref{eq:time-discrete-balance-solute-volume} is solved at a single integration point per element and the current value $v$ are subsequently transferred to the remaining Gauss points.
  {\OR Other choices, such as a three-field formulation, are also possible.}
  The use of the $RT_{0}$ ansatz function, however, may be more appropriate as it yields, both, a conforming discretization and a lower number of degree of freedom per element compared to the standard Lagrange approximation.
  
  Note, however, that in the subsequent parallel simulations, the penalty parameter is smaller by more than an order of magnitude such that the problems described in this section were not observed.
 \begin{figure*}[h]
   \centering
   \begin{overpic}[scale=0.13]
     {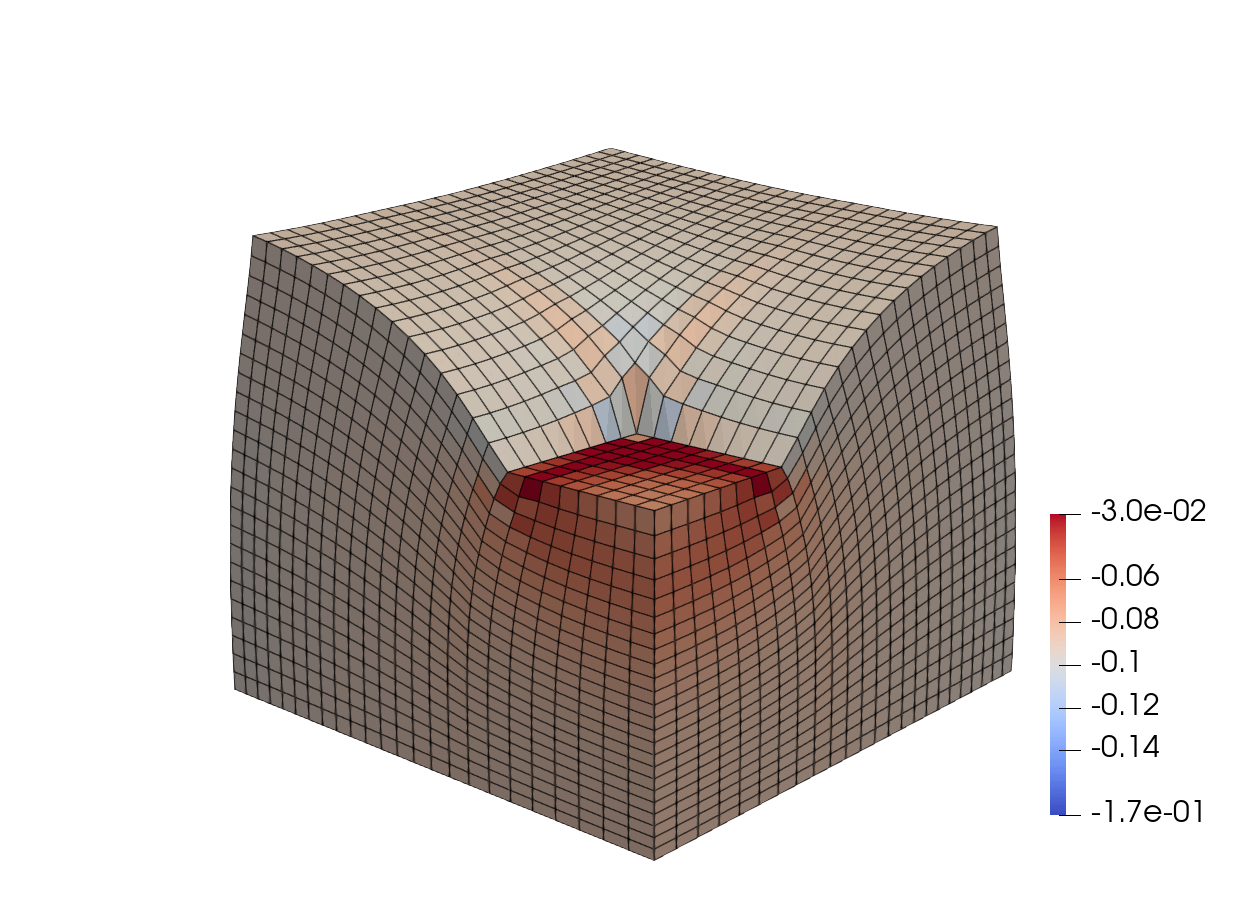}
     \put(5,5){(i)}
   \end{overpic}
   \hspace{-8ex}
   \begin{overpic}[scale=0.13]
     {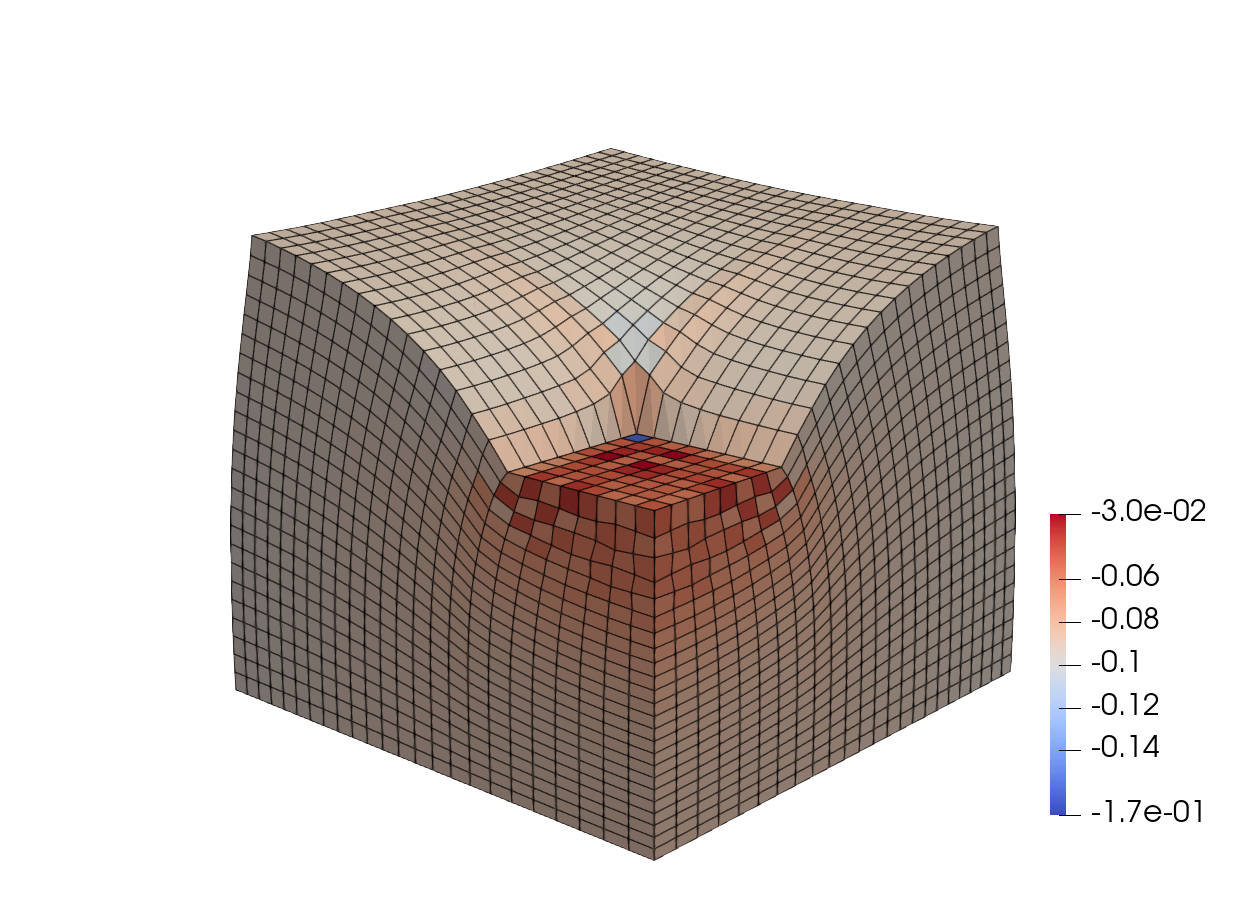}
   \end{overpic}
   \hspace{-8ex}
   \begin{overpic}[scale=0.13]
     {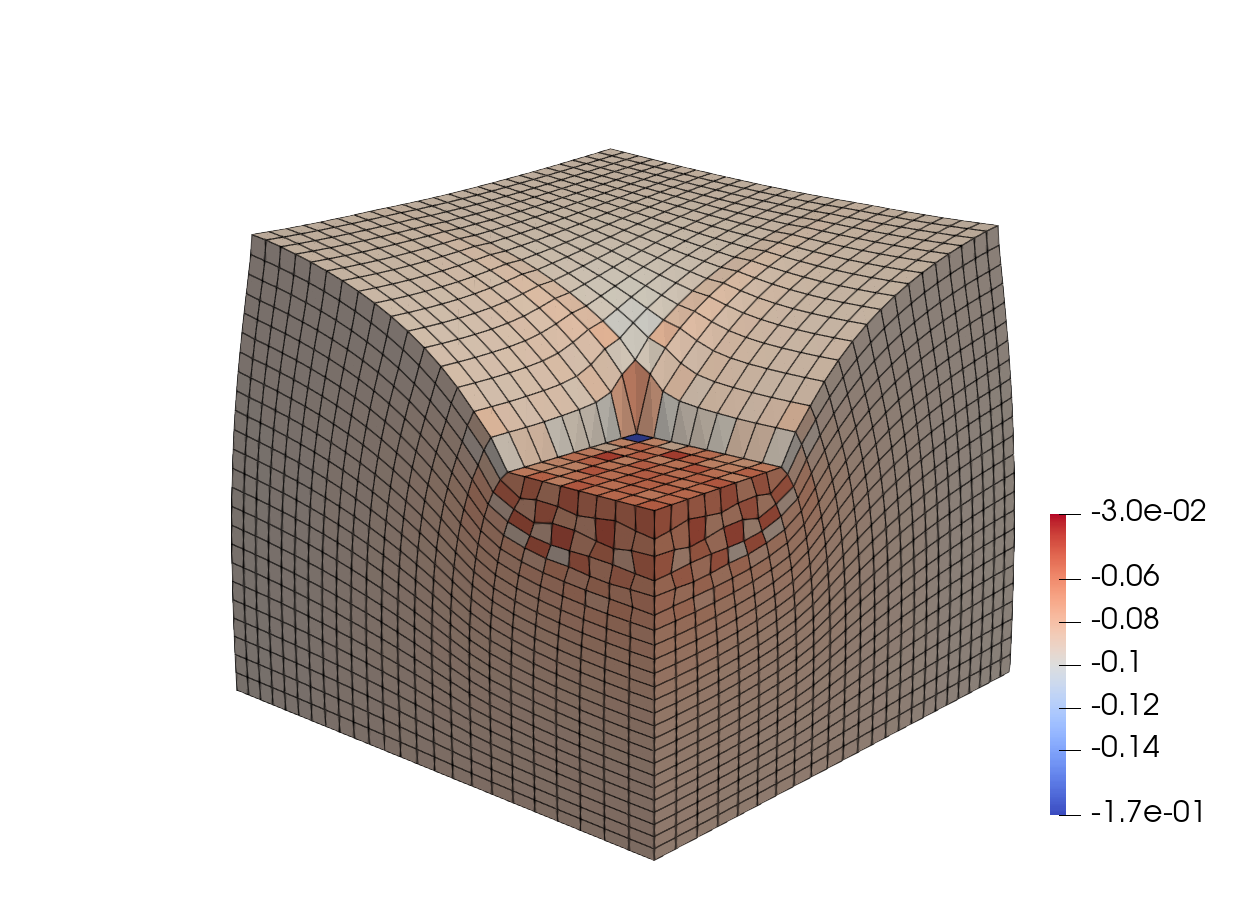}
     \put(95,20){\tiny $\mu$}
   \end{overpic}
   \begin{overpic}[scale=0.13]
     {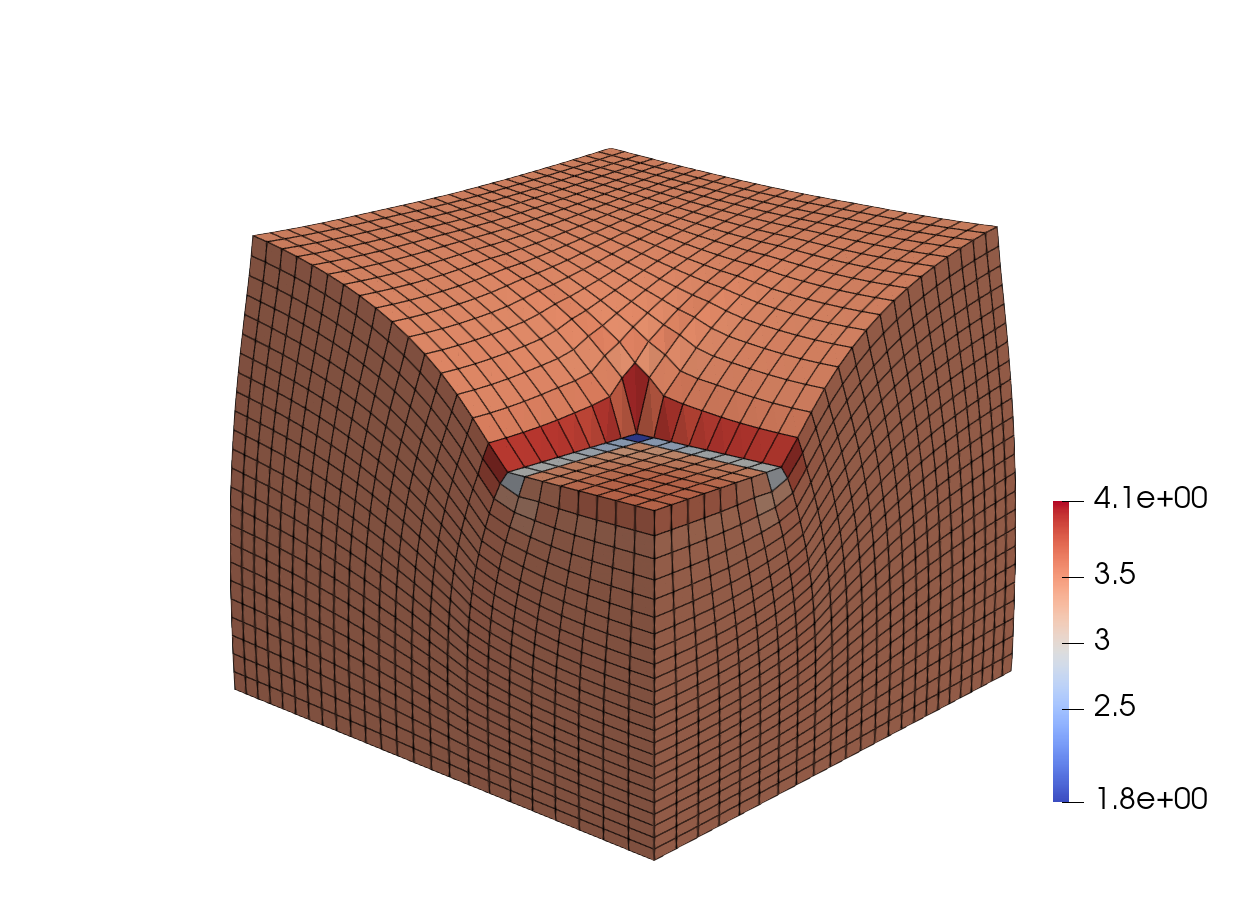}
     \put(5,5){(ii)}
   \end{overpic}
   \hspace{-8ex}
   \begin{overpic}[scale=0.13]
     {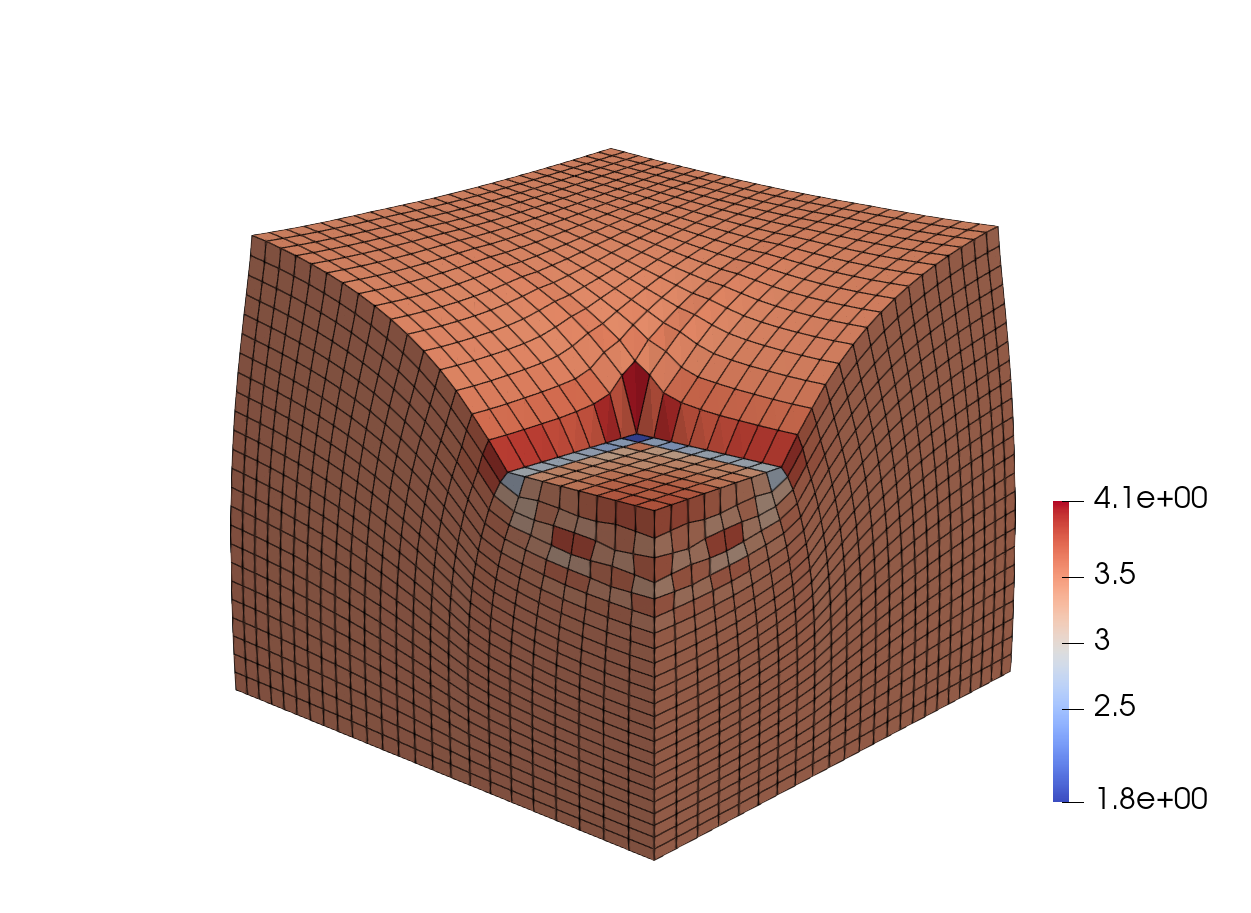}
   \end{overpic}
   \hspace{-8ex}
   \begin{overpic}[scale=0.13]
     {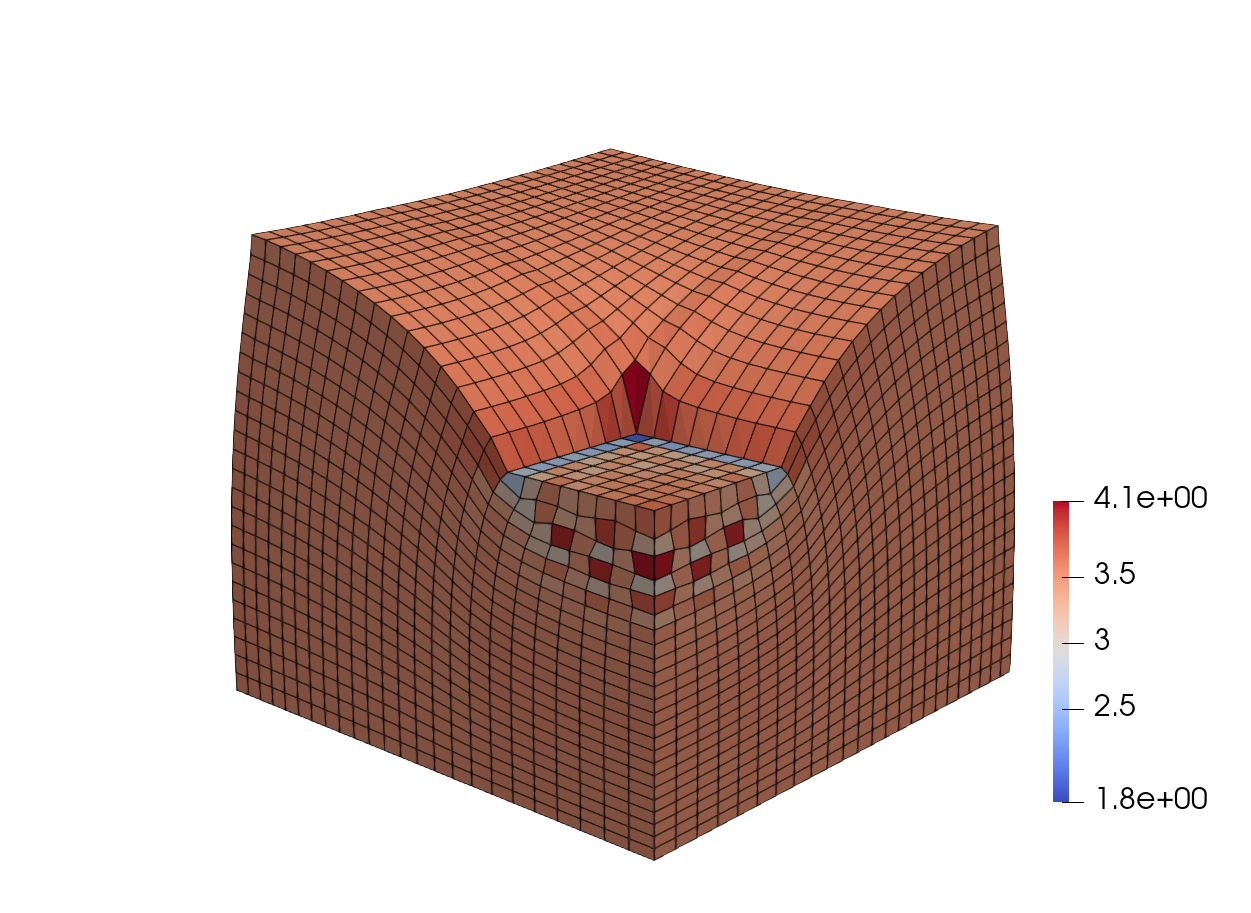}
     \put(95,20){\tiny $v$}
   \end{overpic}
   \begin{overpic}[scale=0.13]
     {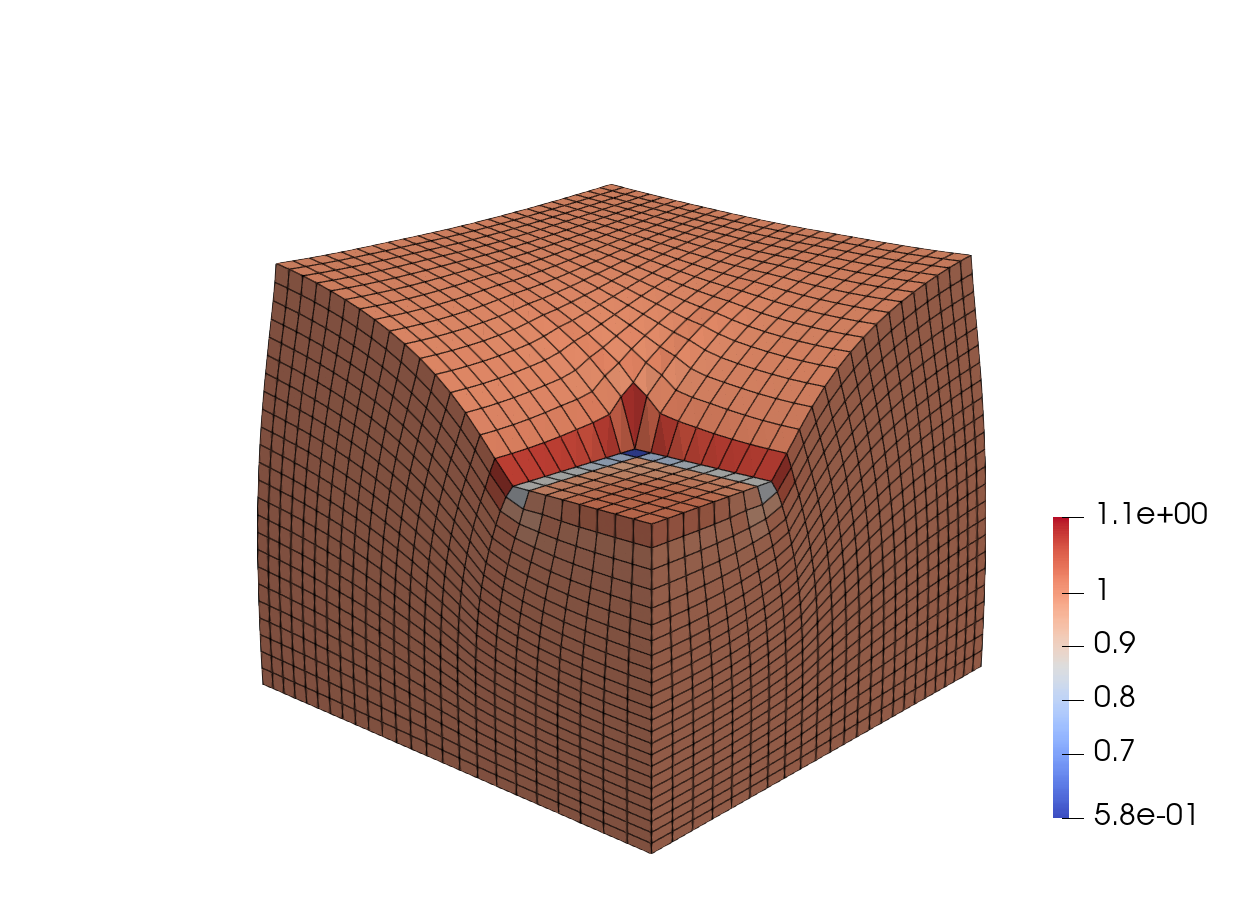}
     \put(5,5){(iii)}
   \end{overpic}
   \hspace{-8ex}
   \begin{overpic}[scale=0.13]
     {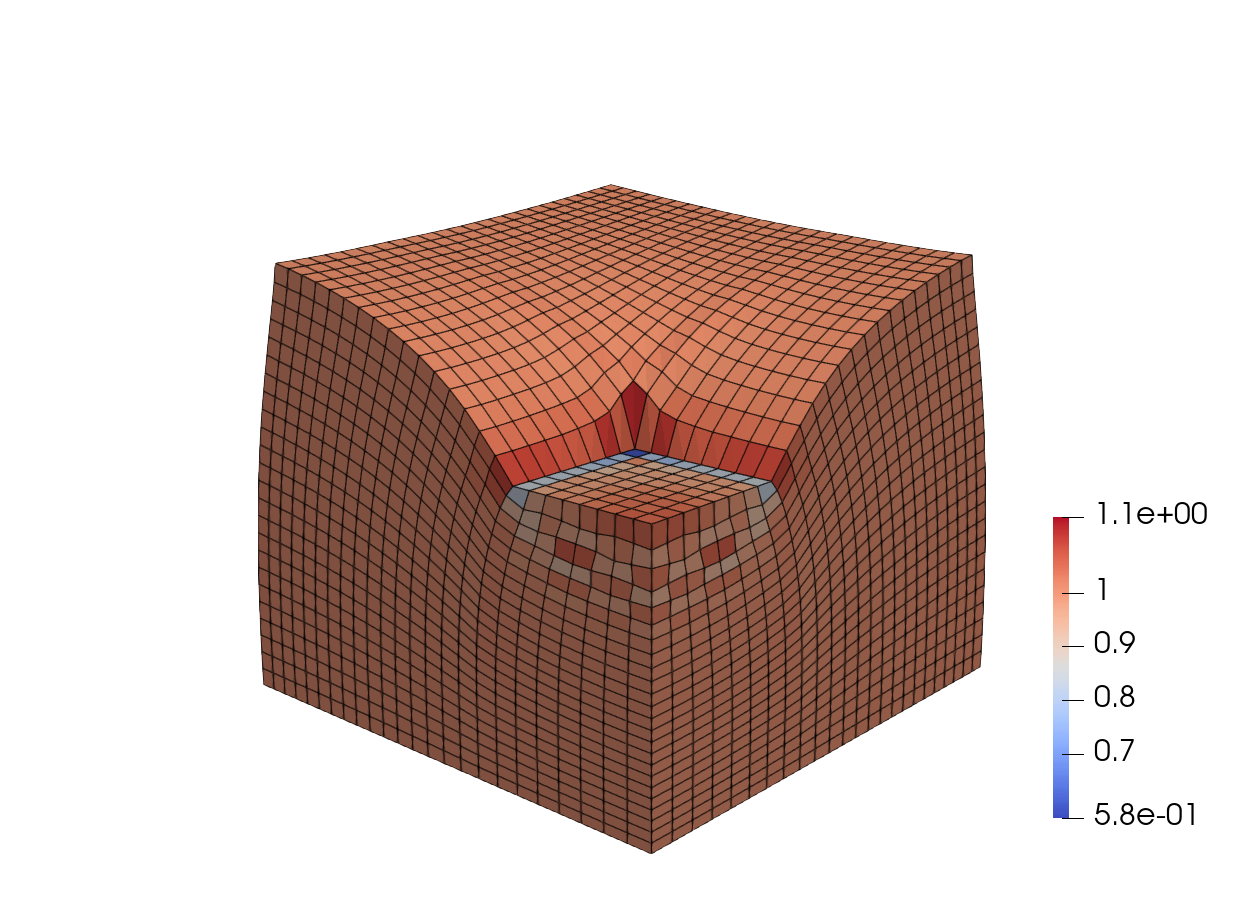}
   \end{overpic}
   \hspace{-8ex}
   \begin{overpic}[scale=0.13]
     {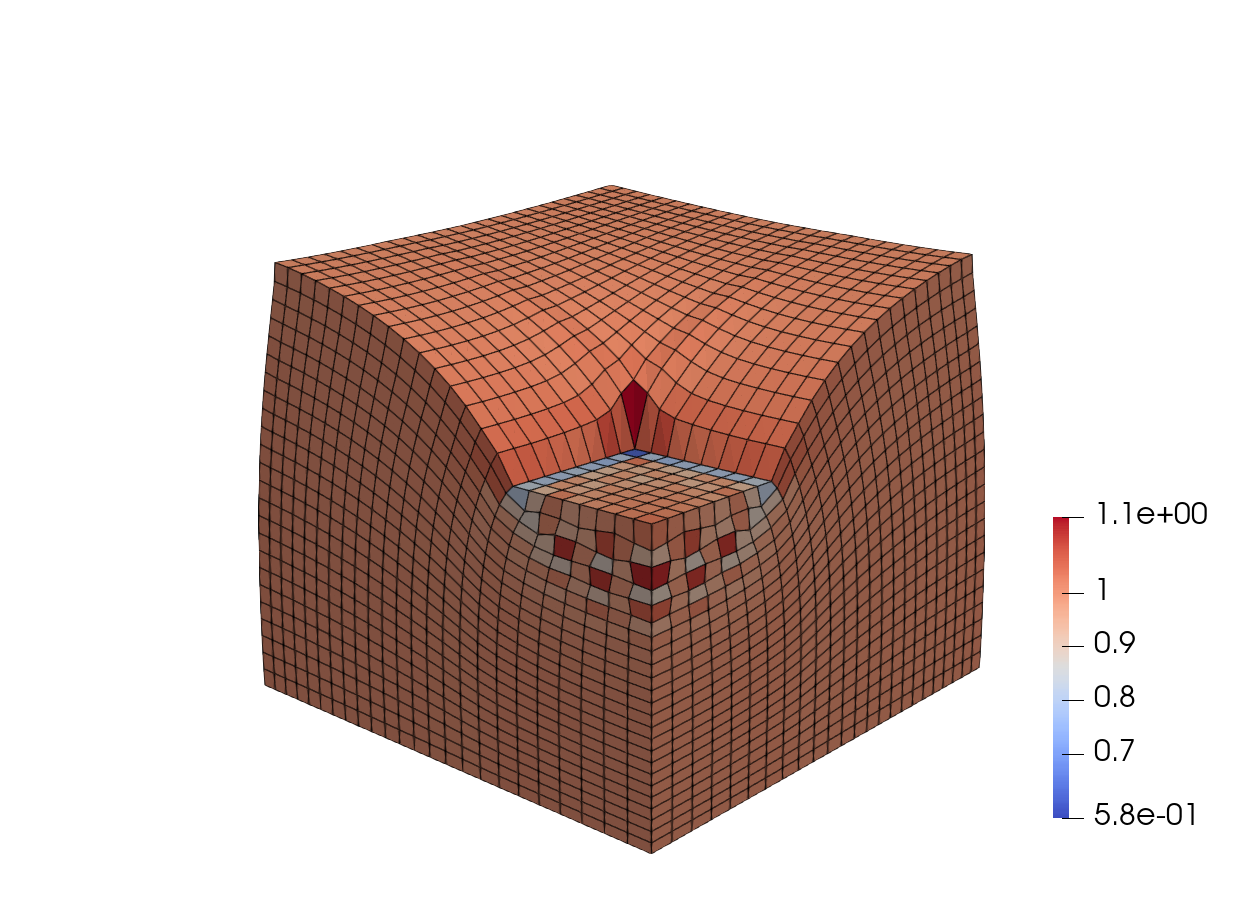}
     \put(95,20){\tiny $J$}
   \end{overpic}
   \caption{Evolution of the spatial distribution of (i) the chemical potential $\mu$, (ii) the swelling volume fraction $v$ and (iii) Jacobian $J$ obtained with the $Q_{1} Q_{1}$ discretization at $t=\unit[1]{s}$, $t=\unit[2]{s}$ and $t=\unit[6]{s}$ (from left to right)}
   \label{fig:evol-spat-distr-field-mech-induced-diffusion-lagr-lagr}
 \end{figure*}
}

\section{Numerical Results}
\label{sec:numer}
\subsection{Performance of the iterative solver}
\label{sec:perf}

{\FR To evaluate {\OR the numerical and parallel performance of the FROSch framework applied to the monolithic system in fully algebraic mode} we consider the boundary value problems described in Section~\ref{sec:free-swell-bvp} and Section~\ref{sec:mech-induced-diffusion-bvp}. We refer to the one in Section~\ref{sec:free-swell-bvp} as the \textit{free-swelling problem} and denote the problem specified in Section~\ref{sec:mech-induced-diffusion-bvp} as the \textit{mechanically induced problem}.}

}

{\FR We compare the use of $Q_1 RT_0$ and $Q_1 Q_1$ ansatz functions
 regarding 
 {\OR the consequences for}
 {\OR the}
 numerical and parallel performance {\OR of the simulation}.

 The different ansatz functions result in different numbers of degrees of freedom per node. For the {\FR $Q_1 Q_1$ ansatz} each node has six degrees of freedom.%
 The usage of {\FR $Q_{1} RT_0$} elements leads to three degree of freedoms per node and one per element face.
 {\OR %
 If not noted otherwise, the construction of the coarse spaces uses the nullspace of the Laplace operator.}
}
{\FR The computing times are always %
{\OR sums}
over all time steps and  Newton steps. We denote the time to assemble the tangent matrix in each Newton step by \textit{Assemble Matrix Time}. 

{\FR By \textit{Solver Time} we denote the time to build the preconditioner (\textit{Setup Time}) and to perform the Krylov iterations (\textit{Krylov Time}).}

{\FR For the triangulation, we executed four refinement cycles on an initial mesh with 27 finite elements resulting in a structured mesh of 110\,592 finite elements. 
}

\subsection{Strong parallel scalability}
\label{sec:strong-scaling}
{\OR For the strong scalability, we consider our problem on the structured mesh of 110\,592 cells, which results in 691\,635 degrees of freedom for the $Q_1 RT_0$ elements and 705\,894 degrees of freedom for the $Q_1 Q_1$ discretization.}
We then increase the number of subdomains {\OR and} cores and expect the computing time to decrease. %
}

\subsubsection{Linear elasticity benchmark problem}
\label{sec:linElas}
To provide a baseline to compare with, we first briefly present strong scaling results for a linear elastic benchmark problem on the unit cube $(0,1)^3$, using homogeneous Dirichlet boundary conditions on the complete boundary and discretized using $Q_1$ elements; see Table~\ref{tab:stronglinelas}.
Here, the 110\,592 finite elements result in only 352\,947 degrees of freedom
since %
the diffusion problem is missing.
{We use a generic right-hand-side vector of ones $(1,\ldots,1)^T$.}

We will use this simple problem as a baseline to evaluate the performance of our solver for our nonlinear coupled problems. Note that due to the homogeneous boundary conditions on the complete boundary, this problem is quite well conditioned, and a low number of Krylov iterations should be expected.

In Table~\ref{tab:stronglinelas}, we see that, using the GDSW coarse space for elasticity (with three displacements but without rotations), we observe numerical scalability, i.e., the number of Krylov iterations does not increase and stays below 30.
Note that this coarse space is algebraic, however, it exploits the knowledge on the numbering of the degrees of freedom.
Other than for FETI-DP and BDDC methods, the GDSW theory does not guarantee numerical scalability for this coarse space missing the three rotations, however, experimentally, numerical scalability has been observed previously for certain, simple linear elasticity problems~\cite{Heinlein:2016:PIT,kupsAlg}.

In Table~\ref{tab:stronglinelas}, the strong scalability is good when scaling from 64 (\unit[28.23]{s}) %
to 216 cores (\unit[7.43]{s}). %
 The {\it Solver Time} increases %
 when scaling from 216 to 512 cores (\unit[41.66]{s}), indicating that coarse problem of size \unit[26\,109] is too large to be solved efficiently by Amesos2 KLU. %
The sequential coarse problem starts to dominate the solver time. %

Note that, as we increase the number of cores and subdomains, the subdomain sizes decrease. 
In our strong parallel scalability experiments, we thus profit from the superlinear complexity of the sparse direct subdomain solvers.

We also provide results for the fully algebraic mode. %
Here, the number of Krylov iterations is slightly higher and increases slowly as the number of cores increase. 
This is not surprising since in fully algebraic mode, we assume a one-dimensional nullspace, which is suitable for Laplace problems.

However, the {\it Solver Time} is comparable for both coarse spaces for 64 (\unit[28.23]{s} vs. \unit[27.35]{s}),  %
125 (\unit[16.25]{s} vs. \unit[13.35]{s}), %
and 216 (\unit[7.43]{s} vs. \unit[5.69]{s}) %
 cores. Notably, the {\it Solver Time} is better for the fully algebraic mode for 512 cores (\unit[41.66]{s} vs \unit[6.59]{s}) %
as a result of the smaller coarse space in the fully algebraic mode (\unit[8\,594] vs. \unit[26\,109]).

In Table~\ref{tab:detaillinElas}, more details of the computational cost are presented. These timings show that for 64, 125, and 216 cores the cost is dominated by the factorizations of the subdomain matrices $K_i$.
Only for 512 cores this is not the case any more.

Interestingly, the fully algebraic mode is thus preferable within the range of processor cores discussed here, although numerical scalability is not achieved.
\begin{table*}
\begin{tabular}{l || c c c |c c c }
\hline 
&\multicolumn{6}{|c}{\bf GDSW}\\
& \multicolumn{3}{c|}{\bf coarse space for elasticity w/o rot. }& \multicolumn{3}{c}{\bf Fully algebraic mode} \\
\# cores & Krylov & Assemble &{\bf Solver} & Krylov & Assemble &{\bf Solver}  \\
         &        & Time     &{\bf Time}&        & Time     &{\bf Time}\\ \hline
64  &32 & \unit[2.90]{s} & {\bf \unit[28.23]{s} }& 39&\unit[2.79]{s}& {\bf \unit[27.35]{s} }\\
125 &32& \unit[1.50]{s} & {\bf \unit[16.25]{s} }&42&\unit[1.48]{s}& {\bf \unit[13.35]{s} }\\
216 &25& \unit[0.85]{s} & {\bf \unit[7.43]{s} }&37&\unit[0.91]{s}& {\bf  \unit[5.68]{s} } \\
512 &30& \unit[0.60]{s} & {\bf \unit[41.66]{s} }&48&\unit[0.57]{s}& {\bf \unit[6.59]{s} } \\ \hline
\end{tabular}
\caption{Strong scaling for the linear elasticity model problem in 3 dimensions using $Q_1$ elements. Dirichlet boundary conditions on the complete boundary.
We operate on a structures mesh with 110\,592 finite elements such that we have 352\,947 degrees of freedom. We use an overlap of two elements. We use the standard GDSW, without rotations in the coarse space. 
}
\label{tab:stronglinelas}
\end{table*}

\begin{table*}
\centering
\renewcommand{\arraystretch}{1.2}
\scriptsize{
\begin{tabular}{l || c c c c}
\hline         
 \multicolumn{5}{c|}{{\bf GDSW Coarse Space}}  \\ \hline
\multicolumn{3}{c}{} &{\bf coarse space for elasticity w/o rot. } & {\bf  Fully algebraic mode}  \\ \hline
\# Cores &  Avg. Size & Max. Size & Comp. $K_i$ &Comp. $K_i$ \\
         &  $K_i$     &  $K_i$    & Time  & Time \\ \hline
     64  & 15\,027    &  18\,789  & 19.69s & 20.06s \\
     125 & 10\,021	   &  13\,824  & 8.53s  & 8.54s \\
     216 & 5\,716.3   &  7\,581   &3.46s   & 3.53s\\
     512 & 4\,672.3	   &  7\,038   &2.19s   & 2.02s \\ \hline
\end{tabular}}
\caption{Detailed cost of the overlapping subdomain problems $K_i$ contained in the {\it Setup Time} in Table~\ref{tab:stronglinelas}.
}
\label{tab:detaillinElas}
\end{table*}

\subsubsection{Free swelling problem}
\label{subsec:free-strong}
{\FR We now discuss the strong scalabilty results for the \textit{free-swelling problem}; see Section~\ref{sec:free-swell-bvp}. Here, the pre-swollen Jacobian $J_0$ is %
chose as $J_0=1.01$. The other material paramater are chosen according to Table~\ref{tab:material-parameter} and Table~\ref{tab:param-bound-cond-free-swelling}.

{\OR For the parallel performance study, we perform two time steps for each test run. In each time step, again, 5 Newton iterations are needed for convergence.

{\OR 
For a numerically scalable preconditioner, we would expect the number of Krylov iterations to be bounded. 
In Table~\ref{tab:free-s}, we observe that we do not obtain good numerical scalability, i.e., the number of iteration increases by 50 percent when scaling from 64 to 512 cores. This can be attributed to the fully algebraic mode, whose coarse space is not quite suitable to obtain numerical scalability; see also Section~\ref{sec:linElas}.
Interestingly, the results are very similar for GDSW and RGDSW with the exception of 512 cores, where the smaller coarse space of the RGDSW method results in slightly better {\it Solver Time}. This is interesting, since the RGDSW coarse space is typically significantly smaller. 
This indicates that the RGDSW coarse space should be preferred in our future works.

The number of iterations is smaller for the $Q_1RT_0$ discretization compared to the $Q_1Q_1$ discretization. Since, in addition, the local subdomain problems are significantly larger when using $Q_1Q_1$ (see Table~\ref{tab:free-detail})
the {\it Solver Times} are better by (approximately) a factor of two when using $Q_1RT_0$ discretizations. %

Strong parallel scalability is good when scaling from 64 to 216 cores. Only incremental improvements are obtained for 512 cores indicating that the problem is too small.
}

If we relate these results to our linear elasticity benchmark problem in Section~\ref{sec:linElas}, we see that with respect to the number of iterations, the (average) number of Krylov iterations is higher by a factor 
{1.5 to 2} 
for the coupled problem compared to the linear elastic benchmark. We believe that this is an  acceptable result.

If we compare the {\it Solver Time}, we need to multiply the {\it Solver Time} in Table~\ref{tab:stronglinelas} by a factor of 10, since 10 linearized systems are solved in the nonlinear coupled problem. 
Here, we see that the {\it Solver Time} is higher by a factor slightly more than 3 when using $Q_1RT_0$ compared to solving 10 times the linear elastic benchmark problem of Section~\ref{sec:linElas}. For $Q_1Q_1$, this factor is closer to 6 or 7. Interestingly, in both cases, this is mostly a result of larger factorization times for the local subdomain matrices (see Table~\ref{tab:free-detail}) and only to a small extent a result of the larger number of Krylov iterations.

\begin{figure}[h]
  \centering
  \includegraphics[width=0.35\textwidth]{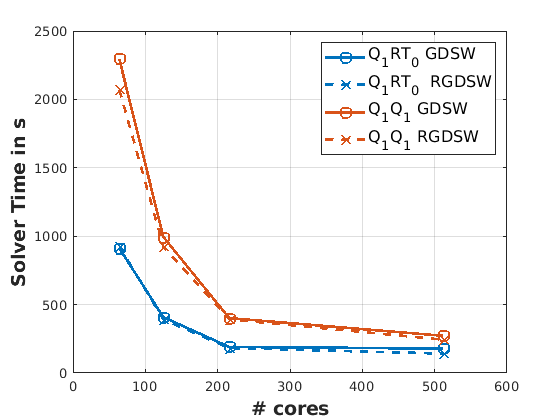}
  \caption{\FR Strong scalability of the \textit{Solver Time (Setup Time + Krylov Time)} for the \textit{free-swelling problem}; see Table~\ref{tab:free-s} for the data.}
  \label{fig:free-strong-1node}
\end{figure}
\begin{table*}
\centering
\renewcommand{\arraystretch}{1.2}
\scriptsize{
\begin{tabular}{l || c c c c |c| c c c c |c}
\hline
\multicolumn{11}{c}{Free swelling problem}\\\hline
\multicolumn{1}{c||}{}&\multicolumn{10}{c}{\bf{FROSch with GDSW coarse space}} \\ \hline
  & \multicolumn{5}{c|}{{\FR $Q_1 RT_0$}} & \multicolumn{5}{c}{{\FR $Q_1 Q_1$}}\\
 & \multicolumn{5}{c|}{5 Newton steps/time step, 2 time steps}& \multicolumn{5}{c}{5 Newton steps/time step, 2 time steps} \\
\#cores & Avg.  & Assemble & Setup & Krylov &{\bf Solver}  & Avg.  & Assemble & {\FR Setup}& Krylov &{\bf Solver}\\
       &  Krylov            & Time    & Time & Time & {\bf Time}&Krylov & Time &Time &Time & {\bf Time}\\ \hline
64 &  49.80 &  \unit[66.02]{s} &  \unit[850.26]{s} &  \unit[62.36]{s} &  \unit[912.62]{s}  &70.90&	\unit[111.81]{s}&  \unit[2\,106.6]{s}& \unit[189.91]{s}& {\bf \unit[2\,296.5]{s}} \\
125&  56.40&  \unit[34.29]{s} & \unit[367.25]{s}  & \unit[40.21]{s} & \unit[407.47]{s} &79.90&	\unit[58.30]{s}&	\unit[900.43]{s}&\unit[93.48]{s} & {\bf \unit[993.91]{s}}\\
216& 46.10&   \unit[19.42]{s} &  \unit[172.10]{s}&   \unit[18.44]{s}&  \unit[190.53]{s}& 62.70& \unit[32.35]{s}&	\unit[362.95]{s}& \unit[38.42]{s}& {\bf \unit[401.37]{s}}\\
512&    73.70&    \unit[9.67]{s} &  \unit[139.44]{s} &   \unit[38.69]{s} &  \unit[178.12]{s} 
&103.40& \unit[15.50]{s}&	\unit[213.45]{s}& \unit[58.27]{s}& {\bf \unit[271.72]{s}} \\ \hline
\multicolumn{1}{c||}{}&\multicolumn{10}{c}{\bf{FROSch with RGDSW coarse space}} \\ \hline
 & \multicolumn{5}{c|}{{\FR $Q_1 RT_0$}} & \multicolumn{5}{c}{{\FR $Q_1 Q_1$}}\\
 & \multicolumn{5}{c|}{5 Newton steps/time step, 2 time steps}& \multicolumn{5}{c}{5 Newton steps/time step, 2 time steps} \\
\#cores & Avg.  & Assemble & Setup & Krylov &{\bf Solver}  & Avg.  & Assemble & {\FR Setup}& Krylov &{\bf Solver}\\
       &  Krylov            & Time    & Time & Time & {\bf Time}&Krylov & Time &Time &Time & {\bf Time}\\ \hline
64 & 48.90&   \unit[65.63]{s} &  \unit[866.80]{s} & \unit[60.80]{s} & \unit[927.60]{s} &68.20&\unit[107.33]{s}&	\unit[1\,944.65]{s}&	\unit[127.08]{s} & {\bf \unit[2071.73]{s}}\\
125& 56.00&   \unit[34.62]{s} & \unit[352.61]{s} & \unit[37.57]{s} &  \unit[390.18]{s}  &78.40&	\unit[57.76]{s}&	\unit[839.54]{s}&	\unit[84.32]{s} & {\bf \unit[923.86]{s}}\\
216& 46.40&   \unit[20.02]{s} &  \unit[165.12]{s} &   \unit[17.38]{s} &  \unit[182.50]{s} &61.00&	\unit[32.65]{s}&	\unit[356.77]{s}&	\unit[35.45]{s}& {\bf \unit[392.22]{s}}\\
512& 72.60&   \unit[9.42]{s} &  \unit[119.23]{s} &   \unit[22.32]{s} &  \unit[141.55]{s} &102.20&	\unit[16.02]{s}&	\unit[201.37]{s}& \unit[42.01]{s} & {\bf \unit[243.38]{s}} \\ \hline
\end{tabular}}
\caption{\FR Strong scalabilty results for the \textit{free-swelling problem} corresponding to Figure~\ref{fig:free-strong-1node}. We operate on a triangulation with 110\,592 finite elements resulting in 691\,635 degrees of freedom for the $Q_1 RT_0$ ansatz functions and 705\,894 degrees of freedom for the $Q_1 Q_1$ ansatz functions. We choose the boundary conditions (i) described in Section~\ref{sec:free-swell-bvp}. We apply the FROSch framework in algebraic mode with the GDSW and the RGDSW coarse space. We peform two time steps with $\Delta t=\unit[0.05]{s}$. By {\it Avg. Krylov} we denote the average number of Krylov iteration over all time- and Newton steps. The time measurements are taken over the whole computation.}
\label{tab:free-s}
\end{table*}
\begin{table*}
\centering
\renewcommand{\arraystretch}{1.2}
\scriptsize{
\begin{tabular}{l || c c c c|c c c c}
\hline         
&  \multicolumn{4}{c|}{{ $Q_1 RT_0$}} & \multicolumn{4}{c}{{$Q_1 Q_1$}} \\ \hline

\multicolumn{3}{c}{} &{\bf GDSW} & {\bf RGDSW} &\multicolumn{2}{c}{} &{\bf GDSW} & {\bf RGDSW} \\ \hline
\# Cores &  Avg. Size & Max. Size & Comp. $K_i$ &Comp. $K_i$ &  Avg. Size & Max. Size & Comp. $K_i$& Comp. $K_i$\\
         &  $K_i$     &  $K_i$    & Time & Time  &  $K_i$     &  $K_i$    & Time & Time \\ \hline
   64	 &    29316	  &	 35634	&\unit[727.26]{s}& \unit[738.45]{s}	 & 31404   &	38334	& \unit[1\,832.0]{s} &\unit[1707.0]{s}\\
   125   &  19461	& 25598	& \unit[303.15]{s} &	 \unit[285.34]{s} &	21044 & 27648 &\unit[761.37]{s} & \unit[700.11]{s}\\
   216   & 11022	& 14028 & \unit[147.47]{s} & \unit[139.18]{s} & 11958 &15162 & \unit[315.06]{s} &\unit[310.17]{s} \\
   512	 &  8937.4	& 12571 & \unit[94.02]{s} & \unit[93.82]{s} &9849.4	&13896 &\unit[153.13]{s}&\unit[149.67]{s} \\ \hline
\end{tabular}}
\caption{Detailed cost of the overlapping subdomain problems $K_i$ contained in the \textit{Setup Time} in Table~\ref{tab:free-s}.}
\label{tab:free-detail}
\end{table*}

\subsubsection{Mechanically induced diffusion problem}
\label{subsec:mech-strong}

\begin{table*}
\begin{tabular}{l||c c c c c}
\hline
Timestep & Newton & Avg. Krylov &Assemble Matrix & {\FR Setup}& Krylov \\
Size     & Steps  &             & Time           & Time   & Time \\
\hline
\multicolumn{6}{c}{{\FR $Q_1 RT_0$}} \\ \hline
0.025&4,4,4,4&132.88&\unit[34.80]{s}&\unit[286.03]{s}&\unit[89.06]{s}\\
0.05&5,5&131.70&\unit[21.37]{s}&\unit[174.48]{s}&\unit[53.44]{s}\\
0.1&5&133.60&\unit[10.43]{s}&\unit[83.46]{s}&\unit[24.47]{s}\\ \hline
\multicolumn{6}{c}{{\FR $Q_1 Q_1$}} \\
\hline
0.025 & 4,4,4,5 &123.94&59.79 s&\unit[1043.5]{s}&\unit[227.19]{s}\\
0.05& 5,5 &125.4&34.02 s&\unit[371.66]{s}&\unit[79.80]{s}\\
0.1 & 5 &123&17.27 s&\unit[190.02]{s}&\unit[39.90]{s}\\\hline
\end{tabular}
\caption{\FR Results for \textit{mechanically induced problem} using 216 cores using different time step sizes $\Delta t$ to reach $t = \unit[0.1]{s}$. We apply FROSch with the GDSW coarse space.}
\label{tab:time step}
\end{table*}

{\FR For the \textit{mechanically induced problem}, we chose a value of $J_0=4.5$ for the pre-swollen Jacobian $J_0$. The other problem parameters are chosen according to Table~\ref{tab:material-parameter}} {\SP and Table~\ref{tab:param-bound-cond-free-swelling}}.

\paragraph{Effect of the time step size}
{\OR Let us note that in our simulations the time step size $\Delta t$ has only a small influence on the convergence of the preconditioned GMRES method.} {\OR Using different choices of the time step $\Delta t$, in Table~\ref{tab:time step} we show the number of Newton and GMRES iterations.
The model problem is always solved on 216 cores until the time $t=\unit[0.1]{s}$ is reached.} %
The number of Newton iterations for each time step slighty differs; see Table~\ref{tab:time step}.
{\OR The small effect of the choice of the time step size on the Krylov iterations is explained by the lack of a mass matrix in the structure part of our model.
However, the diffusion part of the model does contain a mass matrix.
Moreover, time stepping is needed as a globalization technique of the Newton method, i.e., large time steps will result in a failure of Newton convergence.}
A different formulation including a mass matrix for the structural part of the model should be considered in the future as a possibility to improve solver convergence.

\paragraph{Strong scalability for the mechanically induced diffusion}
{\OR We, again, perform  two time steps for each test run. In each time step 5 Newton iterations are needed for convergence.
In Table~\ref{tab:strongMech}, we present results using 64 up to 512 processor cores.

First, we observe that the average number of Krylov iterations is significantly higher compared to Section~\ref{subsec:free-strong}, indicating this problem is significantly harder as a result of the boundary conditions.

Next, we observe that the average number of Krylov iterations is similar for the $Q_1 RT_0$ and the $Q_1 Q_1$ case,
which is different than in Section~\ref{subsec:free-strong}.

In both discretizations, the number average of Krylov iterations increases by about 50 percent for a larger number of cores. This is valid for the GDSW as well as for the RGDSW coase space. 

We also see that the Solver Time for $Q_1 RT_0$ is significantly better in all cases. Indeed, the time for the Krylov iteration ({\it Krylov Time}) as well as the time for the setup of the preconditioner ({\it Setup Time}) is larger for $Q_1Q_1$. The {\it Setup Times} are often drastically higher for $Q_1Q_1$.
As illustrated in Table~\ref{tab:mech-overlap}, this is a result of larger factorization times for the sparse matrices arising from $Q_1Q_1$ discretizations;

\begin{figure*}%
\centering
\includegraphics[width=0.5\textwidth]{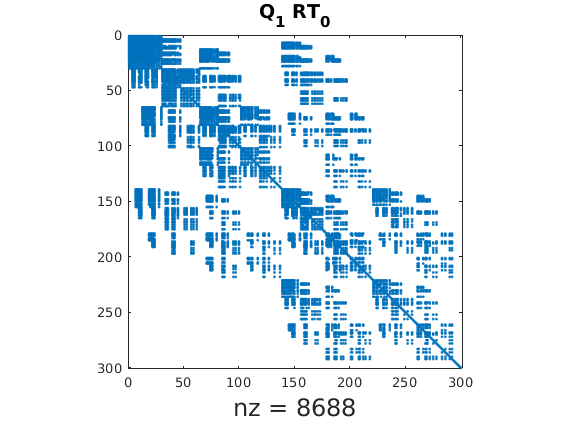}%
\includegraphics[width=0.5\textwidth]{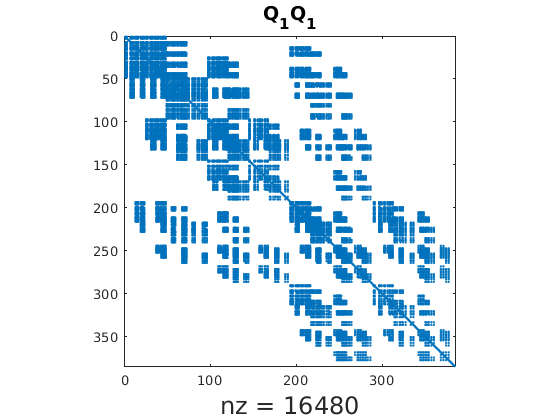}
\caption{Sparsity pattern for the {\FR tangent} matrix of the mechanically induced problem using the {\FR $Q_1 RT_0$} ansatz functions (top) and and the {\FR $Q_1 Q_1$} ansatz functions (bottom). Here, 27 finite elements are employed resulting in 300 degrees of freedom for {\FR $Q_1 RT_0$} and respectively 384  degrees of freedom for {\FR $Q_1 Q_1$}}
\label{fig:spy}
\end{figure*}
To explain this effect, in Figure~\ref{fig:spy} the sparsity patterns for $Q_1 RT_0$ and $Q_1Q_1$ %
are displayed. 
Although the difference is visually not pronounced, 
the number of nonzeros almost doubles using {\FR $Q_1Q_1$}.
Precisely, we have for our example with %
27 finite elements a tangent matrix size of %
300 with 8688 nonzero entries for $Q_1RT_0$, which compares to a a size of %
384 with 16480 nonzero entries for $Q_1 Q_1$.
Therefore, it is not surprising if the factorizations are more computationally expensive using {\FR $Q_1Q_1$}.

In Table~\ref{tab:strongMech}, good strong parallel scalability is generally observed when scaling from 64 to 216 cores. 
Again, only incremental improvements are visible, when using 512 cores. This is, again, an indication that the problem is too small for 512 cores. A three-level approach could also help, here.

If we relate the results to the linear elastic benchmark problem~in Section~\ref{sec:linElas}, we see that the number of Krylov iterations is now larger by a factor of 4 to 6, which is significant. This is not reflected in the {\it Solver Time}, since it is, again, dominated by the local factorization. However, for a local solver significantly faster than KLU, the large number of Krylov iterations could be very relevant for the time-to-solution. %

Next, we therefore investigate if we can reduce the number of iterations by improving the preconditioner.
}

\begin{table*}
\centering
\renewcommand{\arraystretch}{1.2}
\scriptsize{
\begin{tabular}{l || c c c c |c | c c c c| c}
\hline
\multicolumn{11}{c}{\bf Mechanically induced diffusion problem}\\ \hline
\multicolumn{1}{c||}{}&\multicolumn{10}{c}{\bf{FROSch with GDSW coarse space}} \\ \hline
  & \multicolumn{5}{c|}{{\FR $Q_1 RT_0$}} & \multicolumn{5}{c}{{\FR $Q_1 Q_1$}}\\
 & \multicolumn{5}{c|}{5 Newton steps/time step, 2 time steps}& \multicolumn{5}{c}{5 Newton steps/time step, 2 time steps} \\
\#cores & Avg.  & Assemble & Setup & Krylov &{\bf Solver}  & Avg.  & Assemble & {\FR Setup}& Krylov &{\bf Solver}\\
       &  Krylov            & Time    & Time & Time & {\bf Time}&Krylov & Time &Time &Time & {\bf Time}\\ \hline
 64 & 119.90&	\unit[69.09]{s}& \unit[786.75]{s}& \unit[136.05]{s}& {\bf \unit[922.79]{s}} &115.80 &\unit[115.12]{s}&	\unit[1\,970.32]{s}&	\unit[227.09]{s}& {\bf \unit[2197.40]{s}} \\
 125 &137.20&  \unit[37.65]{s}&	\unit[393.31]{s}& \unit[106.14]{s}&{\bf \unit[499.45]{s}} &130.00 &\unit[60.69]{s}&	\unit[1\,378.02]{s}&\unit[259.87]{s}& {\bf \unit[1637.99]{s}}\\
 216 &131.70&	\unit[20.61]{s}& \unit[161.99]{s}& \unit[51.66]{s}&{\bf \unit[213.65]{s}}  &125.40 &\unit[33.66]{s}&	\unit[364.36]{s}& \unit[78.22]{s} & {\bf \unit[442.58]{s}}\\
 512 &180.30&	\unit[9.94]{s}&	\unit[121.26]{s}& \unit[90.62]{s}& {\bf \unit[211.88]{s}} &177.10 & \unit[16.42]{s}&	\unit[218.26]{s}& \unit[97.03]{s} & {\bf \unit[315.29]{s}}\\ \hline
\multicolumn{1}{c||}{}&\multicolumn{10}{c}{\bf{FROSch with RGDSW coarse space}} \\ \hline
& \multicolumn{5}{c|}{{\FR $Q_1 RT_0$}} & \multicolumn{5}{c}{{\FR $Q_1 Q_1$}}\\
 & \multicolumn{5}{c|}{5 Newton steps/time step, 2 time steps}& \multicolumn{5}{c}{5 Newton steps/time step, 2 time steps} \\
\#cores & Avg.  & Assemble & Setup & Krylov &{\bf Solver}  & Avg.  & Assemble & {\FR Setup}& Krylov &{\bf Solver}\\
       &  Krylov            & Time    & Time & Time & {\bf Time}&Krylov & Time &Time &Time & {\bf Time}\\ \hline
 64 & 119.80&\unit[66.19]{s}&\unit[794.16]{s}&\unit[132.48]{s}& {\bf \unit[926.63]{s}}  &116.00& \unit[112.38]{s}&	\unit[1\,962.97]{s}&	\unit[224.79]{s}& {\bf \unit[2187.76]{s}}\\
 125 & 136.60&\unit[39.57]{s}&\unit[411.44]{s}&\unit[102.52]{s} & {\bf \unit[513.96]{s}} &131.60& \unit[60.21]{s}&	\unit[1\,635.75]{s}& \unit[335.15]{s}& {\bf \unit[1970.90]{s}}  \\
 216 & 132.70&\unit[21.80]{s}&\unit[165.22]{s}&\unit[50.10]{s}& {\bf \unit[215.32]{s}}&126.00&	\unit[34.71]{s}&	\unit[357.97]{s}& \unit[73.58]{s} & {\bf \unit[431.55]{s}}\\
 512 & 176.60&\unit[10.31]{s}&\unit[103.88]{s}&\unit[51.26]{s}& {\bf \unit[155.13]{s}}   &177.10& \unit[16.35]{s}&	\unit[205.30]{s}& \unit[75.52]{s}& {\bf \unit[280.82]{s}}\\ \hline
\end{tabular}}
\caption{\FR Strong scalabilty results for the \textit{mechanically induced diffusion problem} corresponding to Figure~\ref{fig:mech-strong-1node}. We operate on a triangulation with 110\,592 finite elements resulting in 691\,635 degrees of freedom for the $Q_1 RT_0$ ansatz functions and 705\,894 degrees of freedom for the $Q_1 Q_1$ ansatz functions. We apply the FROSch framework in fully algebraic mode  with the GDSW and the RGDSW coarse space. We perform two time steps with $\Delta t=\unit[0.05]{s}$. By {\it Avg. Krylov} we denote the average number of Krylov iteration computed over all time- and Newton steps. The time measurements are taken over the complete computation.}
\label{tab:strongMech}
\end{table*}

\begin{table*}
\centering
\renewcommand{\arraystretch}{1.2}
\scriptsize{
\begin{tabular}{l || c c c c|c c c c}
\hline         
&  \multicolumn{4}{c|}{{ $Q_1 RT_0$}} & \multicolumn{4}{c}{{$Q_1 Q_1$}} \\ \hline

\multicolumn{3}{c}{} &{\bf GDSW} & {\bf RGDSW} &\multicolumn{2}{c}{} &{\bf GDSW} & {\bf RGDSW} \\ \hline
\# Cores &  Avg. Size & Max. Size & Comp. $K_i$ &Comp. $K_i$ &  Avg. Size & Max. Size & Comp. $K_i$& Comp. $K_i$\\
         &  $K_i$     &  $K_i$    & Time & Time  &  $K_i$     &  $K_i$    & Time & Time \\ \hline
   64	 &    29316	  &	 35634	&\unit[667.68]{s}& \unit[670.62]{s}	 & 31404   &	38334	& \unit[1739.8]{s} &\unit[1732.00]{s}\\
   125   &  19461	& 25598	& \unit[310.90]{s} &	\unit[330.26]{s} &	21044 & 27648 &\unit[1220.0]{s} & \unit[1481.00]{s}\\
   216   & 11022	& 14028 & \unit[136.91]{s} & \unit[141.31]{s} & 11958 &15162 & \unit[316.70]{s} &\unit[311.32]{s} \\
   512	 &  8937.4	& 12571 & \unit[74.09]{s} & \unit[76.85]{s} &9849.4	&13896 &\unit[153.13]{s}&\unit[150.64]{s} \\ \hline
\end{tabular}}
\caption{Detailed cost of the overlapping subdomain problems $K_i$ contained in the \textit{Setup Time} in Table~\ref{tab:strongMech}}
\label{tab:mech-overlap}
\end{table*}
 
 \begin{figure}[h]
  \centering
  \includegraphics[width=0.35\textwidth]{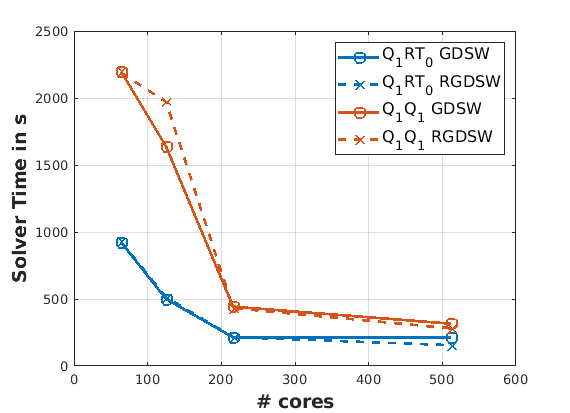}
  \caption{\FR Strong scalability of the \textit{Solver Time} (\textit{Setup Time + Krylov Time}) for the \textit{mechanically induced problem}; see Table~\ref{tab:strongMech} for the data.}
  \label{fig:mech-strong-1node}
\end{figure}

{\OR
\paragraph{Making more use of the problem structure in the preconditioner}
\label{sec:better}
We have observed in Table~\ref{tab:strongMech} that the number of Krylov iterations increased by roughly 50 percent when scaling from 64 to 512 processor cores. We explain this by the use of the fully algebraic mode of FROSch which applies the null space of the Laplace operator to construct the second level. 

To improve the preconditioner, the idea is to use, for the $Q_1Q_1$ discretization,
the three translations 
(in $x$, $y$, and $z$ direction) for the construction of the coarse problem for, both, the structure and the diffusion problem: we use six basis vectors $(1,0,0,0,0,0),\,(0,1,0,0,0,0),\ldots,(0,0,0,0,0,1)$ for construction of the coarse space. Here, the first three components refer to the structure problem and the last three components to the diffusion problem.
Note that the rotations are missing from the coarse space, as usual in this paper, since they would need access to the point coordinates.

In Table~\ref{tab:better}, we see
that, by using this enhanced coarse space, we can reduce the number of Krylov iterations, and, more importantly, we can avoid the increase in the number of iterations.
This means that, experimentally, in Table~\ref{tab:better}, we observe numerical scalability within the range of processor cores considered here.

Note that the larger coarse space resulting from the approach described above is not amortized in terms of computing times for the two-level preconditioner. Therefore, within the range of processor cores considered here, the fully algebraic approach is preferable. Again, a three-level method may help, especially, for the GDSW preconditioner, which has a larger coarse space compared to RGDSW.

A similar approach could be taken for $Q_1RT_0$, i.e., the three translations could be used for the deformation, and the Laplace nullspace could be used for the diffusion. However, this was not tested here.
}

\begin{table}[]
\centering
\renewcommand{\arraystretch}{1.2}
\scriptsize{
\begin{tabular}{l || c c }
\hline         
\multicolumn{3}{c}{\bf Mechanically induced diffusion problem}\\
&  \multicolumn{2}{c}{{$Q_1 Q_1$}} \\ \hline

\multicolumn{1}{c}{} &{\bf GDSW} & {\bf RGDSW} \\
\# Cores &  Avg.   & Avg. \\
         &  Krylov & Krylov\\ \hline
64  &     60.70 & 74.70\\
125 &     60.60 & 73.10\\
216 &     44.80 & 50.10\\
512 &     56.10 & 73.90\\
\end{tabular}
\caption{\OR A lower number of Krylov iterations and better numerical scalability is obtained when a better coarse space is used, spanned by the three translations in $x$, $y$, and $z$ direction, for, both, the structure as well as the diffusion problem.}
\label{tab:better}
}
\end{table}

\subsection{Weak parallel scalability}
\label{sec:weak-scaling}
{\OR For %
the weak parallel scalability, %
we only consider}
the {\FR \textit{free-swelling problem} with type~(ii) Neumann boundary condition as described in Section~\ref{sec:free-swell-bvp}.
The material parameters are chosen as in Section~\ref{subsec:free-strong}.}
{\FR On the initial mesh we perform different numbers of refinement cycles  such that we obtain 512 finite elements per core. }

{\FR For the smallest problem of 8 cores, we have a problem size of 27\,795 which compares to largest problem size of 1\,622\,595 using 512 cores.}

{\FR Hence, %
we increase the problem size {\OR as well as the number of processor cores.} %
{\OR For this setup, in the best case,} the average number of Krylov iterations {\it Avg. Krylov} and the {\it Solver Time} should remain constant.}

{\OR We observe that, within the range of $8$ to $512$ processor cores considered, the number of Krylov iterations grows from $16.1$ to $71.4$ for GDSW and from $14.6$ to $76.6$ for RGDSW. This increase is also reflected in the {\it Krylov Time}, which increases from $\unit[2.51]{s}$ to $\unit[33.52]{s}$ for GDSW and from $\unit[2.06]{s}$ to $\unit[31.03]{s}$ for RGDSW.

However, %
a significant part of the increase in the {\it Solver Time} comes from a load imbalance in the problems with more than 64 cores: the maximum local subdomain size is
$7\,919$ for 8 cores and $14\,028$ for 512 cores; see also Table~\ref{tab:overlapWeak}.
However, Table~\ref{tab:overlapWeak} also shows that the load imbalance does not increase when scaling from 64 to 512 cores, which indicates that the partitioning scheme works well enough.

Again, we see that the coarse level is not fully sufficient to obtain numerical scalability, and more structure should be used, as in Table~\ref{tab:better}, if numerical scalability is the goal.

Note that even for three-dimensional linear elasticity, we typically see an increase of the number of Krylov iterations when scaling from $8$ to $512$ cores, even for GDSW using with the full coarse space for elasticity (including rotations)~\cite[Fig. 15]{Heinlein:2016:PIT}, %
which is covered by the theory. 
In \cite[Fig. 15]{Heinlein:2016:PIT} the number of Krylov iterations only stays almost constant beyond about $2000$ cores. %
However, the increase in the number of iterations is mild for the full GDSW coarse space in \cite[Fig. 15]{Heinlein:2016:PIT} when scaling from $64$ to $512$ cores.

Concluding, we see that using the fully algebraic mode, 
for the range of processor cores considered here, leads to an acceptable method although numerical scalability and optimal parallel scalability is not achieved.

Interestingly, the results for RGDSW are very similar in terms of Krylov iterations as well as the {\it Solver Time} although RGDSW has a significantly smaller coarse space. This advantage is not yet visible here, however, for a larger number of cores RGDSW can be expected to outperform GDSW by a large margin.

\begin{table*}
\centering
\renewcommand{\arraystretch}{1.2}
\scriptsize{
\begin{tabular}{l||c c c c |c|c c c c| c }
\hline &\multicolumn{5}{c|}{\bf FROSch with GDSW coarse space} &\multicolumn{5}{c}{\bf FROSch with RGDSW coarse space} \\ \hline
 & \multicolumn{5}{c}{\OR $Q_1 RT_0$} & \multicolumn{5}{c}{\OR $Q_1 RT_0$} \\
\#Cores & Avg. &Assemble  & {\FR Setup}& Krylov&{\bf Solver} & Avg.    &Assemble  & {\FR Setup}& Krylov&{\bf Solver} \\
        &Krylov&Time      &Time   & Time  &{\bf Time}  &  Krylov & Time            & Time   & Time&{\bf Time} \\ \hline
8  &16.10  &	\unit[17.39]{s}& \unit[46.90]{s}& \unit[2.51]{s} &{\bf \unit[49.41]{s}} &	14.60&	\unit[16.66]{s}& \unit[40.51]{s}& \unit[2.06]{s}&{\bf \unit[42.57]{s}}\\
27 &24.90  &	\unit[19.51]{s}& \unit[81.07]{s}& \unit[6.11]{s}&{\bf \unit[87.18]{s}} & 24.70&	\unit[20.14]{s}&	\unit[83.75]{s}& \unit[5.99]{s}&{\bf \unit[89.74]{s}}\\
64 &34.30  & \unit[20.39]{s}& \unit[147.33]{s}&	\unit[11.92]{s}&{\bf \unit[159.25]{s}}&	34.90&	\unit[20.28]{s}&	\unit[151.13]{s}& \unit[12.14]{s}&{\bf \unit[163.27]{s}} \\
216&52.20 & \unit[21.82]{s}& \unit[162.51]{s}& \unit[20.83]{s}&{\bf \unit[183.34]{s}} &	54.50&	\unit[20.33]{s}&	\unit[160.51]{s}& \unit[19.97]{s}&{\bf \unit[180.58]{s}}\\
512&71.40 & \unit[21.45]{s}&\unit[186.74]{s}&	\unit[33.62]{s}&{\bf \unit[220.36]{s}}&	76.60&	\unit[21.19]{s}&	\unit[183.04]{s}& \unit[31.03]{s}&{\bf \unit[213.07]{s}}\\\hline
\end{tabular}}
\caption{\FR Weak parallel scalability results for the \textit{free swelling problem} with  type (ii) boundary conditions described in Section~\ref{sec:free-swell-bvp} corresponding to Figure~\ref{fig:weakRT}. Each core owns approximately  512 finite elements. We use the $Q_1 RT_0$ ansatz functions and apply FROSch with the GDSW and the RGDSW coarse space. Two time step with $\Delta t = \unit[0.05]{s}$ are performed with each requiering 5 Newton steps. By {\it Avg. Krylov}, we denote the average number of Krylov iteration over all time- and Newton steps. The time measurements are taken over the whole computation.}
\label{tab:weakRT}
\end{table*}

\begin{figure}[]
\includegraphics[width=0.4\textwidth]{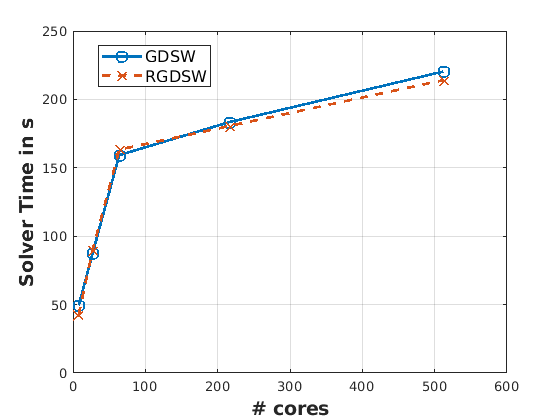}
\caption{\FR Weak parallel scalability of the \textit{Solver Time (Setup Time + Krylov Time)} for the model problem of a \textit{free swelling cube} with  type (ii) boundary conditions described in Section~\ref{sec:free-swell-bvp}. See Table~\ref{tab:weakRT} for the data.}
\label{fig:weakRT}
\end{figure}

Increasing the %
{\FR the number of finite elements assigned to one core}, we expect the weak parallel scalability to improve.}

\begin{table*}
\centering
\renewcommand{\arraystretch}{1.2}
\scriptsize{
\begin{tabular}{l || c c c |c  c }
\hline
         &               &           &          & {\bf GDSW}         & {\bf RGDSW}       \\\hline
         &               &           &          & \multicolumn{2}{c}{\OR $Q_1 RT_0$} \\
\# Cores &  Total Problem& Avg. Size & Max. Size & Comp. $K_i$ & Comp. $K_i$ \\
         &  Size         &  $K_i$    &  $K_i$    & Time        & Time \\ \hline    
  8      & 27\,795       & 7\,504.4  & 7\,919    & \unit[37.13]{s}       & \unit[31.83]{s}\\
  27     & 90\,075       & 9\,139.1  & 12\,207	 & \unit[63.67]{s} 	   & \unit[65.71]{s}\\
  64	 & 209\,187      & 10\,048   & 14\,028   & \unit[124.68]{s}      & \unit[127.66]{s}\\
  216	 & 691\,635      & 11\,022   & 14\,028	 & \unit[137.77]{s}      & \unit[137.05]{s}\\
  512    &1\,622\,595    & 11\,533   &14\,028	 & \unit[158.51]{s}	   & \unit[158.81]{s}\\
  \hline
\end{tabular}}
\caption{\FR Detailed cost of the overlapping subdomain problems $K_i$ for the weak parallel scalabilty results in Table~\ref{tab:weakRT}.}
\label{tab:overlapWeak}
\end{table*}

\subsection{Conclusion and Outlook}
{\FR The FROSch framework has shown a good parallel performance applied algebraically to the fully coupled chemo-mechanic problems. We have compared two benchmark problems {\SP with} different boundary conditions.
The time step size was of minor influence for our specific benchmark problems. }

{\OR
Our GDSW-type preconditioners implemented in FROSch are a suitable choice when used as a monolithic solver for the linear systems arising from the Newton linearization. They perform well when applied in fully algebraic mode even when numerical scalability is not achieved.

Our experiments of strong scalability have shown that, with respect to average time to solve a monolithic system of equations obtained from linearization of the nonlinear monolithic coupled problem discretized by $Q_1 RT_0$, we have to invest a factor of slightly more than three in computing time compared to a solving a standard linear elasticity benchmark problem discretized using the same number of $Q_1$ finite elements. This is a good result, considering that the monolithic system is larger by a factor of almost 2.

Using a $Q_1 Q_1$ discretization, the computing times are much slower. This is mainly a result of a lower sparsity of the finite element matrices.

We have also discussed that, using more structure, we can achieve numerical scalability in our experiments. However, this approach will only be efficient when used with a future three-level extension of our preconditioner.
}

\section{Acknowledgement}
The authors acknowledge the DFG project 441509557
(\url{https://gepris.dfg.de/gepris/projekt/441509557})
within the Priority Program SPP2256 “Variational Methods for Predicting Complex Phenomena in Engineering Structures
and Materials” of the Deutsche Forschungsgemeinschaft (DFG). 

The authors also acknowledge the compute cluster ({DFG project no. 397252409},
\url{https://gepris.dfg.de/gepris/projekt/397252409}) 
of the Faculty of Mathematics
and Computer Science of Technische Universit\"at Bergakademie Freiberg, operated by the Universit\"atsrechenzentrum URZ.

\begingroup
    \setlength{\bibsep}{2pt}
    \setstretch{1}
\bibliography{lit.bib}
\endgroup
\end{document}